\documentclass[10pt,a4paper]{article}
\usepackage[utf8]{inputenc}
\usepackage[T1]{fontenc}
\usepackage{amsmath}
\usepackage{amsfonts}
\usepackage{array}
\usepackage{amssymb}
\usepackage{graphicx}

\usepackage{subfigure}

\usepackage[usenames,dvipsnames]{xcolor} 
\usepackage{cite}
\usepackage[cm]{fullpage}
\usepackage{layout}

\usepackage{mathtools}
\usepackage{amstext}
\usepackage{amsthm}
\usepackage{esint}

\usepackage{babel}

\makeatletter
\theoremstyle{plain}

\theoremstyle{remark}

\theoremstyle{plain}

\makeatother

\usepackage{babel}
\providecommand{\conjecturename}{Conjecture}
\providecommand{\remarkname}{Remark}
\providecommand{\theoremname}{Theorem}

\providecommand{\keywords}[1]{{\textit{Keywords}} #1}

\makeatletter
\newcommand{\thickhline}{%
	\noalign {\ifnum 0=`}\fi \hrule height 1pt
	\futurelet \reserved@a \@xhline	
}
\newcolumntype{"}{@{\hskip\tabcolsep\vrule width 1pt\hskip\tabcolsep}}
\makeatother
\newcommand{\uvec}{\mathbf{u}}
\newcommand{\lap}{\bigtriangleup}

\newcommand{\normal}{\mathbf{n}}

\newcommand{\Domain}{\mathbf{D}}
\usepackage{caption}
\begin{document}
	
	\title{Solving incompressible Navier--Stokes equations on irregular domains and quadtrees by monolithic approach}
	\author{Hyuntae Cho, Yesom Park, and Myungjoo Kang}
	\maketitle
	\begin{abstract}
		We present a second-order monolithic method for solving incompressible Navier--Stokes equations on irregular domains with quadtree grids. A semi-collocated grid layout is adopted, where velocity variables are located at cell vertices, and pressure variables are located at cell centers. Compact finite difference methods with ghost values are used to discretize the advection and diffusion terms of the velocity. A pressure gradient and divergence operator on the quadtree that use compact stencils are developed. Furthermore, the proposed method is extended to cubical domains with octree grids. Numerical results demonstrate that the method is second-order convergent in $L^\infty$ norms and can handle irregular domains for various Reynolds numbers.
		\\
		\\
		\keywords{Incompressible Navier--Stokes, Finite difference method, Quadtree, Octree, Level set method}
	\end{abstract}

	\section{Introduction}
	The motion of an incompressible fluid flow is described by the Navier--Stokes equations:
	\begin{equation}\label{eq:NS}
		\begin{gathered}
			\uvec_t + \uvec \cdot \nabla \uvec =\frac{1}{Re} \bigtriangleup \uvec- \nabla p +\mathbf{f},\\
			\nabla \cdot \uvec =\mathbf{0}.
		\end{gathered}
	\end{equation}
	Here, $\uvec$ is the fluid velocity, $p$ is the pressure, $\mathbf{f}$ is the external force such as gravity, and $Re$ is the Reynolds number. The development of numerical methods for incompressible Navier--Stokes equations is very important in both science and engineering. Applications of incompressible flow, such as subsonic aerodynamics, ship hydrodynamics, and numerical weather prediction, sometimes involve complex geometry and require that small-scale turbulence and vortices be captured. An efficient way to solve the problem is to generate a body-fitted mesh and adaptively refine the mesh where small-scale phenomena are observed. Although various successful studies in recent decades have used body-fitted meshes, we focus on Cartesian grids in this study.
	
	Peskin's immersed boundary method (IBM) \cite{peskin1977numerical} is a popular and effective approach to handling irregular domains on Cartesian grids. However, the numerical delta function of the IBM lowers the order of convergence near the boundary. To address this issue, Leveque and Li \cite{leveque1994immersed} developed the immersed interface method (IIM) for elliptic interface problems, which shows second-order convergence in the $L^\infty$ norm. The IBM also motivated ghost fluid methods (GFMs) \cite{fedkiw1999non,liu2000boundary,kang2000boundary}, which define ghost values at grid points outside the domain. Although the IIM and GFM were developed to handle discontinuity across the interface, the core concept of these sharp capturing methods has been extended to treat various boundary conditions for elliptic equations on irregular domains \cite{gibou2002second}. In addition to the IIM and GFM, various numerical methods are used to treat irregular domains. For example, in the virtual node method \cite{bedrossian2010second,hellrung2012second}, variational formulations on irregular domains embedded in a Cartesian grid are discretized. Coco and Russo developed ghost point finite difference methods \cite{coco2013finite,coco2018second,coco2020multigrid}, which adopt a multigrid method to impose boundary conditions. In a finite volume framework, cut-cell methods are widely used. Originally proposed by Purvis and Bulkhalter \cite{purvis1979prediction}, these methods were used to develop Hodge projections of vector fields \cite{ng2009efficient}. Recent works on finite volume cut-cell methods include the solution of Poisson's equation \cite{bochkov2019solving} and the convection--diffusion equation \cite{barrett2021hybrid} with Robin boundary conditions.
	
	Chorin introduced the projection method \cite{chorin1997numerical} to solve incompressible Navier--Stokes equations on Cartesian grids. It uses Hodge decomposition of vector fields. The intermediate velocities are first computed and then projected to force a discrete divergence-free condition by solving Poisson's equation for the Hodge variable. Although the method developed by Chorin is only first-order accurate, Brown et al. \cite{brown2001accurate} analyzed a second-order projection method. They showed that the boundary conditions of the intermediate velocities are coupled with the Hodge variable, and thus appropriate boundary conditions for the intermediate velocities are necessary to obtain second-order-accurate solutions. 
	For rectangular domains, the boundary conditions for the intermediate velocities and the Hodge variable can be decoupled by using the Hodge variable from previous time steps. However, for irregular domains, it is challenging to impose an accurate boundary condition. Furthermore, especially when the domain moves in time according to the flow, the intermediate velocities and Hodge variable should be solved simultaneously in a single linear system. Thus, the advantage of using the projection method is lost. 
	
	An alternative approach to solving the incompressible Navier--Stokes equations is to discretize the momentum equation while imposing a divergence-free condition in a single linear system for the velocities and pressure. This discretization is called the \it monolithic \rm approach and results in a saddle-point linear system. Although this method creates a larger linear system than projection methods do, the boundary conditions for the velocities are no longer coupled with other variables. Similarly, the jump conditions of two-phase incompressible flow can be treated more easily than they can by projection methods, and recent works \cite{cho2020fully,chen2018direct,schroeder2014second} have reported success in obtaining convergence in $L^\infty$ norms. For a single-phase flow problem on irregular domains, a very recent study by Coco \cite{coco2020multigrid} solved the equations using the monolithic approach. The author used a ghost-point finite difference method to impose boundary conditions and obtained second-order convergence of the velocity and divergence of the velocity in the $L^1$ and $L^\infty$ norms. These monolithic methods exhibit convergence with nontrivial analytical solutions; note, however, that discretized linear systems are rank-deficient (the pressure is unique up to a constant). In \cite{cho2020fully,schroeder2014second}, the right-hand side was projected onto the range space of the linear system by formulating the kernel of its transpose. By contrast, Coco \cite{coco2020multigrid} augmented the approach with an additional scalar unknown. Although Yoon et al. \cite{yoon2016solving} remarked that
	finding a unique solution in the orthogonal complement of the kernel is equivalent to the augmentation strategy, the range space of the resulting linear system was not specified by Coco \cite{coco2020multigrid}. 
	
	Staggered marker-and-cell (MAC) grids \cite{harlow1965numerical} have been widely used to solve the incompressible Navier--Stokes equations \eqref{eq:NS} numerically. In a MAC grid, the velocities are located on cell faces, whereas the pressures or Hodge variables are located at cell centers. MAC grids do not exhibit the checkerboard instability that occurs in a collocated grid, and they enable a simple discretization for the incompressibility condition. However, there are obstacles to developing numerical formulations on MAC grids for quad/octree data structures.
	The extension of numerical methods designed for uniform grids to the quad/octree grid is not trivial, and it requires additional effort. 
	Early works on solving incompressible flow on quad/octree grids focused on the Euler equations, where the viscosity term in \eqref{eq:NS} was neglected. Popinet \cite{popinet2003gerris} introduced a second-order projection method on a \it graded \rm octree grid, that is, a grid where the size ratio of adjacent cells must be 2:1. This work was later extended by Losasso et al. \cite{losasso2004simulating} to relax the restrictions on adaptive refinement. Although Popinet \cite{popinet2003gerris} discretized the advection term in the Eulerian approach using an upwind method, Losasso et al. \cite{losasso2004simulating} used a semi-Lagrangian method by extrapolating the velocities located on the cell faces to the cell vertices. For cases where the viscosity terms are not neglected, Olshanskii et al. \cite{olshanskii2013octree} and Guittet et al. \cite{guittet2015stable} developed second-order projection methods on a graded and non-graded octree MAC grid, respectively.
	The method of Olshanskii et al. \cite{olshanskii2013octree} requires a relatively large stencil for discretizing the viscous term. By contrast, Guittet et al. \cite{guittet2015stable} developed a finite volume method that requires a smaller stencil; however, a Voronoi mesh was generated with respect to the cell faces.

	Min and Gibou \cite{min2006second} developed a second-order projection method for incompressible flow on quad/octree grids other than MAC grids, where the velocities and Hodge variables are collocated at cell vertices. The main advantage of the collocated grid on vertices is that second-order finite differences with compact stencils can be used. As in \cite{guittet2015stable}, Hodge decomposition that guarantees stability was proposed. However, the grid layout in \cite{min2006second} cannot be used in the monolithic approach because there is no boundary condition for the pressure, whereas homogeneous Neumann boundary conditions are imposed for the Hodge variable in the projection method.

	In this paper, we introduce a numerical method for solving the incompressible Navier--Stokes equations on an irregular domain with non-graded quadtree structure and on a cubical domain with non-graded octree structure. Unlike many other methods on quad/octree grids, the proposed discretization is based on the monolithic approach. The velocities are collocated at cell vertices, and the pressure variable is stored at cell centers, which is also called the Arakawa B grid layout \cite{arakawa1977computational}. Because the Arakawa B grid layout is adopted, the effective and simple finite differences in \cite{min2006second} on the quadtree can be used. Furthermore, we propose gradient and divergence operators that map from cell centers to cell vertices and from cell vertices to cell centers, respectively. These operators are extended to treat an irregular boundary, and numerical results demonstrate that the proposed discretization is second-order accurate. The rest of the paper is structured as follows. Section 2 describes a level set method of capturing irregular domains and the proposed grid layout on a quadtree data structure. Section 3 explains the numerical method and discretization in detail. Numerical experiments and results are presented in Section 4.

	\section{\label{sec:GridGeneration} Computational domain and grid layout}

	\subsection{Level set method}
	We are interested in solving the incompressible Navier--Stokes equations \eqref{eq:NS} in an irregular domain $\Omega$ within a rectangular computational domain $\Domain=\left[a,b\right]\times\left[ c,d\right]$. The interface between the fluid and the solid is captured by the level-set method \cite{osher1988fronts}. The interface is described by the zero level set of a level set function $\phi: \Domain \to \mathbb{R}$ as
	$
	\Gamma= \left\{ \mathbf{x}\in\Domain | \phi(\mathbf{x},t)=0\right\}.
	$
	In particular, we define the fluid domain as
	$
	\Omega= \left\{ \mathbf{x}\in\Domain | \phi(\mathbf{x},t)\leq 0\right\},
	$
	and thus the solid domain is $\Domain \setminus \Omega$.
	
	Among infinitely many level sets that represent the interface as a zero level set, a signed distance function has the desirable property that $|\nabla \phi|=1$, and thus $\phi$ does not exhibit the instability caused by a steep or shallow gradient. Therefore, the level set function is reinitialized to become a signed distance function by solving the following
	eikonal equation: 
	\[
	\phi_{\tau}+\text{sgn}\left(\phi^{0}\right)\left(\left\Vert \nabla\phi\right\Vert -1\right)=0,
	\]
	where $\text{sgn}$ denotes the signum function, and $\phi^{0}$ is the
	initial level set function. In particular, we adopt the subcell resolution technique of Russo and Smereka \cite{russo2000remark} to discretize the equation and the second-order discretization of Min and Gibou \cite{min2007second}. An advantage of the level set method is that $\phi$ can be used to compute geometric quantities such as the normal vector and curvature according to the following formula:
	\[
	\mathbf{n}=\frac{\nabla\phi}{\left|\nabla\phi\right|}, \quad \kappa= \nabla \cdot \frac{\nabla\phi}{\left|\nabla\phi\right|}.
	\]
	
	\subsection{Quadtree data structure}
\begin{figure}	
		\begin{centering}
			\includegraphics[width=0.4\textwidth]{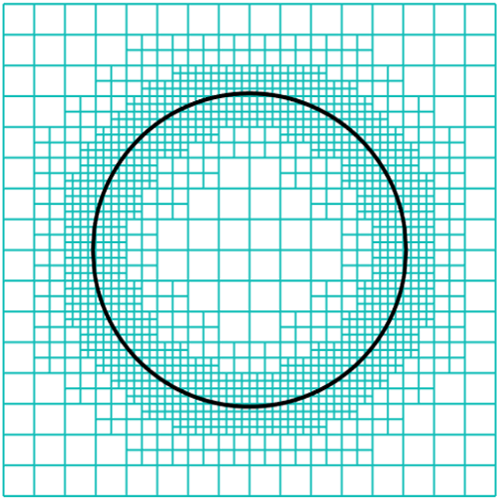}
			\caption{ An example of the quadtree mesh over a circular domain. The grid has 
				a maximum resolution of min level$=4$ and max level$=7$.\label{fig:quadtree}}
			\par\end{centering}
	\end{figure}

	For adaptive refinement, we define a root cell as the computational domain $\Domain$ and recursively split each cell
	into quadrants until the desired level of refinement is achieved;
	that is, a parent cell is divided into four children cells of equal size. Here, the level of the cell is defined as the number of refinements used for its construction. The notation $n/m$ describes the level a quadtree grid, where the minimum level is $n$, and the maximum level is $m$. Note that a quadtree grid of level $n/n$ is equivalent to a $2^n \times 2^n$ uniform Cartesian grid.
	In this study, we consider non-graded quadtree grids, and thus the size ratio of adjacent cells is not restricted.
	When a partial differential equation on an irregular domain is solved on the Cartesian grid, the order of the local truncation error near the boundary is usually lower than that of interior grid points. Therefore, as suggested in many works \cite{min2007second,theillard2019sharp,guittet2015stable}, we split a cell when its size exceeds
	its distance to the boundary. More precisely, a cell $C$ is refined
	if the following condition is satisfied:
	\[
	\underset{v\in\text{vertices}(C)}{\text{min}}\mid\phi(v)\mid\leq L\sqrt{2}\Delta x,
	\]
	where $L$ refers to the Lipschitz constant of the level set function
	$\phi$, and $\Delta x$ denotes the size of cell $C$ so that $\sqrt{2}\Delta x$
	indicates the size of the cell diagonal. As we reinitialize the level
	set function that allows $\phi$ to become a signed distance function,
	we take the Lipschitz constant to be $L=1.2$. Figure \ref{fig:quadtree} shows an example of a quadtree grid generated in a circular domain. This locally refined
	grid structure enables a storage reduction compared to the uniform
	Cartesian grid, whose size is the same as that of the finest cell of the quadtree.

	\subsection{Grid layout}
\begin{figure}		
		\begin{centering}
			\mbox{
				\subfigure[]{\includegraphics[width=0.4\textwidth]{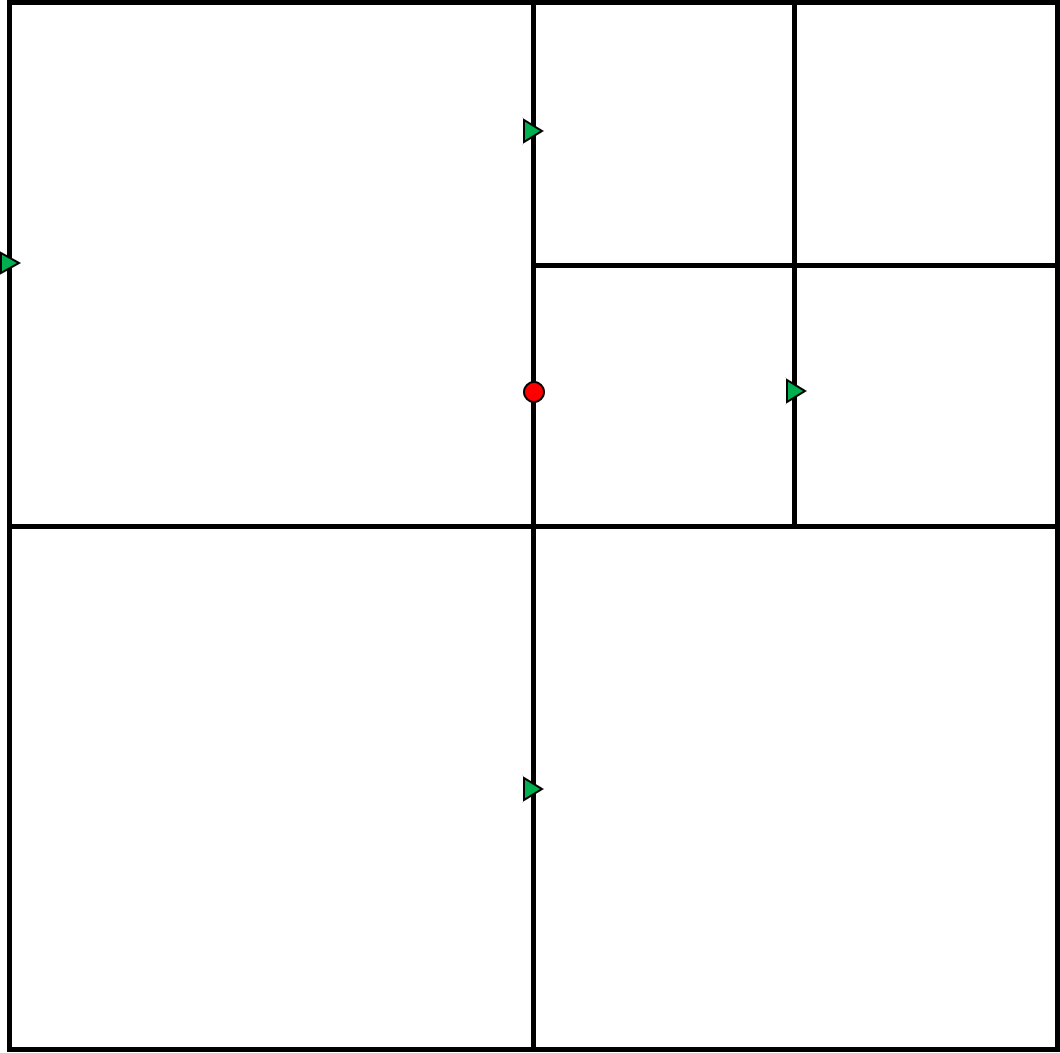}}
				$\ \ \ \ \ \ \ \ \ \ $ 
				\subfigure[]{\includegraphics[width=0.4\textwidth]{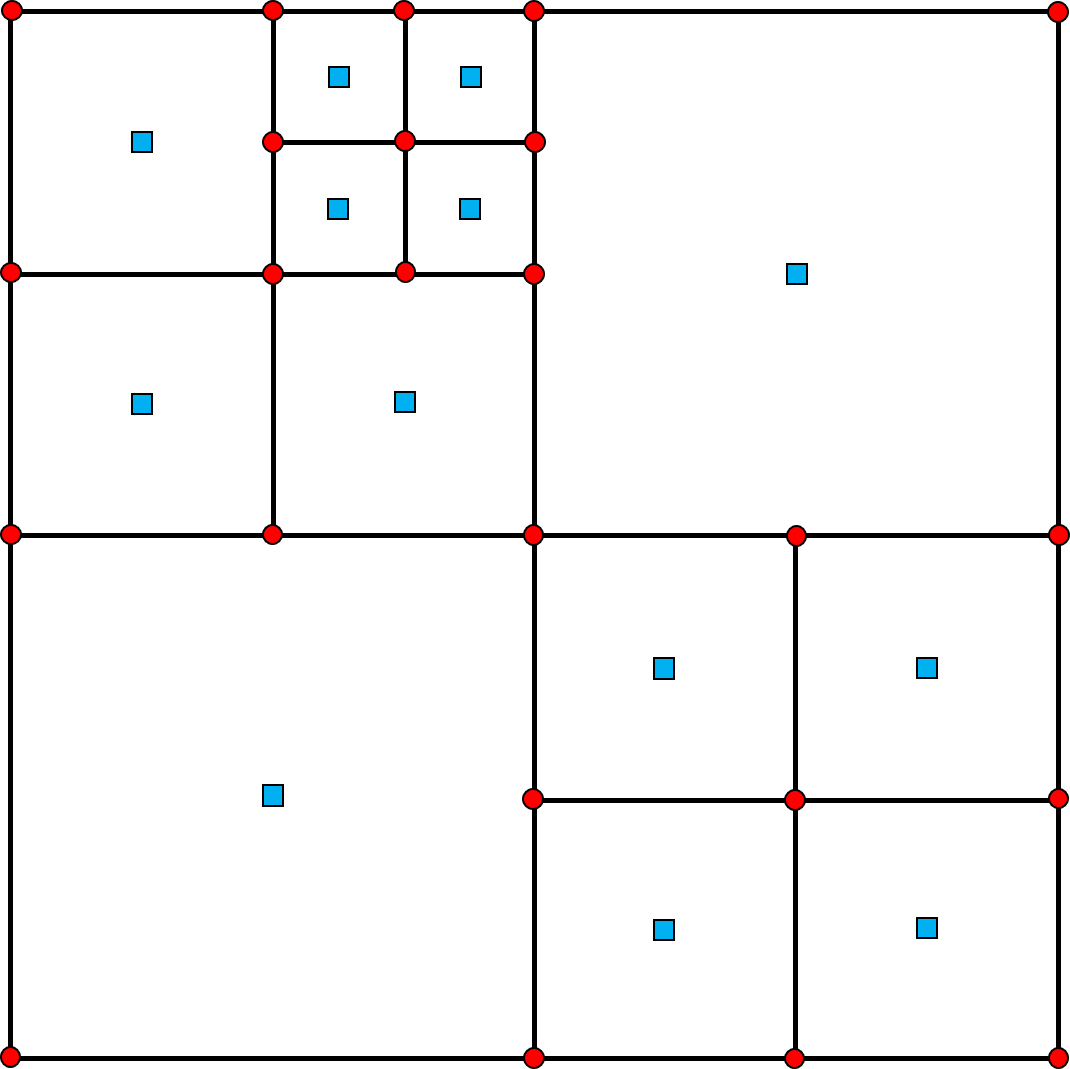}}
			}
			\caption{(a) Adjacent compact stencil in MAC grid.
				(b) General layout of the Arakawa B-grid on the quadtree mesh. Velocity variables are defined at the grid corners(red circles) 
				whereas the pressure variables are stored at the center of cells(blue squares). \label{fig:macquadtree}}
			\par\end{centering}
	\end{figure}

	A staggered MAC grid has typically been used \cite{guittet2015stable,olshanskii2013octree} to solve the incompressible Navier--Stokes equations \eqref{eq:NS} on quad/octree data structures. In a two-dimensional domain, the horizontal and vertical velocities are located on the horizontal and vertical cell faces, respectively, and the pressure variables are located at cell centers. One major difficulty of using a MAC grid on a quad/octree data structure is the discretization of the advection and viscous terms. For example, consider the discretization at the cell face marked with a circle in Figure \ref{fig:macquadtree}(a). The cell face marked with a triangle represents the $x$ component velocity $u$, which intersects the target cell face or is a face of the same cell. Within this compact stencil, $u_x$ and $u_{xx}$ cannot be discretized with errors of $o(\Delta x^2)$ and $o(\Delta x)$, respectively. To address this problem, Olshanskii et al. \cite{olshanskii2013octree} used wider stencils to discretize the advection and diffusion terms. By contrast, Guittet et al. \cite{guittet2015stable} addressed these difficulties by adopting a semi-Lagrangian scheme for the advection term and a finite volume approach to treat the viscous term. However, the quantities defined on the faces were interpolated to nodal values to employ the semi-Lagrangian scheme, and a Voronoi mesh with respect to the cell faces was constructed for the diffusion term. Furthermore, to ensure the stability, Olshanskii et al. developed an explicit filter that modifies the advection term near the coarse-to-fine grid interface, and Guittet et al. derived a discrete divergence operator that uses faces from different cells.
	
	Recognizing these difficulties, we use an Arakawa B grid \cite{arakawa1977computational} on the quad/octree data structure for spatial discretization, as shown in Figure
	\ref{fig:macquadtree}(b). That is, the velocity components are located at cell vertices, and the pressure variables are located at cell centers. One advantage of this grid layout is that the advection and diffusion terms can be discretized using compact stencils with a consistent truncation error. The Arakawa B grid layout is equivalent to the $Q1/P0$ finite element method and is known to exhibit checkerboard instability. By adopting a stabilizing technique used in the finite element community, we develop a second-order monolithic method with reduced checkerboard instability for solving \eqref{eq:NS} in non-graded quad/octree grids in the following section.

	\section{\label{sec:NumericalMethod}Numerical method}
	
	In this section, we present a finite-difference-based second-order-accurate monolithic scheme for solving the incompressible Navier--Stokes equations (\ref{eq:NS})
	on the non-graded quadtree structure described in Section \ref{sec:GridGeneration}.
	Time is discretized using the semi-Lagrangian method combined with the second-order backward difference formula introduced by Xiu et al. \cite{xiu2001semi}.
	By implicitly discretizing the incompressibility condition, we obtain the following temporal discretization for (\ref{eq:NS}):
	\begin{equation}\label{eq:saddle_system}
		\begin{cases}
			\frac{3\mathbf{u}^{n+1}-4\mathbf{u}_{d}^{n}+\mathbf{u}_{d}^{n-1}}{2\triangle t}&=-\nabla p^{n+1}+\frac{1}{Re}\triangle\mathbf{u}^{n+1}+f^{n+1},\\
			\nabla\cdot \mathbf{u}^{n+1}&=S(p^{n+1}).
		\end{cases}
	\end{equation}
	Here, the subscript $d$ indicates departure points in the semi-Lagrangian scheme, and $S$ is a stabilizing operator used to ease the checkerboard instability of the pressure.
	In the next section, we describe the spatial discretization and then discuss the induced linear system briefly.
	We point out that the
	suggested discretization uses only the nodal values of adjacent cells.
	
	\subsection{Advection term}
	
	The advection term in (\ref{eq:NS}) is treated by a semi-Lagrangian
	method, which discretizes the Lagrangian derivative of the velocity
	field in time instead of the Eulerian derivative. The semi-Lagrangian
	scheme consists of tracking a fluid particle backward by time integration
	of a characteristic equation and recovering the values at the departure
	point using suitable high-order interpolation schemes. 
	The
	departure points are traced by a second-order Runge--Kutta method as
	\begin{align*}
		\mathbf{x}_{d}^{n} & =\mathbf{x}_{ij}^{n+1}-\triangle t\mathbf{u}^{n+\frac{1}{2}}\left(\mathbf{x}_{ij}^{n+1}-\frac{\triangle t}{2}\mathbf{u}^{n}\left(\mathbf{x}_{ij}^{n+1}\right)\right),
	\end{align*}
	and
	\begin{align*}
		\mathbf{x}_{d}^{n-1} & =\mathbf{x}_{ij}^{n+1}-2\triangle t\mathbf{u}^{n}\left(\mathbf{x}_{ij}^{n+1}-\triangle t\mathbf{u}^{n}\left(\mathbf{x}_{ij}^{n+1}\right)\right).
	\end{align*}
	The velocity $\mathbf{u}^{n+\frac{1}{2}}$ at the intermediate time
	step is approximated by the second-order extrapolation 
	\[
	\mathbf{u}^{n+\frac{1}{2}}\left(\mathbf{x}\right)=\frac{3}{2}\mathbf{u}^{n}\left(\mathbf{x}\right)-\frac{1}{2}\mathbf{u}^{n-1}\left(\mathbf{x}\right).
	\]
	Because the points $\mathbf{x}_{ij}^{n+1}$ and $\mathbf{x}_{d}^{n}$
	usually do not coincide with grid nodes,
	appropriate spatial interpolations with nodal values are required to obtain the velocity values at these points. As
	the domain evolves with time, these points may lie outside the
	domain at preceding time steps, and neighboring nodes may
	be external points, making
	interpolation impossible. Here we design an interpolation scheme
	to compensate for these problems using boundary ghost points. To describe
	the interpolation procedure formally, we use the following notation.
	Suppose we would like to approximate $\mathbf{u}^{n}$ at a point $\mathbf{x}$
	and denote by $\mathbf{C}_{0}$ the leaf cell in $\Omega^{n}$ containing
	$\mathbf{x}$. 
	
	\begin{figure}[h]
		
		\centering{}
		\subfigure[]{\includegraphics[width=0.45\textwidth]{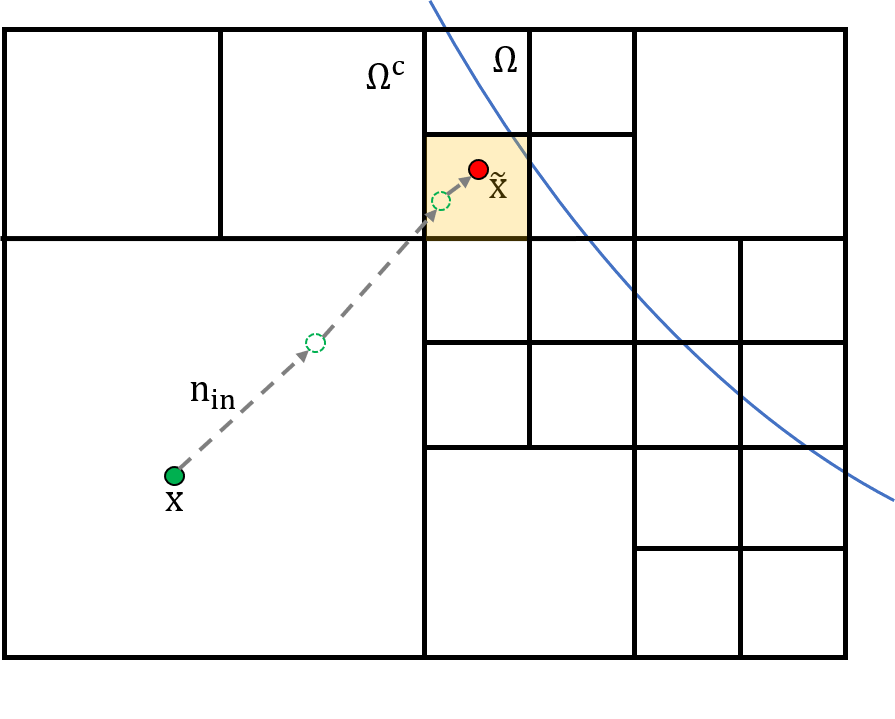}}
		\subfigure[]{\includegraphics[width=0.45\textwidth]{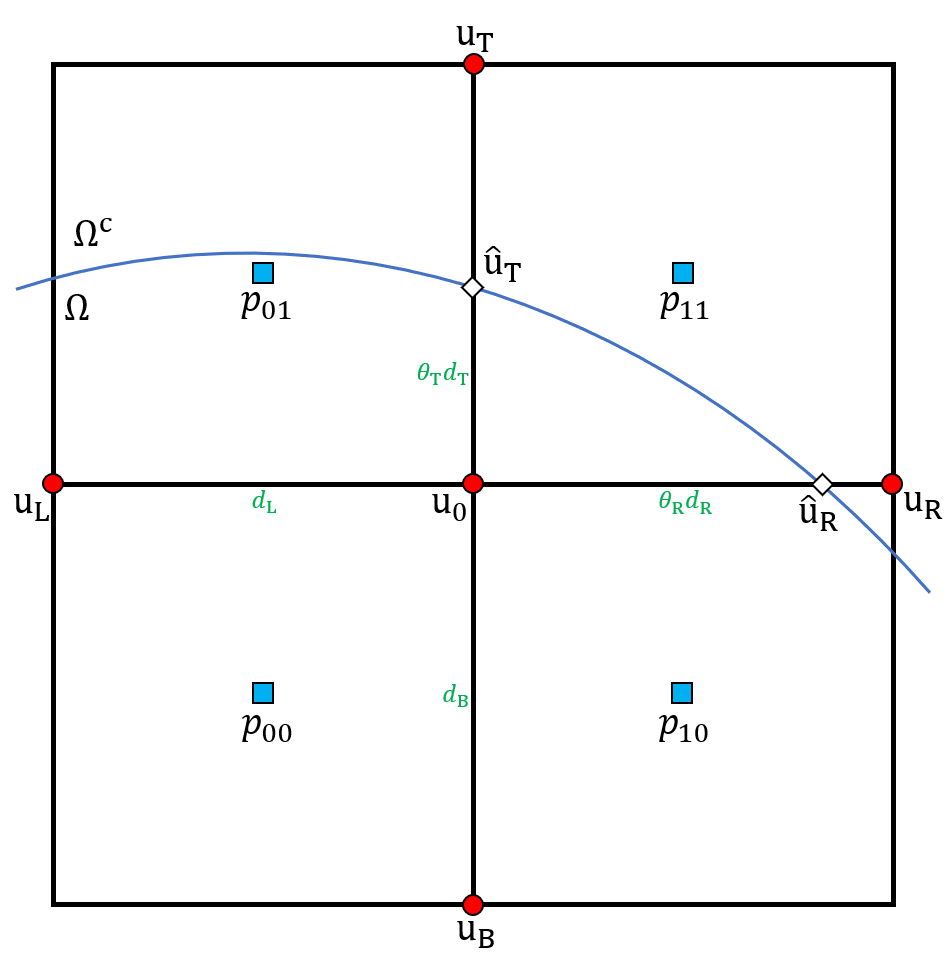}}
		\caption{Grid illustrations for discretizing advection term: (a) Iterative procedure of finding cell for interpolation, (b) neighboring cells of a node near the boundary.}\label{fig:advection_discretization}
	\end{figure}
	When all the vertices of the cell $\mathbf{C}_0$ belong to $\Omega$, $\mathbf{u}^{n}\left(\mathbf{x}\right)$
	is approximated by the quadratic essentially non-oscillatory (ENO) interpolation. Let $\mathbf{x}_{00},\mathbf{x}_{10},\mathbf{x}_{01},\mathbf{x}_{11}$ denote the lower left, lower right, upper left, and upper right vertices at the four corners of the cell, and $u_{00},u_{10},u_{01},u_{11}$ represent the value of $u$ at these four vertices, respectively. In addition, let $\left(\theta_x, \theta_y\right)= \frac{\mathbf{x}-\mathbf{x}_{00}}{\Delta x}$ for the cell size $\Delta x$. Then the quadratic ENO interpolation is defined as 
	\[
	\begin{aligned}	\uvec(x, y)&=(1-\theta_x)(1-\theta_y)\uvec_{00} +(1-\theta_x) \theta_y \uvec_{01}    +  \theta_x(1-\theta_y)\uvec_{10} +\theta_x\theta_y\uvec_{11} \\
		& - \frac{\theta_x(1-\theta_x)\Delta x^2}{2}D^{min}_{xx}\uvec- \frac{\theta_y(1-\theta_y)\Delta x^2}{2}D^{min}_{yy}\uvec,\end{aligned}
	\]
	where  
	\[D^{min}_{xx}\mathbf{u}= \text{minmod} \left(D_{xx} \uvec_{ij}\right) \]
	for $i,j=0,1$. Here, $minmod$ denotes a function that returns the value that has the minimum absolute value, and $D_{xx} u_{ij}$ is a discretization of $u_{xx}$ at $\mathbf{x}_{ij}$ according to the second-order finite differences in \cite{min2006second ,min2007second}. For the detailed formula, see Section \ref{subsec:Diffusion}.

	When one of the vertices does not belong to $\Omega$, interpolation is performed as follows.
	\begin{enumerate}
		\item Replace $\mathbf{C}_0$ with another cell that has at least one vertex belonging to $\Omega$.
		\item Approximate the velocity using the moving least-squares (MLS) method.
	\end{enumerate}
	When $\mathbf{C}_{0}$ has a vertex that belongs to $\Omega$, we skip step 1 above. By contrast, when there is no vertex that belongs to $\Omega$, we adopt the bisection algorithm introduced in \cite{coco2020multigrid}. First, we set $\mathbf{x}^0 = \mathbf{x}$. Then we repeat the sequence
	\[ \mathbf{x}^{k+1}=\mathbf{x}^k - \frac{\triangle x}{2}\mathbf{n}(\mathbf{x}^k)\] 
	until the cell containing $\mathbf{x}^k$ has a vertex that belongs to $\Omega$. After the iteration terminates, we set $\mathbf{C}_0$ to be the cell containing $\mathbf{x}^k$. This process is schematically illustrated in Figure \ref{fig:advection_discretization}(a).
	
	After the cell $\mathbf{C}_0$ is selected, we apply the MLS method. A linear polynomial $P_u$ with basis $\{1,x,y\}$ is constructed so as to minimize
	\[\sum_i w(\mathbf{x},\mathbf{x}_i )\left(P_u(\mathbf{x_i})-u_i\right)^2,\]
	where $w(\mathbf{x},\mathbf{y})=e^{-|\mathbf{x}-\mathbf{y}|^2}$ is a Gaussian weight function. To construct a sample point $\{ \left(\mathbf{x}_i,u_i\right)\}$, we first set $\mathbf{x}_0=\left(x_0,y_0\right)=argmax_{\mathbf{x}\in vertices(\mathbf{C}_0)} -\phi(\mathbf{x})$. Let $\mathbf{C}_j$ denote a cell such that $\mathbf{x}_0\in \partial \mathbf{C}_j$. Then we define $\{ \left(\mathbf{x}_i\right)\}=\mathbf{X}_{in} \cup \mathbf{X}_{bd}$, where
	\[\mathbf{X}_{in}=\{\mathbf{y}| \mathbf{y}\in Vertices(\mathbf{C}_j)\cap \Omega  \}\]
	and
	\[\mathbf{X}_{bd}=\{\mathbf{y}| \mathbf{y}\in \mathbf{C}_j \cap \Gamma  ,\  \left(\mathbf{x_0}-\mathbf{y}\right) \mathbin{\!/\mkern-5mu/\!} \mathbf{e}_i \text{ for some } i=1,2 \}\]
	for the standard basis $e_1=(1,0),e_2=(0,1)$. Note that $\mathbf{X}_{in}$ consists of neighboring grid points that belong to $\Omega$, and $\mathbf{X}_{bd}$ consists of the points that belong to $\Gamma$, where the line segment with the point and $\mathbf{x}_0$ lies along a Cartesian axis. In practice, we compute $\mathbf{X}_{bd}$ using the techniques in the GFM \cite{liu2000boundary,kang2000boundary}. For example, consider the situation illustrated in Figure \ref{fig:advection_discretization}(b). The grid point $\mathbf{X}_R$ is not in the fluid domain $\Omega$. Let $\phi_R$ and $\phi_0$ denote the values of the level set function $\phi$ at $\mathbf{x}_R$ and $\mathbf{x}_0$, respectively. Then $\hat{\mathbf{x}}_R=(x_0 +\theta_R d_R, y_0)$ is in $\mathbf{X}_{bd}$, where $\theta_R$ is computed as
	\[
	\theta_{R}=\frac{\mid\phi_{0}\mid}{\mid\phi_{0}\mid+\mid\phi_{R}\mid}.
	\]
	Sample points $\{ \left(\mathbf{x}_i,u_i\right)\}$ are then defined naturally; $u_i$ denotes the value of $u$ on $\mathbf{x}_i$ when $\mathbf{x}_i \in \Omega$, and $u_i$ is given by the boundary condition when $\mathbf{x}_i \in \Gamma$. After finding the polynomial $P_u$, we simply approximate 
	\[\uvec(\mathbf{x} ) \approx P_u(\mathbf{x}). \]
	A complete description of the MLS method is given in \cite{levin1998approximation,lancaster1981surfaces}.

	\subsection{\label{subsec:Diffusion}Diffusion term}
	
\begin{figure}
		\subfigure[]{\includegraphics[width=0.46\textwidth]{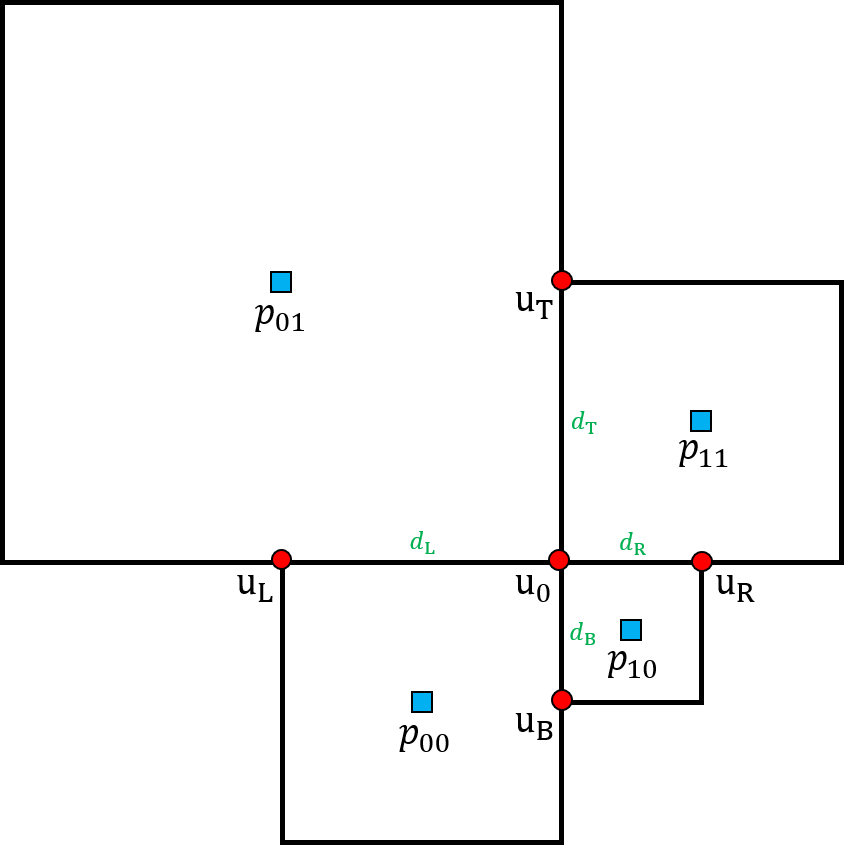}}$\ \ $
		\subfigure[]{\includegraphics[width=0.56\textwidth]{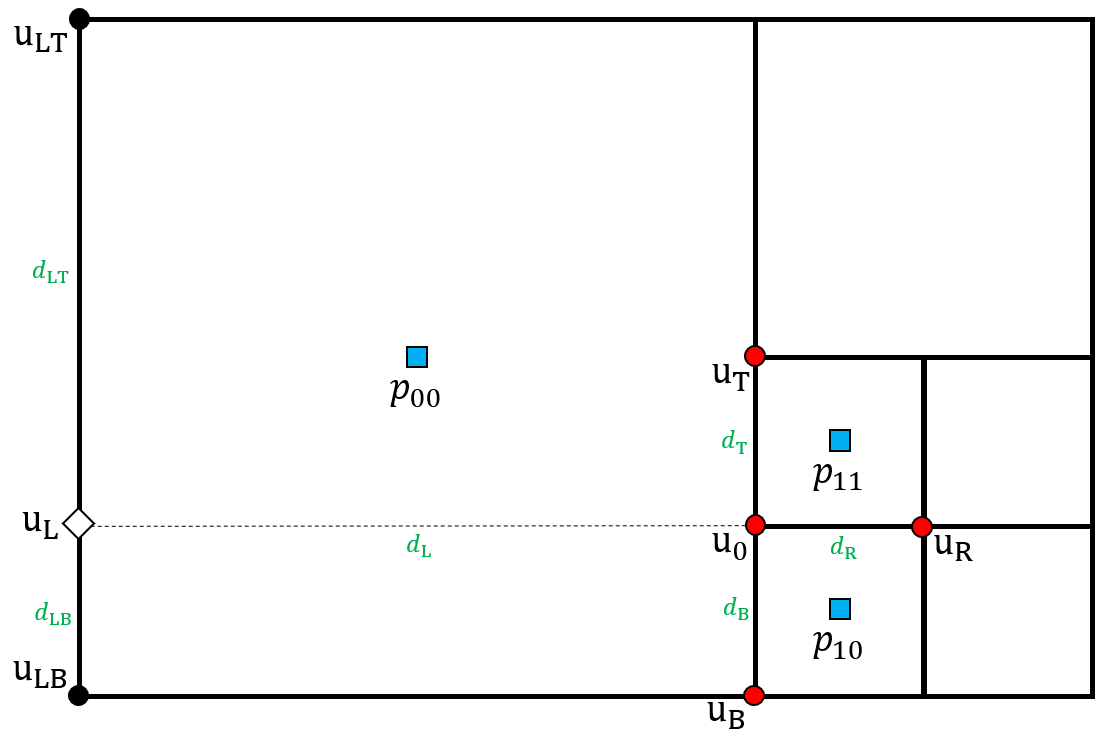}}
		\caption{\label{fig:stencil_stokes}Neighboring nodes of a regular node (a)
			and a T-junction node (b) in quadtree. The dashed line is imaginary.}
	\end{figure}

	The Laplacian velocity operator is discretized using the method introduced by Min et al. \cite{min2006supra}.
	We refer to a regular node as a node for which neighboring
	nodes in all the Cartesian directions exist; see Figure \ref{fig:stencil_stokes}(a). Let $\mathbf{u}_{R},\ \mathbf{u}_{L},\ \mathbf{u}_{T},\mathbf{u}_{B}$ be the values of $\mathbf{u}$ at the nodes immediately adjacent to $\mathbf{u}_{0}$, and
	$d_{R},\ d_{L},\ d_{T},d_{B}$ denote the distances between $\mathbf{x}_{0}$
	and the neighboring nodes to the right, left, top, and bottom, respectively. Then the Laplacian of $\mathbf{u}$ at the regular node $\mathbf{x}_{0}$ is discretized as
	\[
	\triangle\mathbf{u}_{0}=\left[\frac{\mathbf{u}_{R}-\mathbf{u}_{0}}{d_{R}}-\frac{\mathbf{u}_{0}-\mathbf{u}_{L}}{d_{L}}\right]\frac{2}{d_{R}+d_{L}}+\left[\frac{\mathbf{u}_{T}-\mathbf{u}_{0}}{d_{T}}-\frac{\mathbf{u}_{0}-\mathbf{u}_{B}}{d_{B}}\right]\frac{2}{d_{T}+d_{B}}.
	\]

	Nonuniform Cartesian grids produce hanging nodes, which appear only within the region of ﬁner resolution.
Referring to figure \ref{fig:stencil_stokes}(b), the Laplacian of $\mathbf{u}$ at a T-junction node is discretized by
	\[
	\triangle\mathbf{u}_{0}=\left[\frac{\mathbf{u}_{R}-\mathbf{u}_{0}}{d_{R}}-\frac{\mathbf{u}_{0}-\mathbf{u}_{L}}{d_{L}}\right]\frac{2}{d_{R}+d_{L}}+\left(1-\frac{d_{LT}d_{LB}}{2d_{L}\left(d_{R}+d_{L}\right)}\right)\left[\frac{\mathbf{u}_{T}-\mathbf{u}_{0}}{d_T}-\frac{\mathbf{u}_{0}-\mathbf{u}_{B}}{d_{B}}\right]\frac{2}{d_{T}+d_{B}}
	\]
	for 
	\[\uvec_L = \frac{d_{LB}}{d_{LB}+d_{LT}}\uvec_{LT}+ \frac{d_{LT}}{d_{LB}+d_{LT}}\uvec_{LB}.\]

	When the node is close to the boundary $\Gamma$, some of the
	adjacent nodes used in the discretization may lie outside
	the domain. For example, the right-hand neighbor $\mathbf{x}_R$ may not belong to $\Omega$. In this case, as explained in the previous section, we find the boundary point $\hat{\mathbf{x}}_R$ by computing $\theta_R$. Then, by replacing $d_R$ with $\theta_R d_R$ and $u_R$ with the boundary condition, the Laplacian can be discretized at the node close to the boundary.

	\subsection{\label{subsec:PressureGrad}Pressure gradient}
	
	We now consider the discretization of the pressure gradient operator. In the Arakawa
	B grid structure, the pressure variables are located at cell centers. Pressure nodes are usually not aligned horizontally or vertically with velocity nodes. Therefore, when $p_x$ is discretized at the node $\mathbf{x}_0$, the intermediate value $p_L,p_R$ and corresponding points $\mathbf{x}_L,\mathbf{x}_R$, which have the same $y$ coordinate as $\mathbf{x}_0$, are computed and used in the discretization.

\begin{figure}
		\subfigure[]{\includegraphics[width=0.46\textwidth]{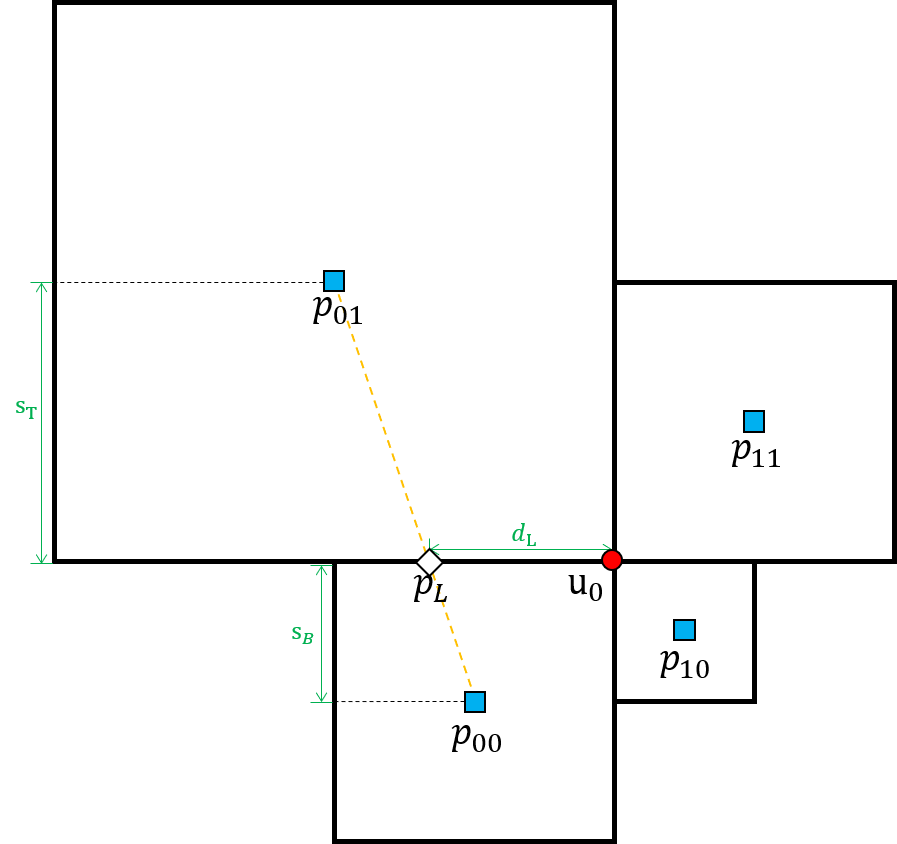}}$\ \ $
		\subfigure[]{\includegraphics[width=0.52\textwidth]{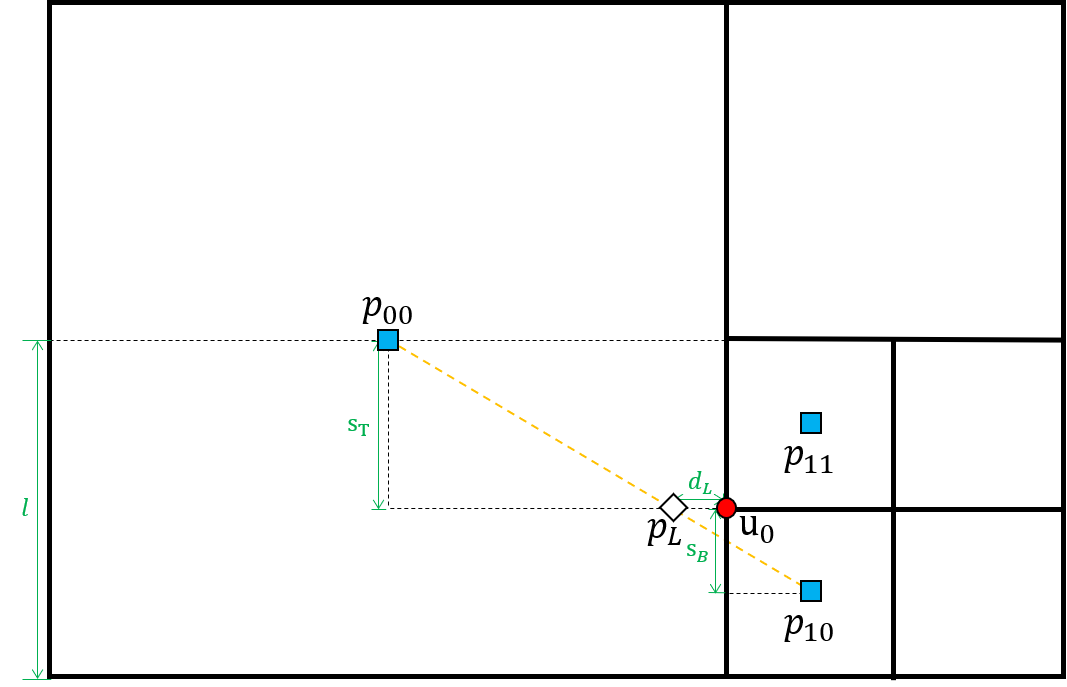}}
		\caption{\label{fig:grid_pressure}Neighboring nodes of a regular node (a)
			and a T-junction node (b) for discretization of the pressure gradient.
			Dashed lines are imaginary.}
	\end{figure}

	First, consider the case when $\mathbf{x}_0$ is a regular node. Referring to Figure \ref{fig:grid_pressure}(a), the pressure variables $p_{00},p_{01}$ on the left side of $\mathbf{x}_0$ are linearly interpolated to define 
	\[
	p_{L}\coloneqq\frac{s_{B}p_{01}+s_{T}p_{00}}{\theta_{T}+\theta_{B}}.
	\]
	The distance $d_L$ between $\mathbf{x}_L$ and $\mathbf{x_0}$ is then  
	\[
	d_{L}\coloneqq2\frac{s_{T}s_{B}}{s_{T}+s_{B}}.
	\]
	We can compute $p_R$ and $d_R$ in a similar manner. Then $p_x$ is discretized as
	\[
	p_{x}=\frac{p_{R}-p_{L}}{d_{R}+d_{L}}.
	\]
	
	For the T-junction node $\mathbf{x}_{0}$ in the cell configuration 
	shown in Figure \ref{fig:grid_pressure}(b), the preceding
	discretization cannot be directly used to approximate $p_{L}$ because there is only one adjacent cell to the left of $\mathbf{x}_0$. In Figure \ref{fig:grid_pressure}(b), the $y $ coordinate of the node with $p_{00}$ is greater than that of $\mathbf{x_0}$. Then, the line segment connecting the two pressure nodes with $p_{00}$ and $p_{10}$ intersects the line $y=y_0$ on the left side of $\mathbf{x_0}$. Therefore, $p_{00},p_{10}$ and its nodes are linearly interpolated to define $p_L$ and $d_L$ as
	\[
	p_{L}\coloneqq\frac{s_{B}p_{10}+s_{T}p_{00}}{\theta_{T}+\theta_{B}},
	\]
	and
	\[
	d_{L}\coloneqq\frac{\left(l-s_{T}\right)s_{B}}{s_{T}+s_{B}}.
	\]
	When the $y $ coordinate of the node with $p_{00}$ is less than that of $\mathbf{x_0}$, $p_{11}$ and its node are used instead of $p_{10}$.
	
	The discretization of $\frac{\partial p}{\partial y}$ is defined similarly. For nodes near the boundary, whenever one of the vertices of the cell belongs to $\Omega$, ghost pressure nodes are defined at the cell. As a result, the discretization described above can be directly applied to the nodes near the boundary.
	\subsubsection{Alternative discretization in finite volume approach} \label{subsubsec:pressure_fv}
	There is a relatively simple approach to discretizing the pressure gradient operator by the finite volume method. It does not require that an intermediate value be computed, and the extension to octrees is natural and simple. As shown in Figure \ref{fig:stencil_stokes}, we define the control volume $V$ of the grid node $\mathbf{x}_0=(x_0,y_0)$ as $V=[x_0-\frac{d_L}{2},x_0+\frac{d_R}{2}]\times[y_0-\frac{d_L}{2},y_0+\frac{d_T}{2}]$. Note that the control volume is consistent with the discretization of the Laplacian discussed in Section \ref{subsec:Diffusion}. Then the divergence theorem yields 
	\[\int_V p_x  d\mathbf{x}=\int_{\partial V } \left(p,0\right) \cdot \normal_V ds,\] where $\normal_V$ is the outward vector normal to $\partial V$.
	Let $C_i$ denote the cell intersecting the control volume $V$ and $p_i$ denote a pressure variable on the cell $C_i$. By interpreting the pressure variable as a piecewise constant function on the cell, the integral can be approximated in the finite volume sense as
	\begin{equation*}
		\int_{\partial V } \left(p,0\right) \cdot \normal_V ds=\sum_i   \left(p_i,0\right) \cdot \normal_V |\partial V \cap C_i|.
	\end{equation*}
	Then, $p_x$ is discretized by dividing this approximation by the area of the control volume $V$. For example, $\nabla p$ in Figure \ref{fig:stencil_stokes}(a) is discretized by
	\begin{equation*}
		\begin{aligned}
			p_x= \frac{4}{(d_R+d_L)(d_T+d_B)} \left(  (p_{11}-p_{01}) \frac{d_T}{2} +(p_{10}-p_{00}) \frac{d_B}{2}\right), 	\\
			p_y= \frac{4}{(d_R+d_L)(d_T+d_B)} \left(  (p_{11}-p_{10}) \frac{d_R}{2} +(p_{01}-p_{00}) \frac{d_L}{2}\right) .	
		\end{aligned}
	\end{equation*}
	On the hanging node shown in Figure \ref{fig:stencil_stokes}(b), the discrete pressure gradient is defined as
	\begin{equation*}
		\begin{aligned}
			p_x= \frac{4}{(d_R+d_L)(d_T+d_B)} \left(  (p_{11}-p_{00}) \frac{d_T}{2} +(p_{10}-p_{00}) \frac{d_B}{2}\right), 	\\
			p_y= \frac{4}{(d_R+d_L)(d_T+d_B)} \left(  (p_{11}-p_{10}) \frac{d_R}{2} +(p_{00}-p_{00}) \frac{d_L}{2}\right) .	
		\end{aligned}
	\end{equation*}
	On the nodes near the boundary, the control volume is modified as described in Section \ref{subsec:Diffusion}. If $\mathbf{x}_R\notin \Omega$, $\theta_R$ is computed such that $\mathbf{x}_0 + (\theta_R d_R,0)\in \Gamma$, and $d_R$ is replaced with $\theta_R d_R$. Note that the finite volume approach is motivated by \cite{purvis1979prediction}. Although it is much simpler than the finite difference approach, it is found numerically that the order of convergence is reduced to first order at higher Reynolds number. The degradation of the order of accuracy is reported in Table \ref{tab:pressure_grad}. For this reason, the finite difference approach is used in most of the numerical experiments. 
	
	\subsection{Divergence-free equation}
\begin{figure}
		\centering{}\includegraphics[width=0.9\textwidth]{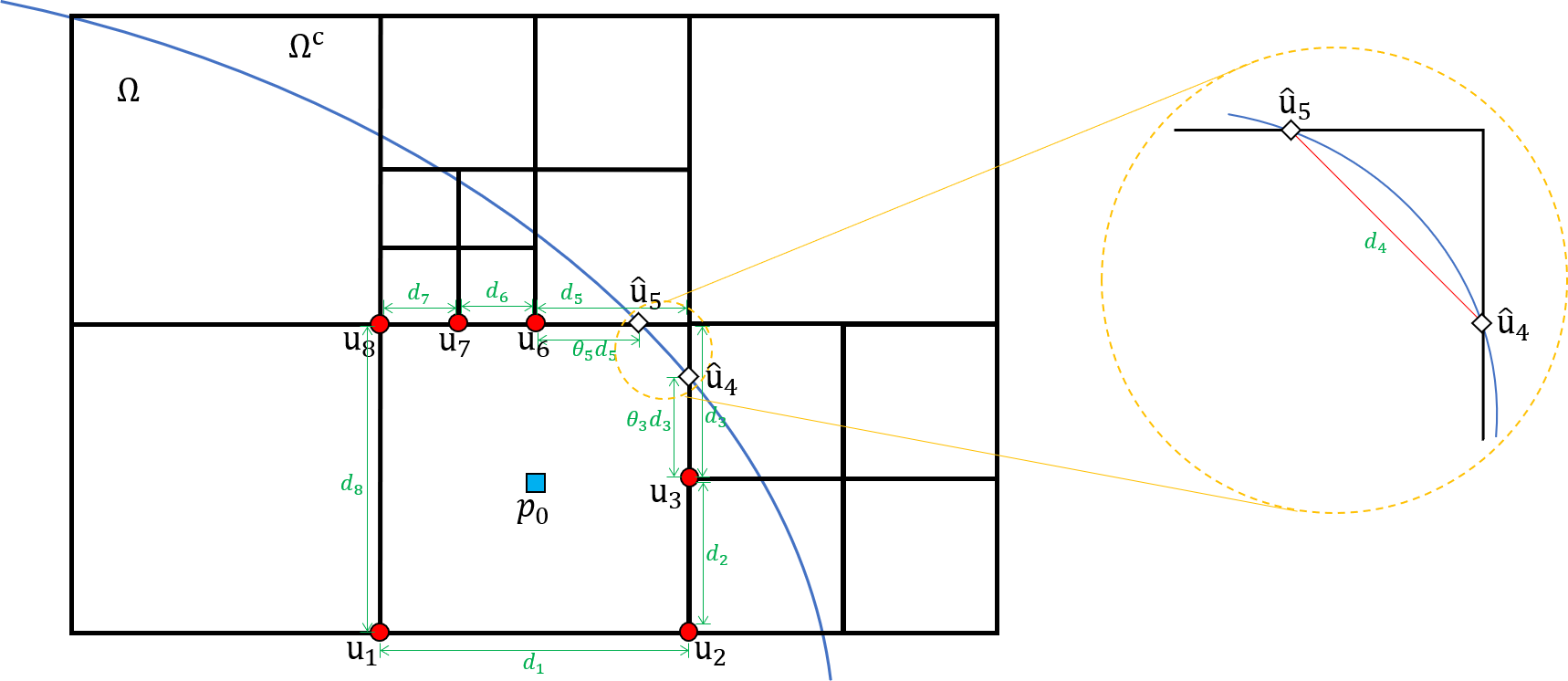}\caption{\label{fig:div}Fictitious cell configuration used in the discretization
			of divergence operator.}
	\end{figure}

	The divergence of the velocities is discretized at each cell. Let $\mathbb{\mathbf{C}}_{0}$ denote a
	cell, and define a set $\mathbf{X}=\{\mathbf{x}_i\}_{i=1,\ldots,n}\subset \partial \mathbf{C}_0$ as 
	\[
	\mathbf{x} \in \mathbf{X} \iff \mathbf{x}\in \partial \mathbf{C}_0 \text{ is a grid node } \quad \text{or} \quad \mathbf{x} \in \Gamma \cap \partial \mathbf{C}_0,
	\] 
	where $\mathbf{x}_i$ are ordered in the counterclockwise direction. Then $\partial \left(\mathbf{C}_0 \cap \Omega\right)$ is discretely approximated as the union of line segments $\Gamma_i$:
	\[
	\partial \left(\mathbf{C}_0 \cap \Omega\right) \approx \sum_{i=1}^n \Gamma_i, 
	\]
	where $\Gamma_{i}$ is a line segment connecting $\mathbf{x}_i$ and $\mathbf{x}_{i+1}$. Here, $\mathbf{x}_{n+1}=\mathbf{x}_0$. The normal vector $\normal_i$ corresponding to $\Gamma_i$ is obtained by $90^\circ$ clockwise rotation of $\left(\mathbf{x}_{i+1}-\mathbf{x}_i\right)/\|\mathbf{x}_{i+1}-\mathbf{x}_i\|$.
	Then the divergence of the velocities is discretized in the finite volume sense using the divergence theorem:
	\begin{equation}\label{eq:div_discretization}
		\begin{aligned}
			\intop_{\mathbf{C}_0 \cap \Omega}\nabla\cdot\mathbf{u}dxdy & =\intop_{\partial\left(\mathbf{C}_0 \cap \Omega\right)}\mathbf{u}\cdot \normal ds=\Sigma_{i=1}^{F}\intop_{\Gamma_{i}}\mathbf{u}\cdot \normal_{i}ds =\Sigma_{i=1}^{n} \left(\frac{\uvec_{i}+\uvec_{i+1}}{2}\cdot \normal_{i}\right)|\Gamma_i|.
		\end{aligned}
	\end{equation}
	Precise discretization is demonstrated
	by considering the general cell configuration shown in Figure \ref{fig:div}.
	The set $\mathbf{X}$ consists of eight points; $\mathbf{x}_1,\mathbf{x}_2,\mathbf{x}_3,\mathbf{x}_6,\mathbf{x}_7,\mathbf{x}_8$ are grid nodes that belongs to $\Omega$, and $\mathbf{x}_4,\mathbf{x}_5$ are boundary points.
	As shown in the previous sections, $\theta_3,\theta_5$, which define the position of the boundary, are computed by linear interpolation of $\phi$. Consequently, the divergence
	operator is discretized as follows:
	
	\begin{align*}
		\frac{1}{\left|\mathbf{C}_{0}\right|}\intop_{\Omega_{0}}\nabla\cdot\mathbf{u}dxdy & =\frac{1}{\left|\mathbf{C}_{0}\right|}\Sigma_{i=0}^{7}\left(\mathbf{u}_{i}\cdot \normal_{i}\right)\mid\Gamma_{i}\mid\\
		& =\frac{1}{\left|\mathbf{C}_0\right|}\Sigma_{i=0}^{7}\left(\frac{\mathbf{u}_{i}+\mathbf{u}_{i+1}}{2}\cdot \normal_{i}\right)d_{i}\\
		& =\left\{ \left(\frac{d_{2}}{2}u_{2}+\frac{d_{2}+\theta_{3}d_{3}}{d}u_{3}+\frac{\theta_{3}d_{3}}{2}\hat{u}_{4}\right) +d_{4}\frac{\hat{\mathbf{u}}_{4}+\hat{\mathbf{u}}_{5}}{2}\cdot \normal_{4} \right.\\
		& +\left(\frac{\theta_{5}d_{5}}{2}\hat{v}_{5}+\frac{\theta_{5}d_{5}+d_{6}}{2}v_{6}+\frac{d_{6}+d_{7}}{2}v_{7}+\frac{d_{7}}{2}v_{8}\right) -\left(\frac{d_{8}}{2}u_{8}+\frac{d_{8}}{2}u_{1}\right) \\
		& \left. -\left(\frac{d_{1}}{2}v_{1}+\frac{d_{1}}{2}v_{2}\right)\right\} \frac{1}{\left|\mathbf{C}_{0} \right|.}
	\end{align*}
	In the last expression, parentheses indicate boundary integrand
	values in the counterclockwise direction from the right edge.

	Note that the discretization of the divergence-free equation \eqref{eq:div_discretization} on a uniform Cartesian grid is equivalent to the discretization of the bilinear-constant velocity--pressure finite element method. It is also known as the $Q1-P0$ element and exhibits pressure checkerboard instability. To reduce the instability, stabilizing techniques have been studied in the finite element community. For example, the penalty method \cite{hughes1979finite} reduces the instability by decoupling the velocity and pressure variables. Other examples include local and global stabilizers \cite{bochev2006stabilization,silvester1990stabilised}, which add an artificial term to the divergence-free equation, and the pressure gradient projection method \cite{blasco2001space,codina2000stabilized}, which was inspired by fractional step methods. Among many stabilizers, we extend the global stabilizer in \cite{silvester1990stabilised} to the quadtree grid as follows. Using the same notation as in \eqref{eq:div_discretization}, let $\mathbf{C}_i$ be a neighboring cell to $\mathbf{C}_0$ such that $\mathbf{C_0}\cap \mathbf{C_i}=\Gamma_i$, and let $p_i$ denote the pressure variable in the cell $\mathbf{C}_i$. Then the stabilizer is defined as 
	\[ 	
	S( p) = \epsilon \frac{1}{|\mathbf{C}_0|}\Sigma_{i=1}^{n} \left(\left(p_i - p_0\right)\cdot \normal_{i}|\Gamma_i|\right)
	\]
	for $\epsilon=\Delta x_s^2$. In Section \ref{subsubsec:irregulardomain}, the need for the stabilizer is described.

\begin{figure}[h]
		\centering{}
		\subfigure[ Checkerboard type kernel with rank 2]{\includegraphics[width=0.4\textwidth]{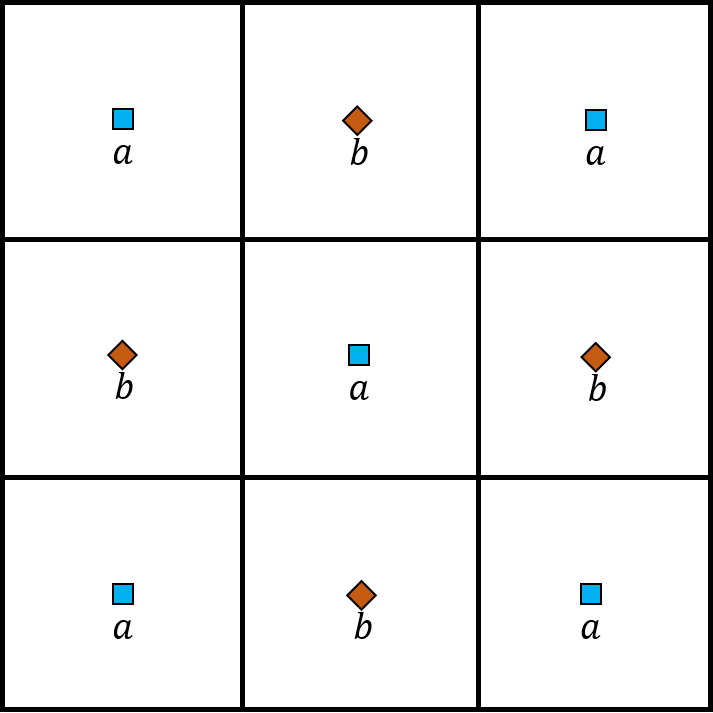}}
		$\ \ \ \ \ \ \ \ \ \ $ 
		\subfigure[ Constant pressure kernel induced by discretization at the hanging node]{\includegraphics[width=0.4\textwidth]{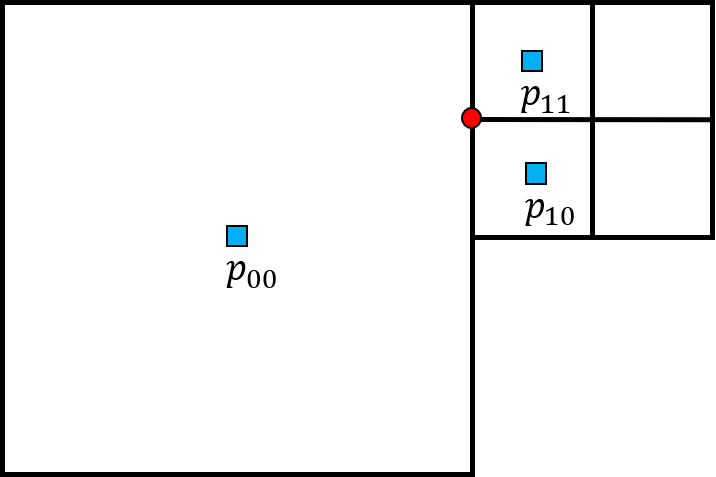}}
		\caption{Illustrations of pressure variable that belongs to kernel space of gradient operator.}\label{fig:kernel}
	\end{figure}

	\paragraph{Remarks}
	Note that the presented pressure gradient operator is rank-2 deficient on a uniform Cartesian grid. For any real numbers $a,b\in \mathbb{R}$, if the pressure variables are given as illustrated in Figure \ref{fig:kernel}(a), $\nabla p=\mathbf{0}$ at the interior vertices. However, for quadtree grids, this problem is resolved. The discretization $\frac{\partial p}{\partial y}=0$ at the vertex marked with a red circle in Figure \ref{fig:kernel}(b) yields $p_{11}=p_{10}$; consequently, the checkerboard pressure distribution shown in Figure \ref{fig:kernel}(a) no longer belongs to the kernel space unless $a=b$. Therefore, on nonuniform Cartesian grids, the $p\to \nabla p$ operator is rank-1 deficient. To show that the divergence operator $\nabla \cdot$ is rank-1 deficient, let us denote the divergence operator by $\mathbf{D}$, and let $\mathbf{M}$ be a diagonal matrix that maps a pressure variable to pressure variable and is defined as
	\[\mathbf{M}_{ii}= |\mathbf{C}_i|, \]
	where $|\mathbf{C}_i|$ is the area of the cell $\mathbf{C}_i$. It can be shown that $\left(\mathbf{M} \circ \mathbf{D}\right)^T$, like the pressure gradient operator, is rank-2 deficient on a uniform grid and is reduced to a rank-1 deficient matrix on a quadtree grid. Note that the kernel of $\left(\mathbf{M} \circ \mathbf{D}\right)^T$ is spanned by the vector of all those at cell centers. Because the kernel space of the transpose matrix has been determined, we project the right-hand side onto the range space before solving the linear system. Let $r_i$ denote a variable located at the cell center on the right-hand side. We define $\hat{r}_i$ as 
	\[\hat{r}_i=r_i - \frac{1}{|\mathbf{C}_i|}\frac{\sum_j r_j | \mathbf{C}_j |} {\sum_j |\mathbf{C}_j |}.\]
	Then $r_i$ is replaced by $\hat{r}_i$.
	
	Let $\mathbf{G}$ denote the linear operator $p \to \left(\int_V p_x , \int_V p_y\right)$ in the finite volume discretization introduced in Section \ref{subsubsec:pressure_fv}. Then it can be shown that the discrete gradient operator is the negative adjoint of the discrete divergence, that is, $-\mathbf{G}=\left(\mathbf{M} \circ \mathbf{D}\right)^T$. However, this property does not hold for the proposed discretization of the pressure gradient by the finite difference method. Despite the desirable matrix property, we do not adopt the finite volume approach because of the lower order of convergence at high Reynolds number.

	\subsection{Extension to three spatial dimensions}
	The proposed method can be extended to three dimensions on cubical domains with octree data structure. The discretization of the advection and diffusion terms can be naturally extended to three spatial dimensions as in \cite{min2006second,min2006supra}. However, the discretization of the pressure gradient and incompressibility condition require modifications. 
	
	We first describe the discretization of the pressure gradient. Consider the discretization of the derivative of the pressure in the $x$ direction at a node $\mathbf{x}_{0}$. As shown in Section \ref{subsec:PressureGrad}, a standard approach is to compute the intermediate values $p_{R}$ and $p_{L}$ with their locations $\mathbf{x}_{R}$ and $\mathbf{x}_{L}$, respectively, which have the same $y$ and $z$ coordinates as $\mathbf{x}_{0}$. In contrast to the two-dimensional case, three or four cells are used to approximate the intermediate values. If four cells are chosen, bilinear interpolation on the $yz$ plane is applied to approximate the intermediate values. If three cells are used, linear interpolation is applied instead of bilinear interpolation. For example, suppose we have chosen three cells whose center is $\left(x_i,y_i,z_i\right)$ for $i=1,2,3$, and let us denote by $p_i$ the pressure at the corresponding cell. Linear polynomials $P_p$ and $P_x$ are constructed on the $yz$ plane so that
	\begin{equation*}
		\begin{aligned}
			P_p (y_i, z_i)=p_i,\quad 	P_x(y_i,z_i)=x_i.
		\end{aligned}
	\end{equation*}
	Then $p_R$ and $x_R$ are approximated as 
	$	p_R=P_p (y_0, z_0)$ and $x_R=P_x(y_0,z_0)$.

	The remaining difficulty is the suitable choice of cells, which depends on the number of adjacent cells in the positive $x$ direction with respect to $\mathbf{x}_0$, that is, cells for which the maximum $x$ coordinate is larger than $x_0$.
	
\begin{figure}
		\centering{}
		\subfigure[Two cells are located in positive $x-$direction]{\includegraphics[width=0.4\textwidth]{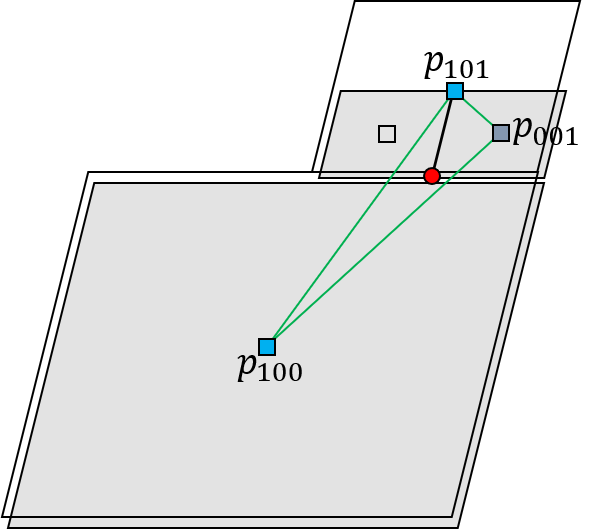}}
		$\ \ \ \ \ \ \ \ \ \ $ 
		\subfigure[One cell is located in positive $x-$direction]{\includegraphics[width=0.4\textwidth]{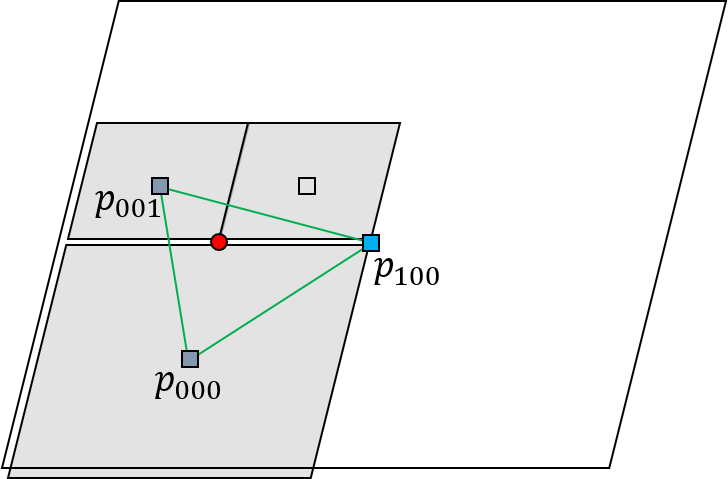}}
		\caption{\label{fig:pressure3d} General cell configurations for three cases considered for the discretization of the pressure gradient. These represent projected grids onto $yz$-plane. Cells in negative $x-$directions are illustrated by shaded regions and the red dot is the evaluation node $\mathbf{x}_{0}$.}
	\end{figure}

	When there are three or four cells in the positive $x$ direction, the cells are used to construct linear or bilinear polynomials, respectively. However, when there are only two cells or one cell in the positive $x$ direction, we use the cell information from the negative $x$ direction in the linear interpolation. Specifically, we select cells so that the point $\left(y_0,z_0\right)$ lies within a triangle with vertices $(y_1,z_1),(y_2,z_2)$, and $(y_3,z_3)$. Figure \ref{fig:pressure3d} shows an example.
	
\begin{figure}
		\centering{}\includegraphics[width=0.8\textwidth]{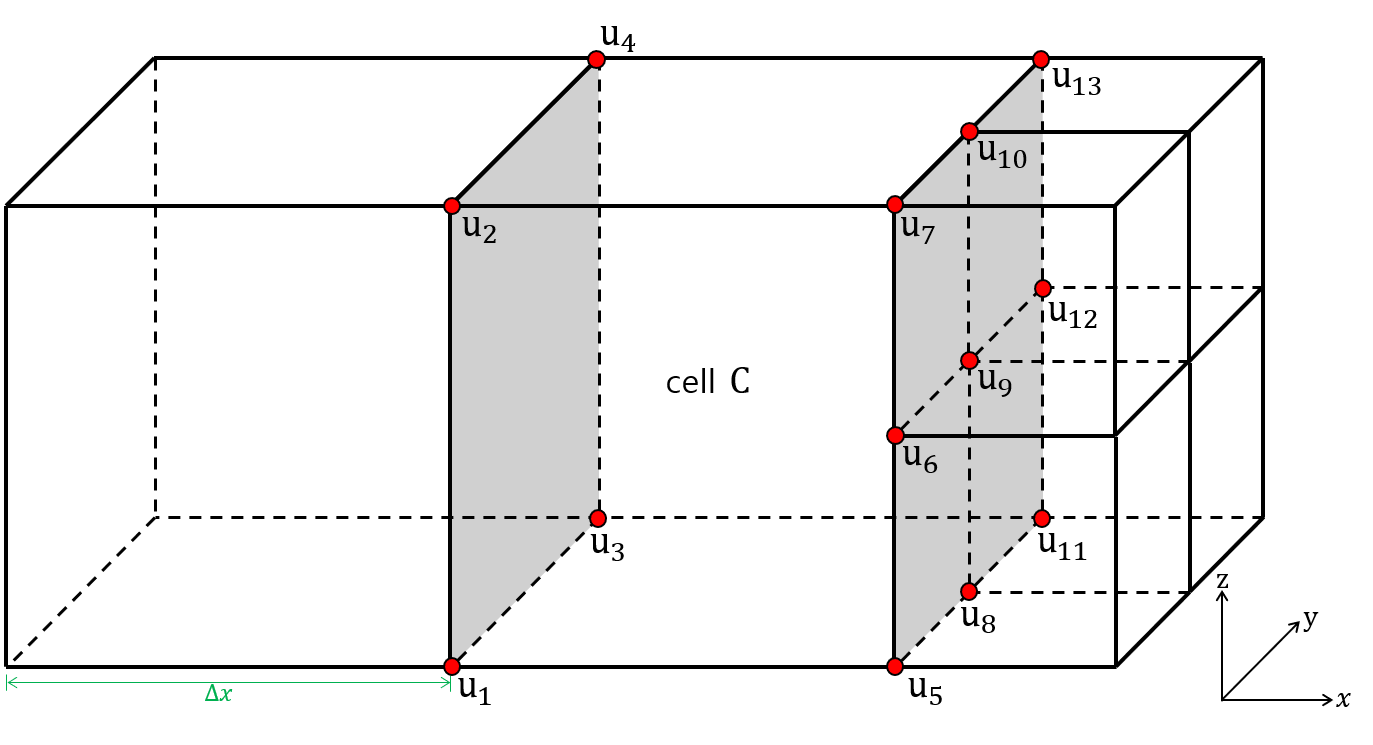}\caption{\label{fig:div_3d}General configuration of an octree grid used for the discretization of the divergence operator in the cell center of $\mathbf{C}$.}
	\end{figure}

	Now, consider the discretization of the divergence operator. As in two spatial dimensions, we use the divergence theorem to approximate the divergence operator in a cell $\mathbf{C}$. In three dimensions, the line fraction is replaced with a face fraction, and thus the line integration becomes an area integration. The discrete form of the divergence operator is given by
	\[
	\nabla\cdot\mathbf{u}= \frac{1}{\mid\mathbf{C}\mid}\Sigma_{F\in\mathcal{F}\mathbf(C)}\mid F\mid\left(\mathbf{u\cdot n}\right)_{F},
	\]
	where $\mathcal{F}\left(\mathbf{C}\right)$ denotes the set of all faces $F$ of the cell $\mathbf{C}$, and $\mathbf{n}$ is the outward unit normal to each face. $\mid\mathbf{C}\mid$ represents the volume of the cell $\mathbf{C}$, whereas $\mid\mathbf{F}\mid$ is the area of the cell face. Moreover, $\left(\mathbf{u\cdot n}\right)_{F}$ is discretized by the average of $\mathbf{u}$ at the four corner vertices of the face.
	As shown in Figure \ref{fig:div_3d}, the discretization of the area integration with the normal vector parallel to the $x$ axis is calculated as
	\begin{align*}
		\mid\mathbf{C}\mid\frac{\partial u}{\partial x} & \approx 
		\left(\frac{\triangle x}{2}\right)^{2}\frac{u_{5}+u_{6}+u_{8}+u_{9}}{4}+\left(\frac{\triangle x}{2}\right)^{2}\frac{u_{6}+u_{7}+u_{9}+u_{10}}{4}\\
		& + \left(\frac{\triangle x}{2}\right)^{2}\frac{u_{8}+u_{9}+u_{11}+u_{12}}{4} + \left(\frac{\triangle x}{2}\right)^{2}\frac{u_{9}+u_{10}+u_{12}+u_{13}}{4} \\
		& -\triangle x^{2}\frac{u_{1}+u_{2}+u_{3}+u_{4}}{4}.
	\end{align*}
	
	\section{Numerical results}
	
	In this section, we present the results of several numerical experiments to demonstrate
	the performance of the proposed method. We begin with numerical tests
	that confirm the second-order accuracy of the proposed method. Next, we provide numerical evidence that illustrates the performance
	of the method for standard benchmark problems. We implemented
	all the numerical experiments in C++ on a personal computer. To solve
	the saddle-point system (\ref{eq:saddle_system}) generated by our
	discretization, we used the generalized minimal residual method with an incomplete
	LU preconditioner.
	
	\subsection{Accuracy tests}
	
	To estimate the order of accuracy of the suggested method,
	we compare approximated numerical solutions to an analytic solution.
	For a two-dimensional example, consider a single-vortex problem whose exact solutions are given
	by
	\begin{equation}\label{eq:single_vortex}
		\begin{aligned}
			u(x,y,t) & =-\cos x\sin y\cos t,\\
			v(x,y,t) & =\cos y\sin x\cos t,\\
			p(x,y,t) & =-\frac{1}{4}\cos^{2}t\left(\cos2x+\cos2y\right).
		\end{aligned}
	\end{equation}
	The corresponding forcing term $f=\left(f_{1}, f_{2}\right)$ is given by
	\begin{align*}
		f_{1}(x,y,t) & =\cos x\sin y\left(\sin t-2\frac{1}{Re}\cos t\right),\\
		f_{2}(x,y,t) & =-\cos y\sin x\left(\sin t-2\frac{1}{Re}\cos t\right).
	\end{align*}
	The computational domain is taken to be $\Omega=\left[-\frac{\pi}{2},\frac{\pi}{2}\right]^{2}$,
	and the terminal time is set to $t=\pi$. 
	
	\subsubsection{\label{subsec:Random}Random quadtree}
	
\begin{figure}
		\centering{}
		\subfigure[]{\includegraphics[width=0.4\textwidth]{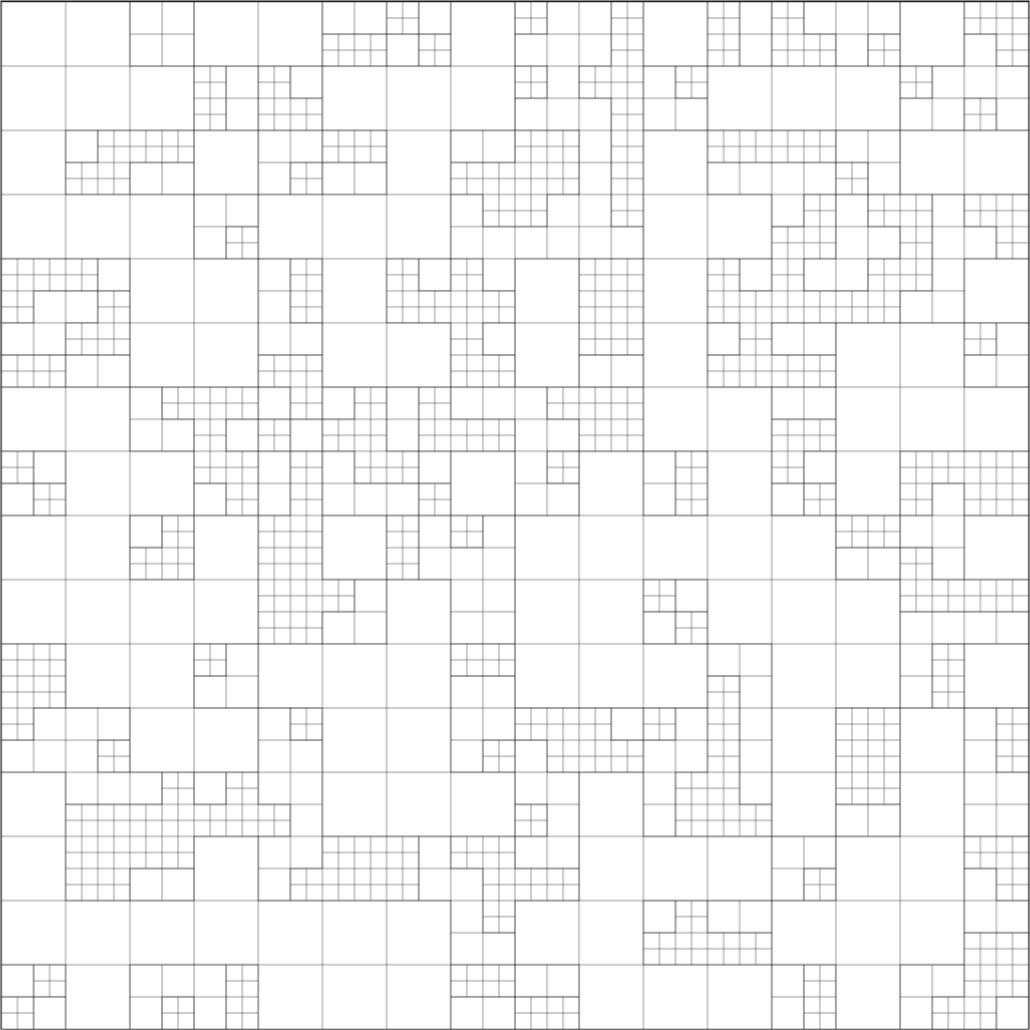}}$\quad\quad\quad$
		\subfigure[]{\includegraphics[width=0.4\textwidth]{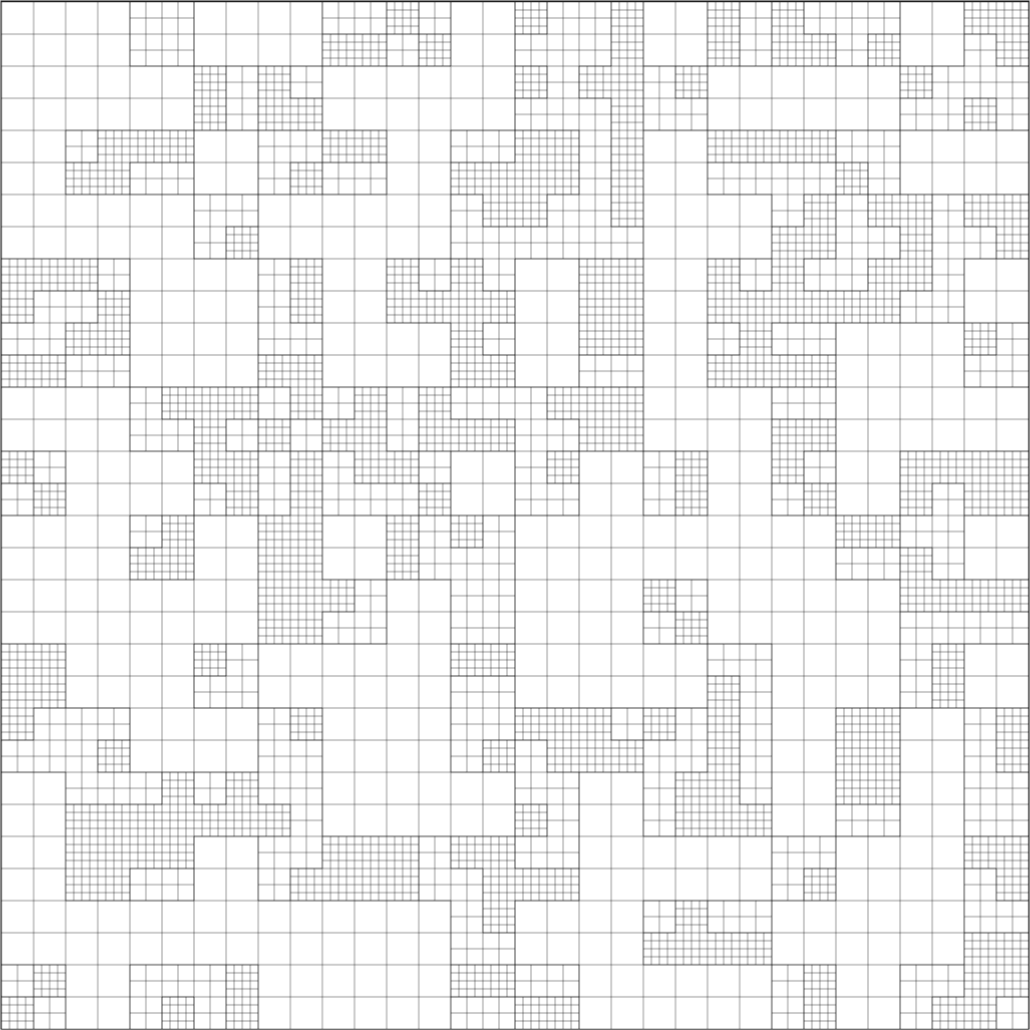}}
		\caption{A example of random quadtree of level (a) $4/6$ and refined (b) $5/7$ used in example \ref{subsec:Random}.} \label{fig:random_grid}
	\end{figure}

\begin{figure}
		\centering{}
		\mbox{
			\subfigure[]{\includegraphics[width=0.3\textwidth]{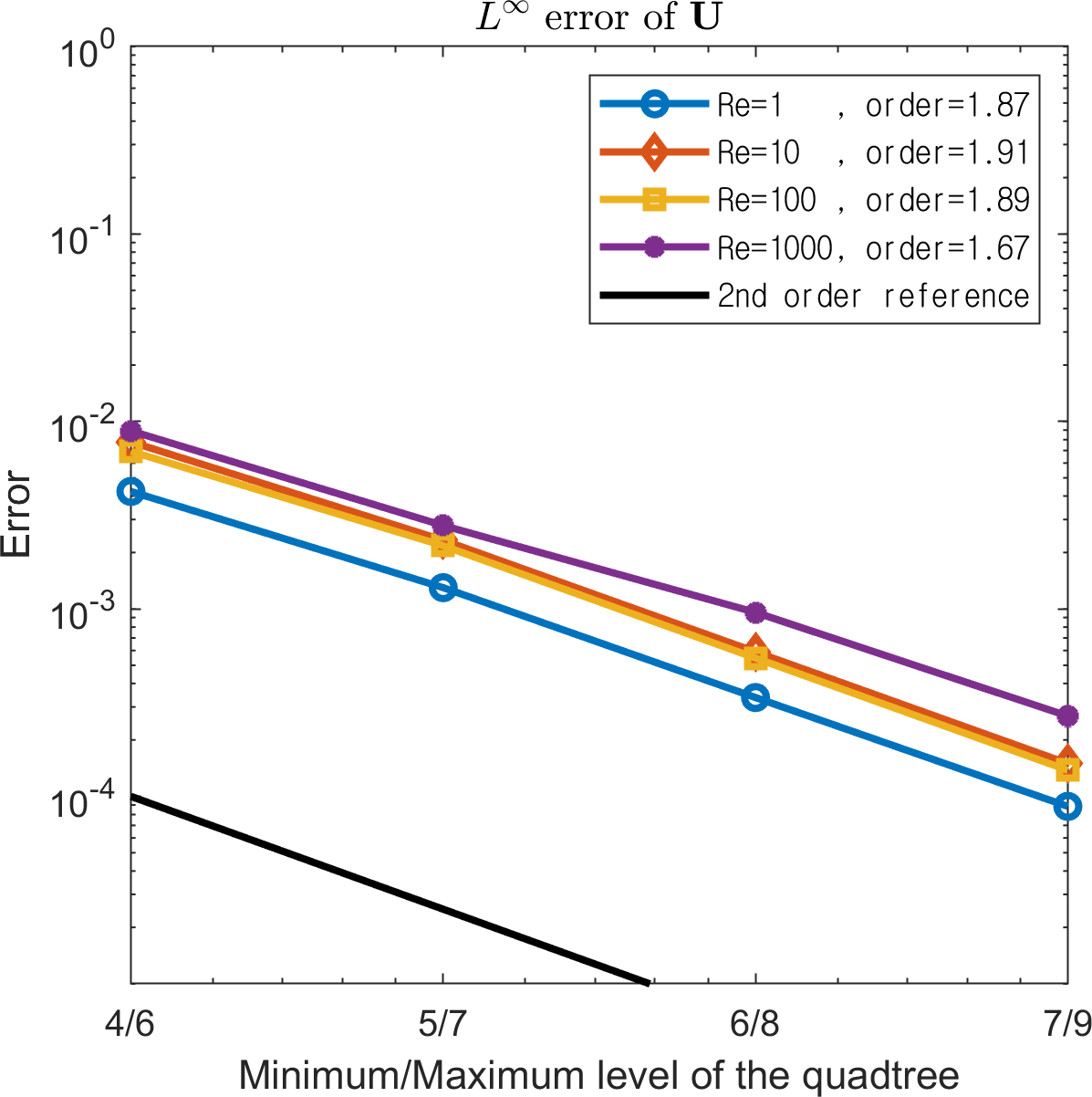}}$\ $
			\subfigure[]{\includegraphics[width=0.3\textwidth]{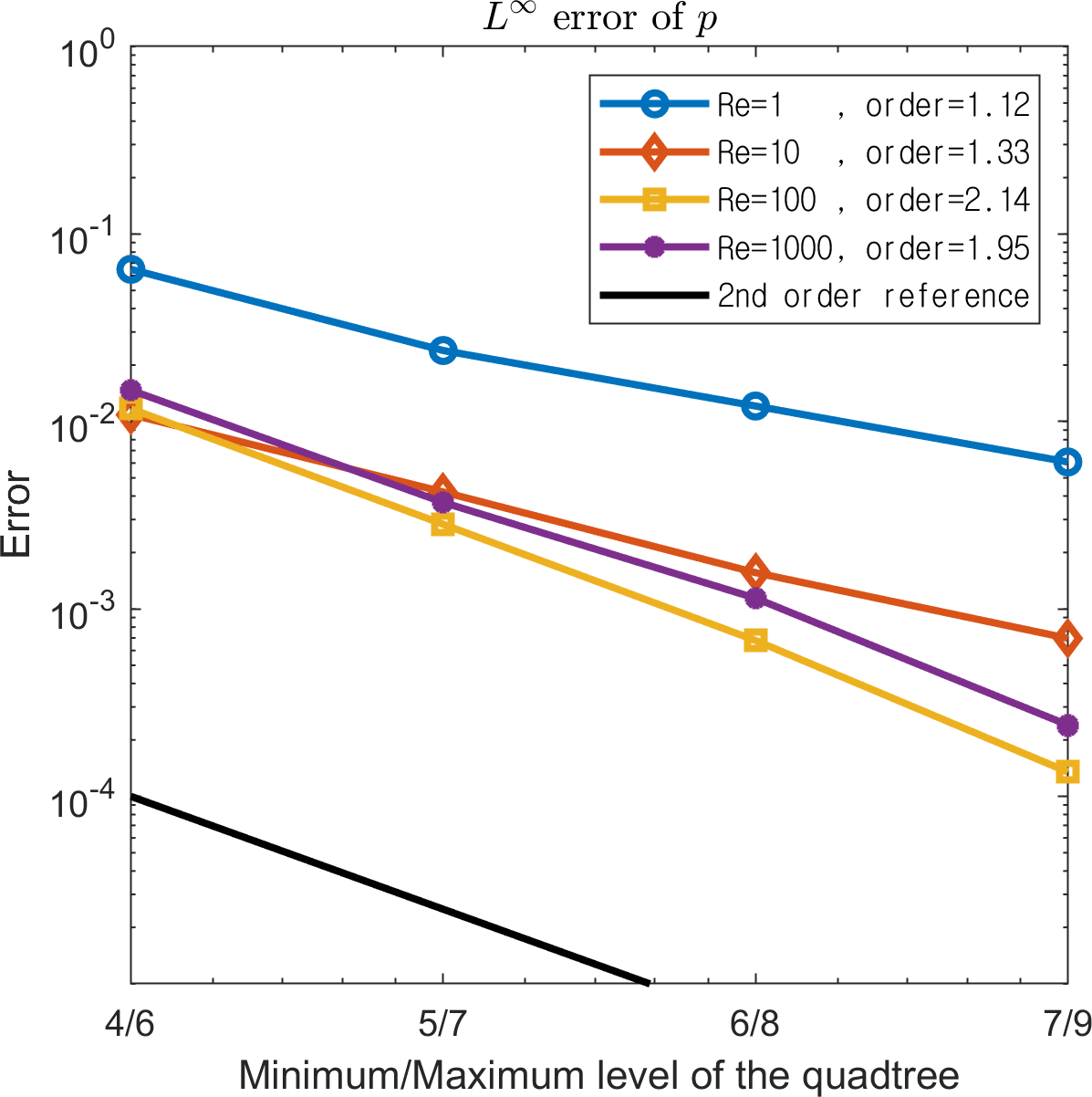}}$\ $
			\subfigure[]{\includegraphics[width=0.3\textwidth]{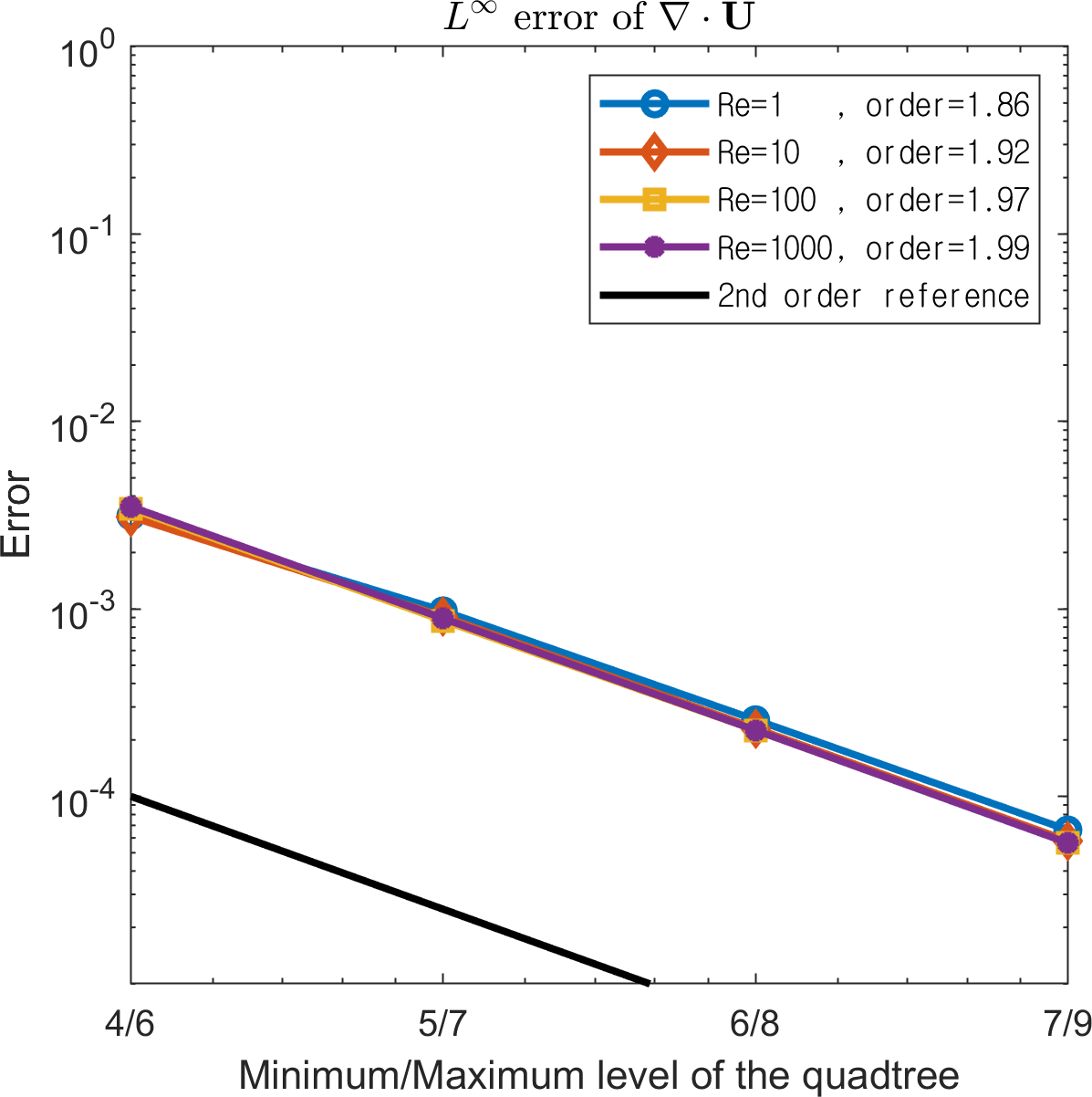}}
		}\\
		\centering{}
		\mbox{
			\subfigure[]{\includegraphics[width=0.3\textwidth]{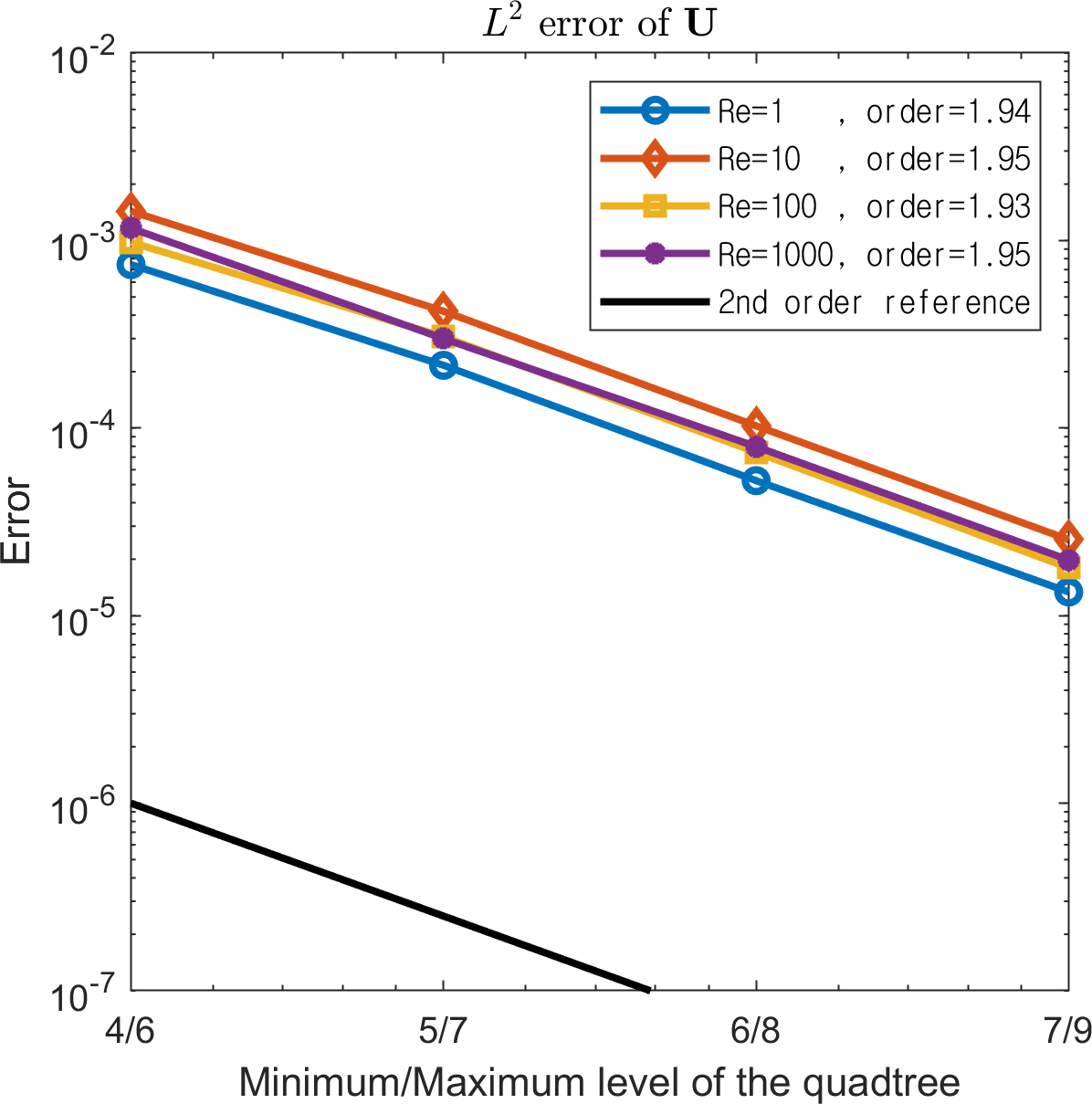}}$\ $
			\subfigure[]{\includegraphics[width=0.3\textwidth]{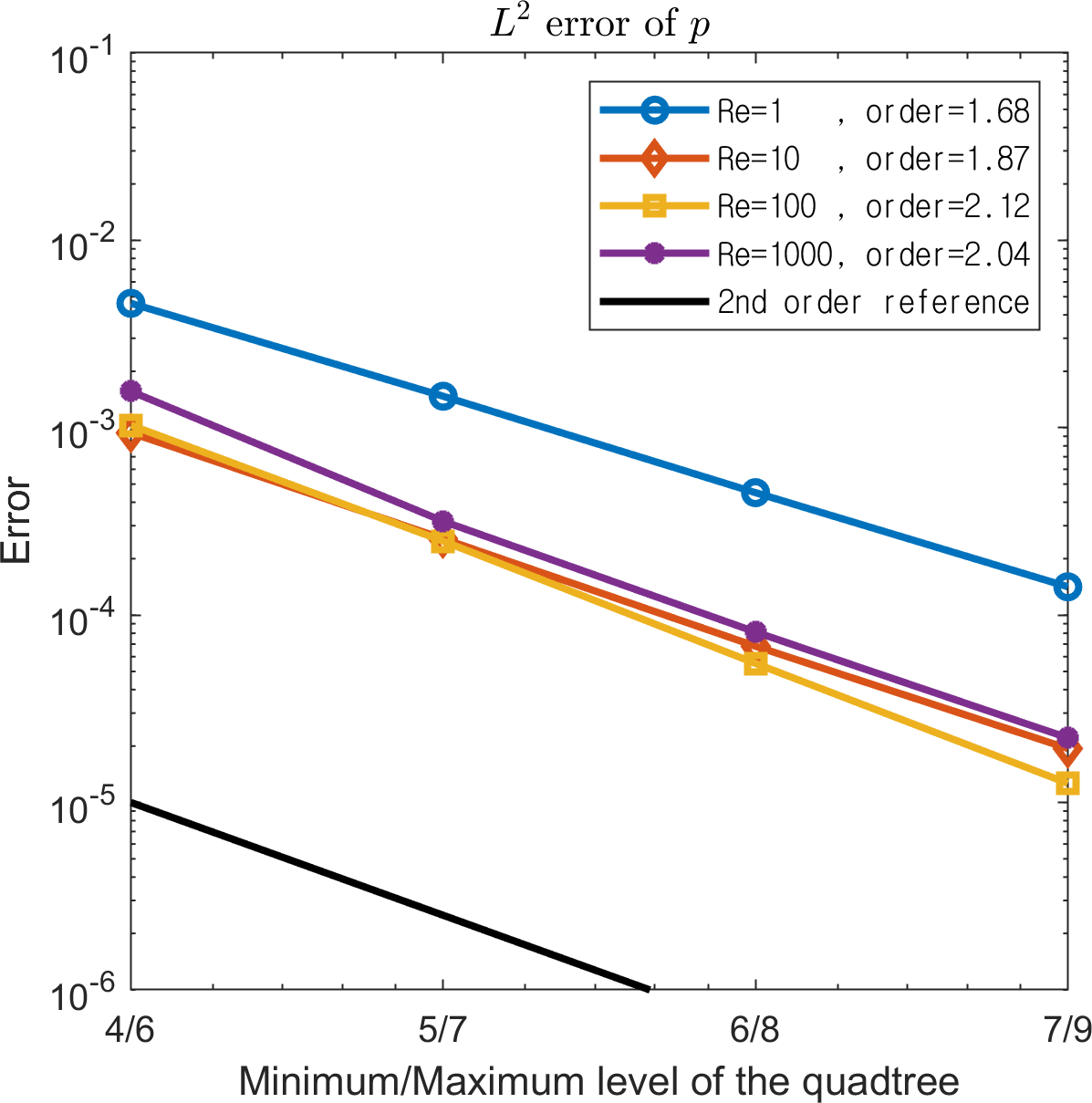}}$\ $
			\subfigure[]{\includegraphics[width=0.3\textwidth]{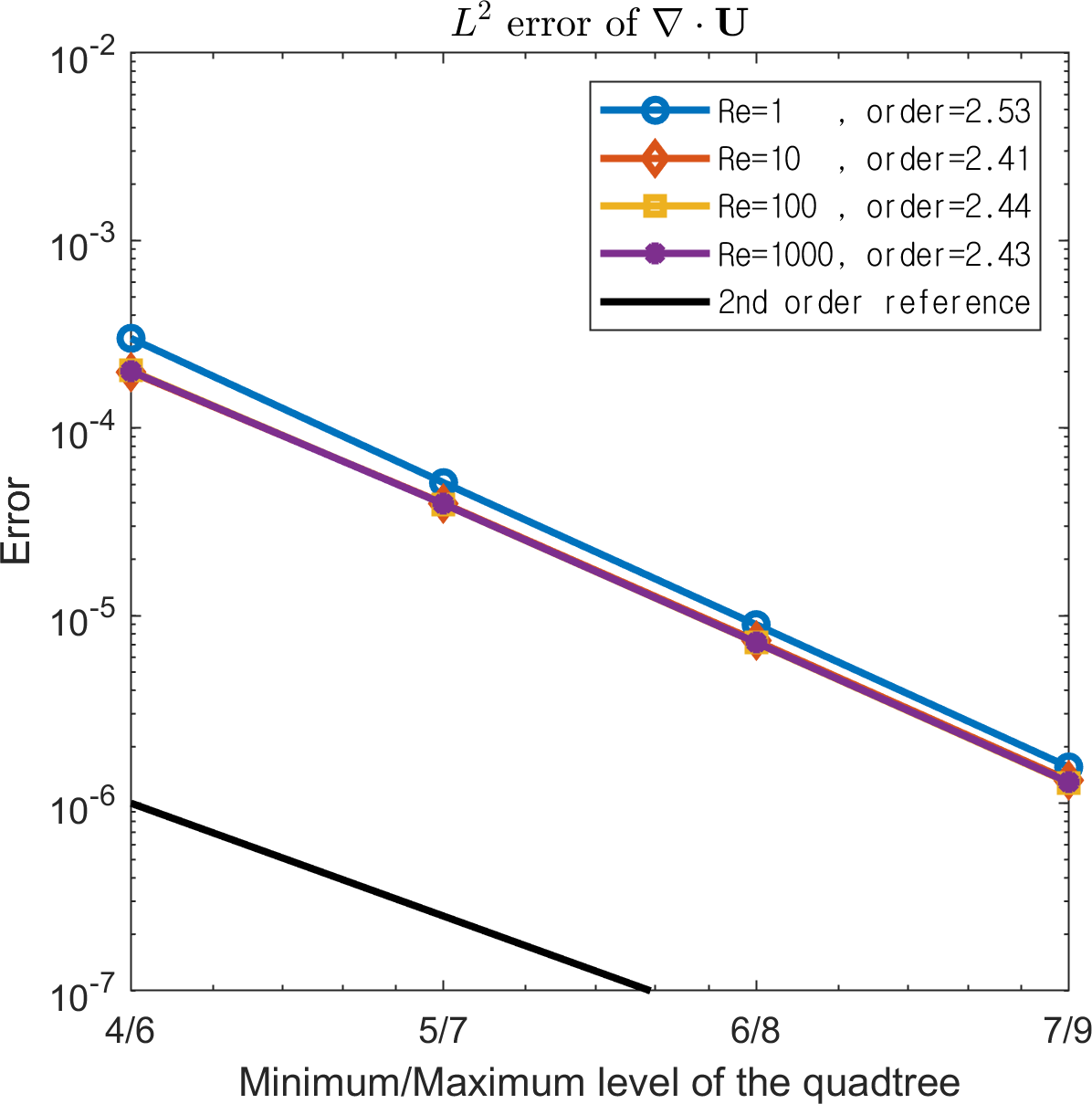}}
		}
		\caption{$L^\infty$ and $L^2$ errors of the velocity, pressure and divergence on random qaudtree grid.}\label{fig:random_quadtree_error}
	\end{figure}

	To demonstrate the robustness of the proposed method, we first describe numerical tests
	on randomly generated non-graded quadtrees of levels ranging from $4/6$ to $7/9$ for a wide range of Reynolds numbers between $1$ and $1000$. For each Reynolds number, a test was conducted in four randomly generated grids with $\Delta t= \Delta x_s$. Figure \ref{fig:random_grid} shows an example of the randomly generated quadtree of level $4/6$ and a refined grid of level $5/7$. The average numerical errors on four randomly generated grids and the best-fit order of convergence are presented in Figure \ref{fig:random_quadtree_error}. The numerical results indicate that the proposed method is second-order accurate for the fluid velocity regardless of the grid formulation and Reynolds number. Second-order convergence of the pressure of the velocity in $L^\infty$ norms was not observed at relatively low Reynolds number; however, second-order convergence in $L^2$ norms was observed.

	\begin{figure}

	\end{figure}

	\subsubsection{Irregular domains}\label{subsubsec:irregulardomain}
	
	To demonstrate the second-order accuracy in various irregular domains,
	we conduct a numerical experiment with the same exact solutions as \eqref{eq:single_vortex} with $Re=1000$ on the four domains $\Omega= \{(x,y)\in \mathbb{R}^2| \phi(x,y)\leq 0\}$ in $\left[-\frac{\pi}{2},\frac{\pi}{2}\right]^2$:\\
	
\begin{figure}
		\centering
		\mbox{
			\subfigure[circle]{\includegraphics[width=0.23\textwidth]{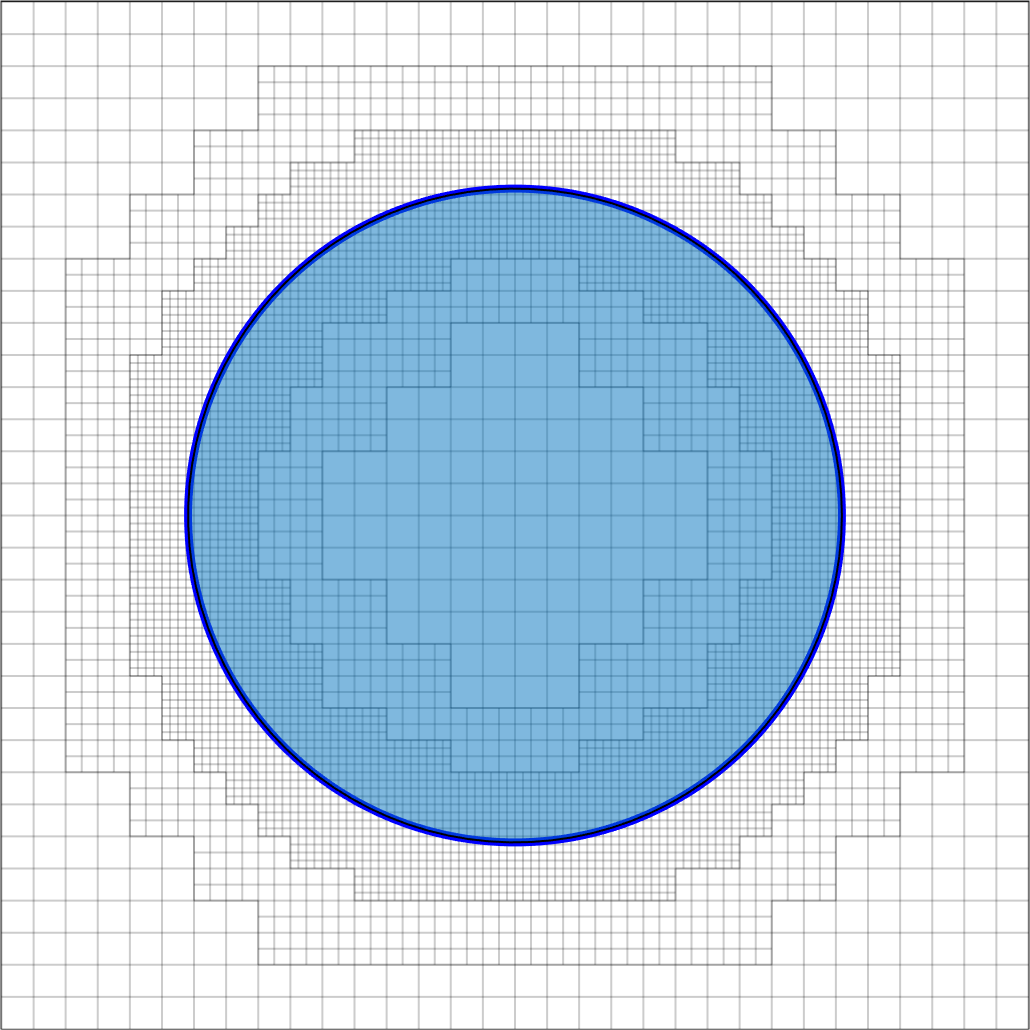}}
			\subfigure[ellipse]{\includegraphics[width=0.23\textwidth]{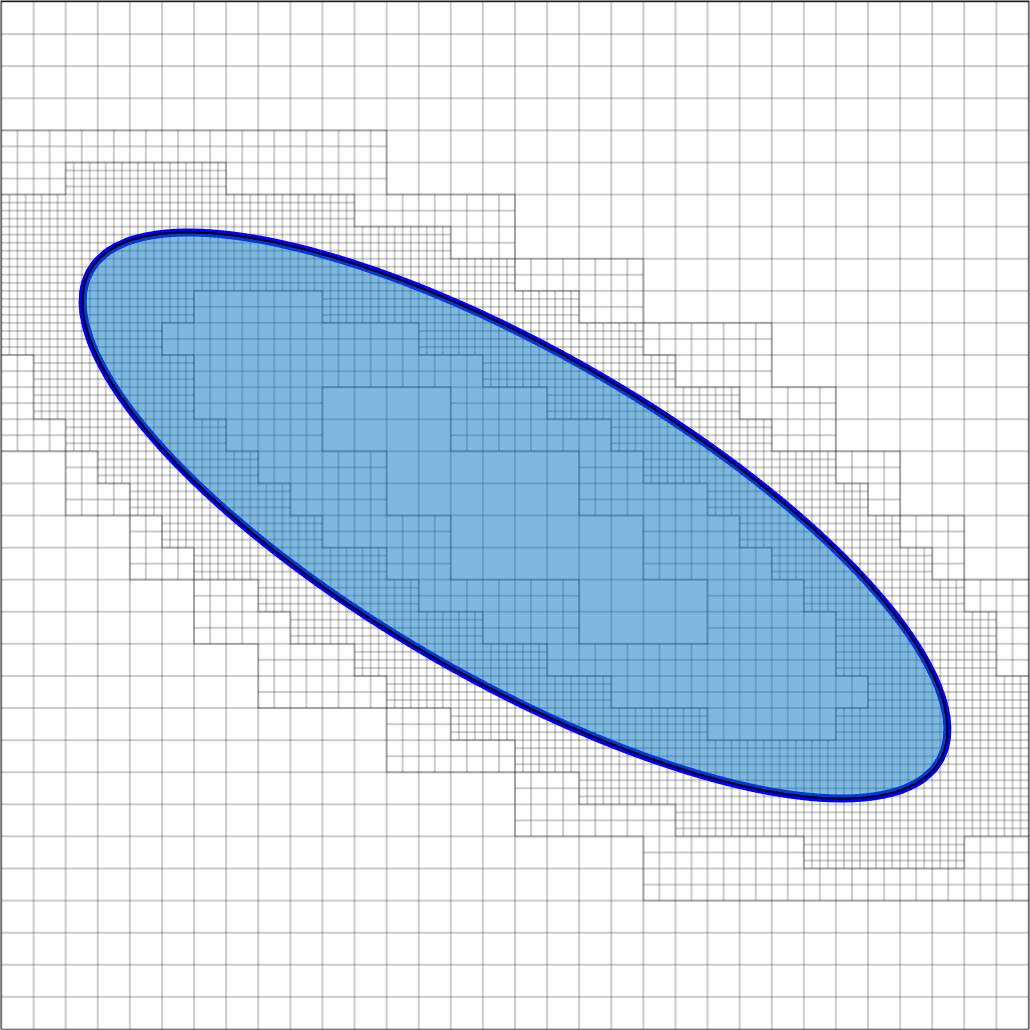}}
		}
		\centering
		\mbox{
			\subfigure[flower]{\includegraphics[width=0.23\textwidth]{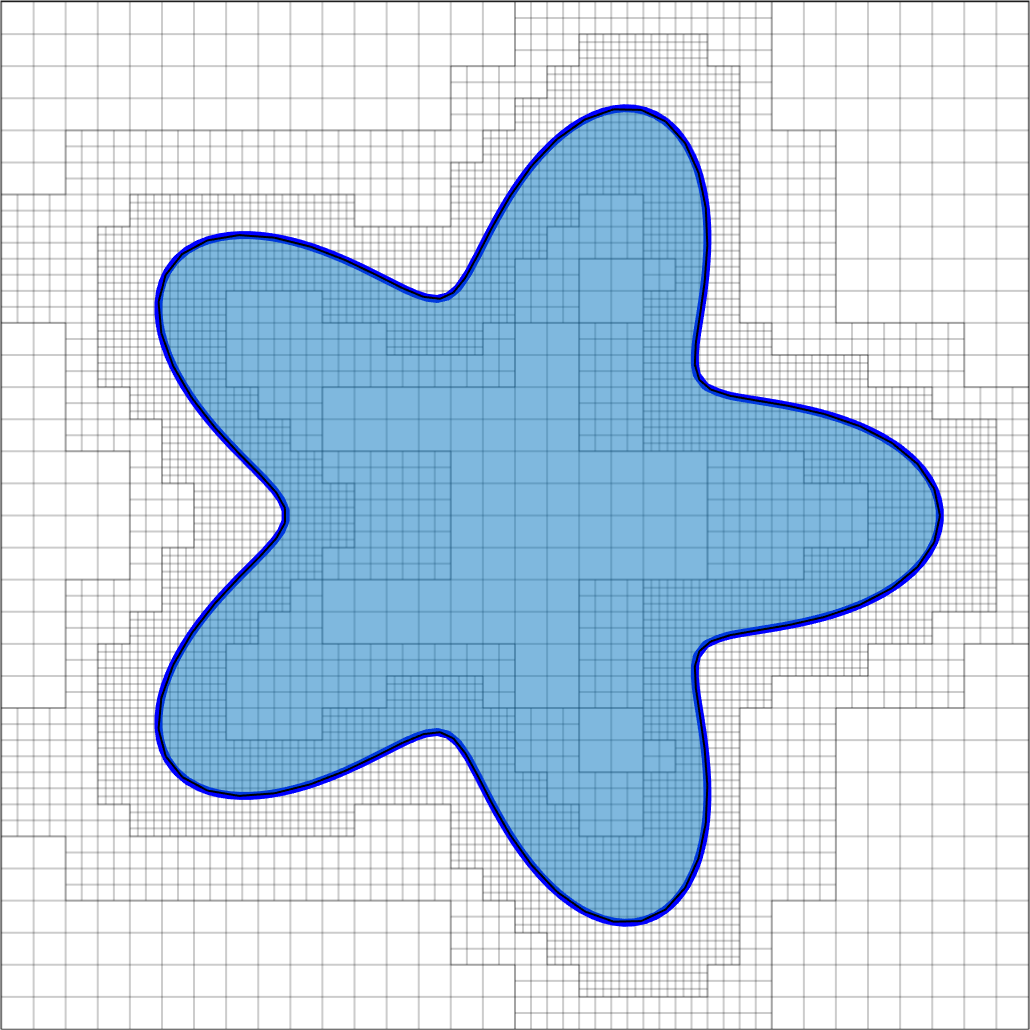}}
			\subfigure[moving flower]{\includegraphics[width=0.23\textwidth]{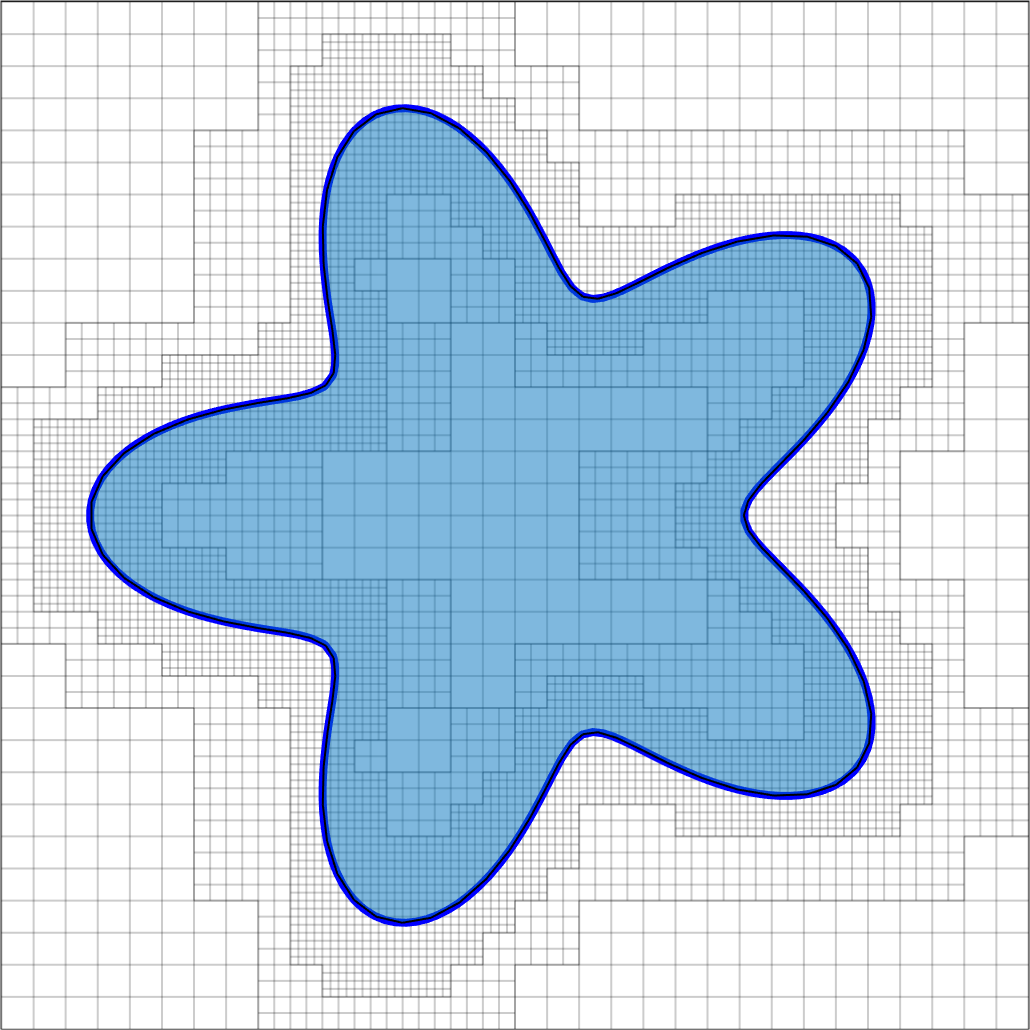}}
		}
		\caption{Various test domains of example \ref{subsubsec:irregulardomain}}\label{fig:TestDomains}
	\end{figure}

\begin{figure}
		\centering{}
		\mbox{
			\subfigure[]{\includegraphics[width=0.3\textwidth]{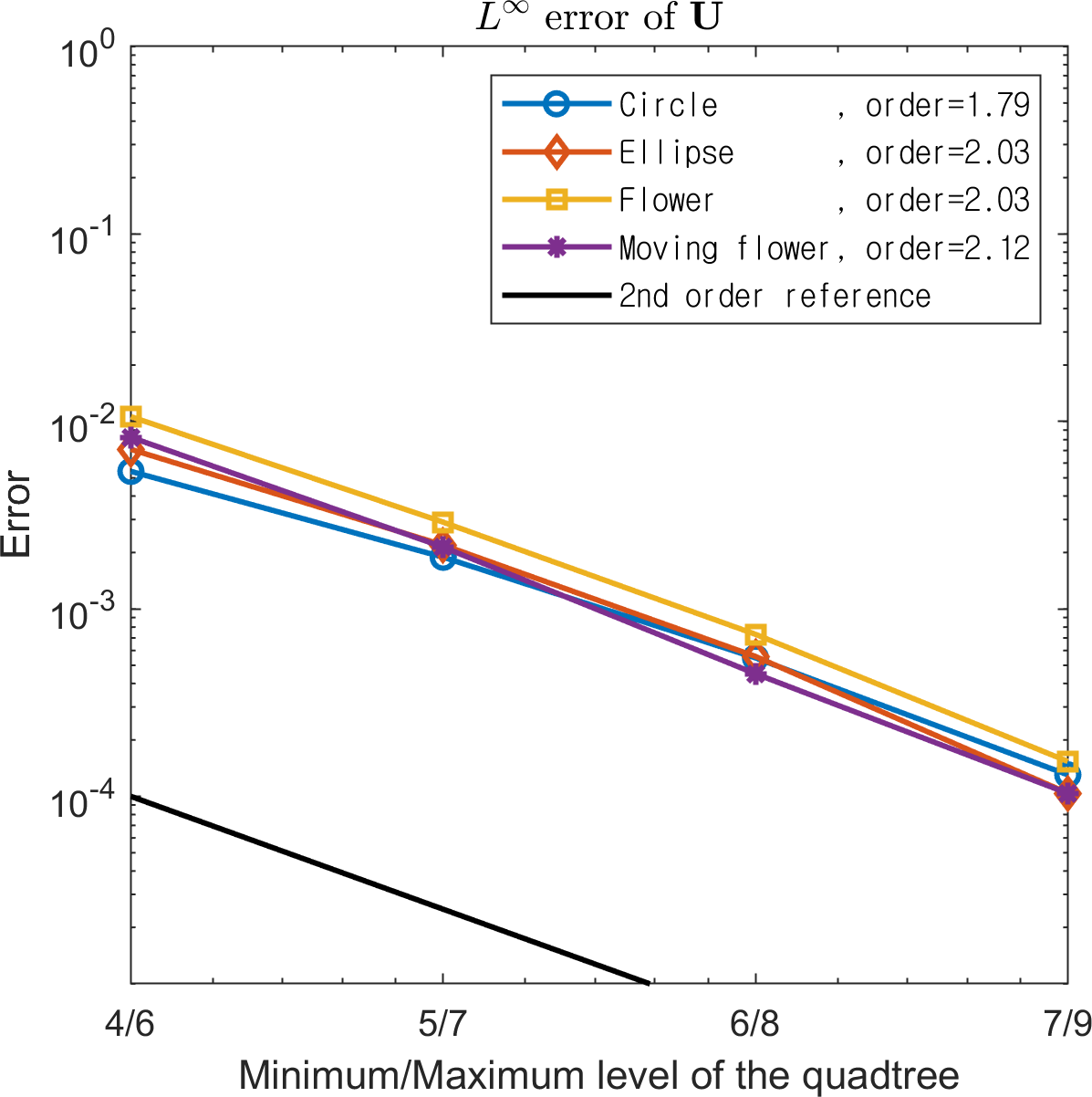}}	
			\subfigure[]{\includegraphics[width=0.3\textwidth]{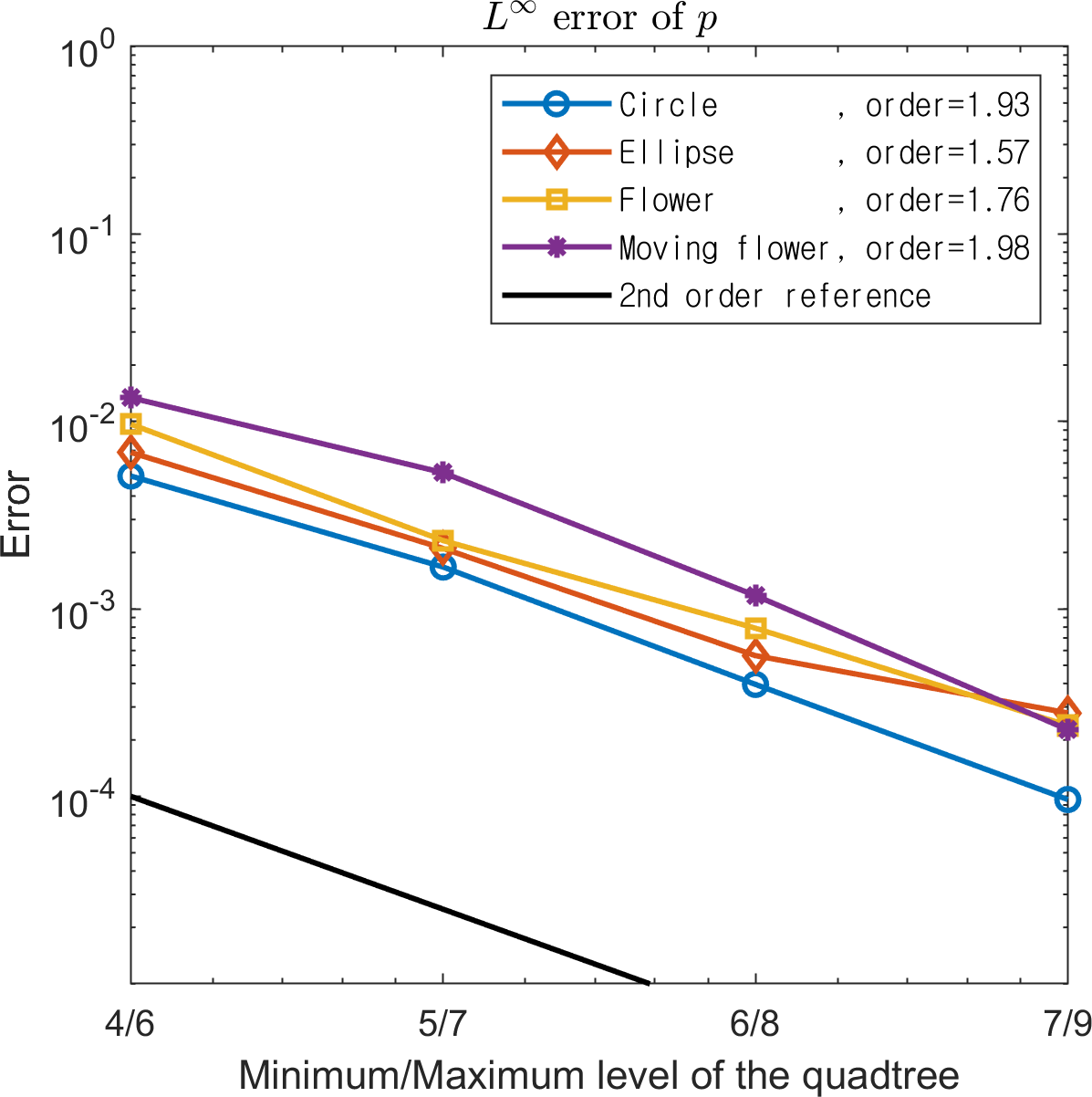}}
			\subfigure[]{\includegraphics[width=0.3\textwidth]{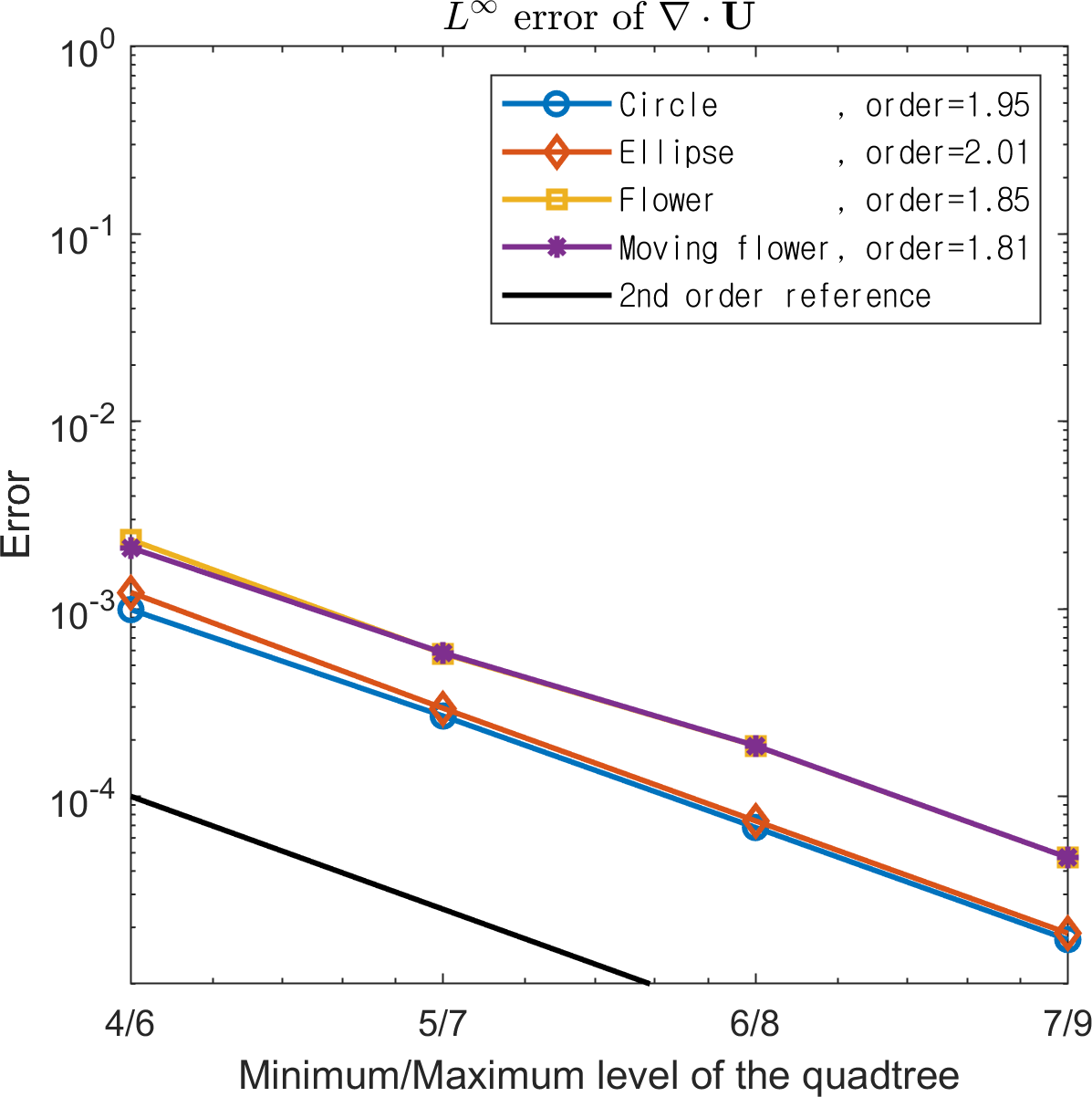}} 
		}\\
		\centering{}
		\mbox{
			\subfigure[]{\includegraphics[width=0.3\textwidth]{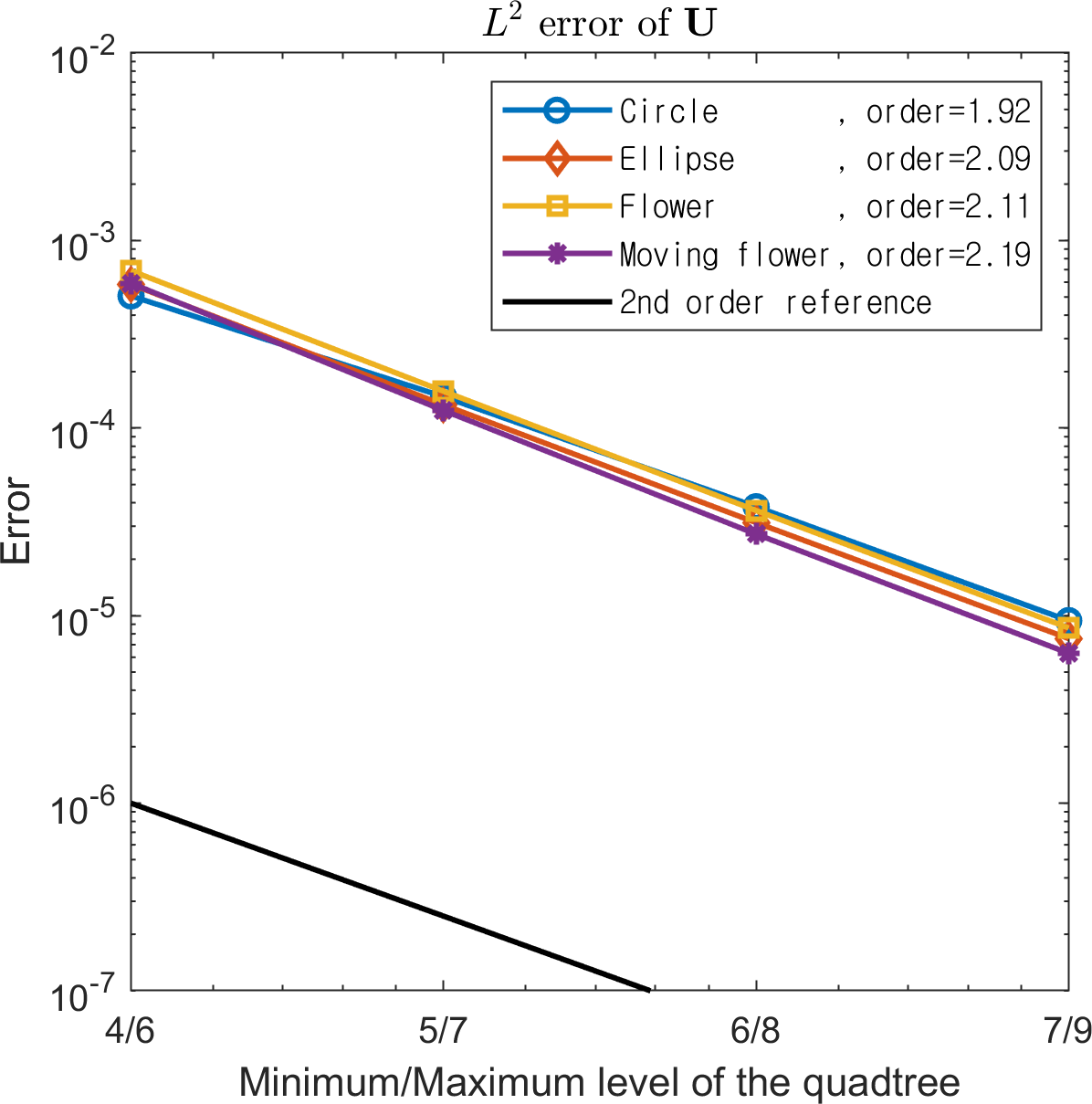}}
			\subfigure[]{\includegraphics[width=0.3\textwidth]{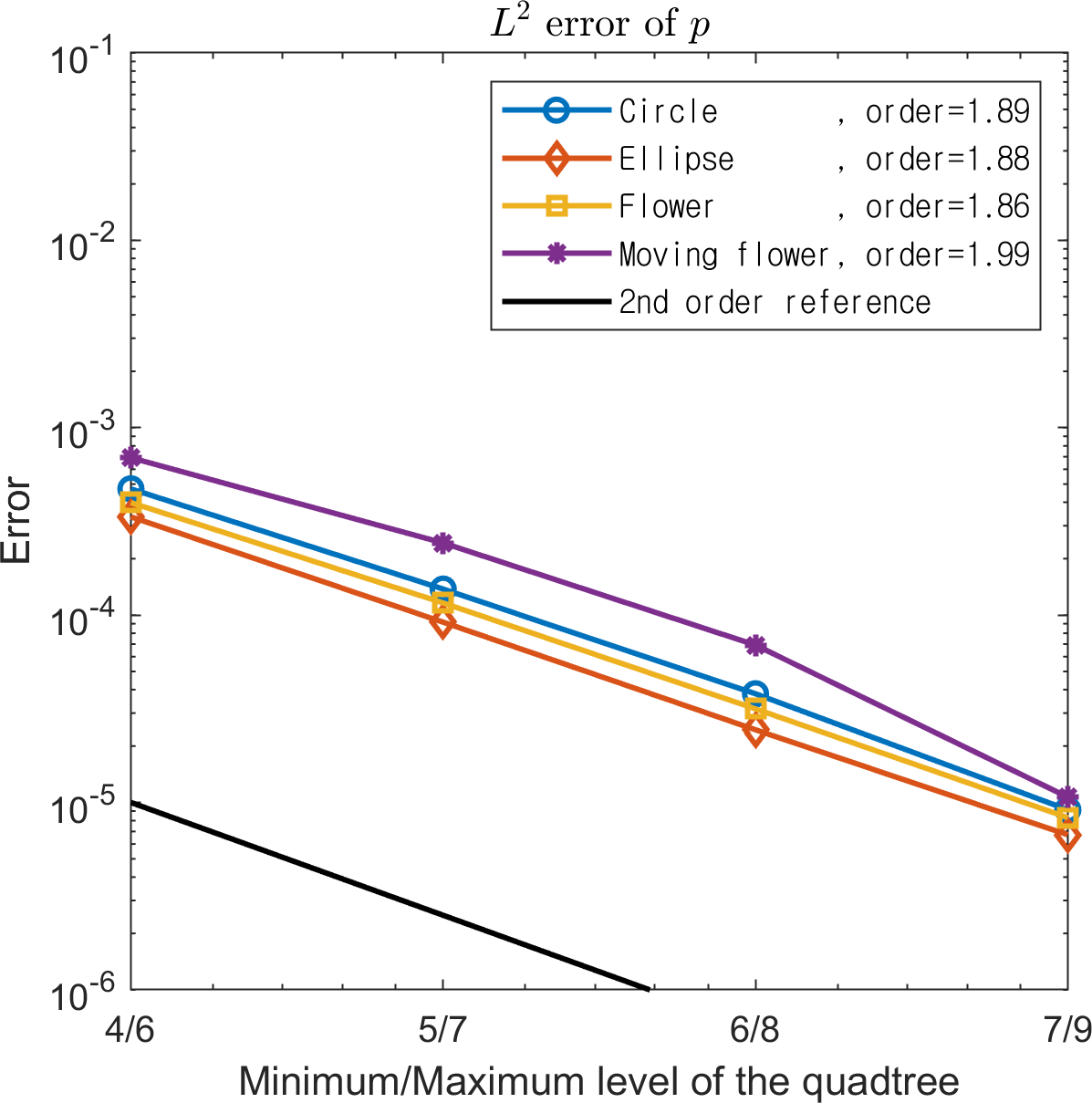}}
			\subfigure[]{\includegraphics[width=0.3\textwidth]{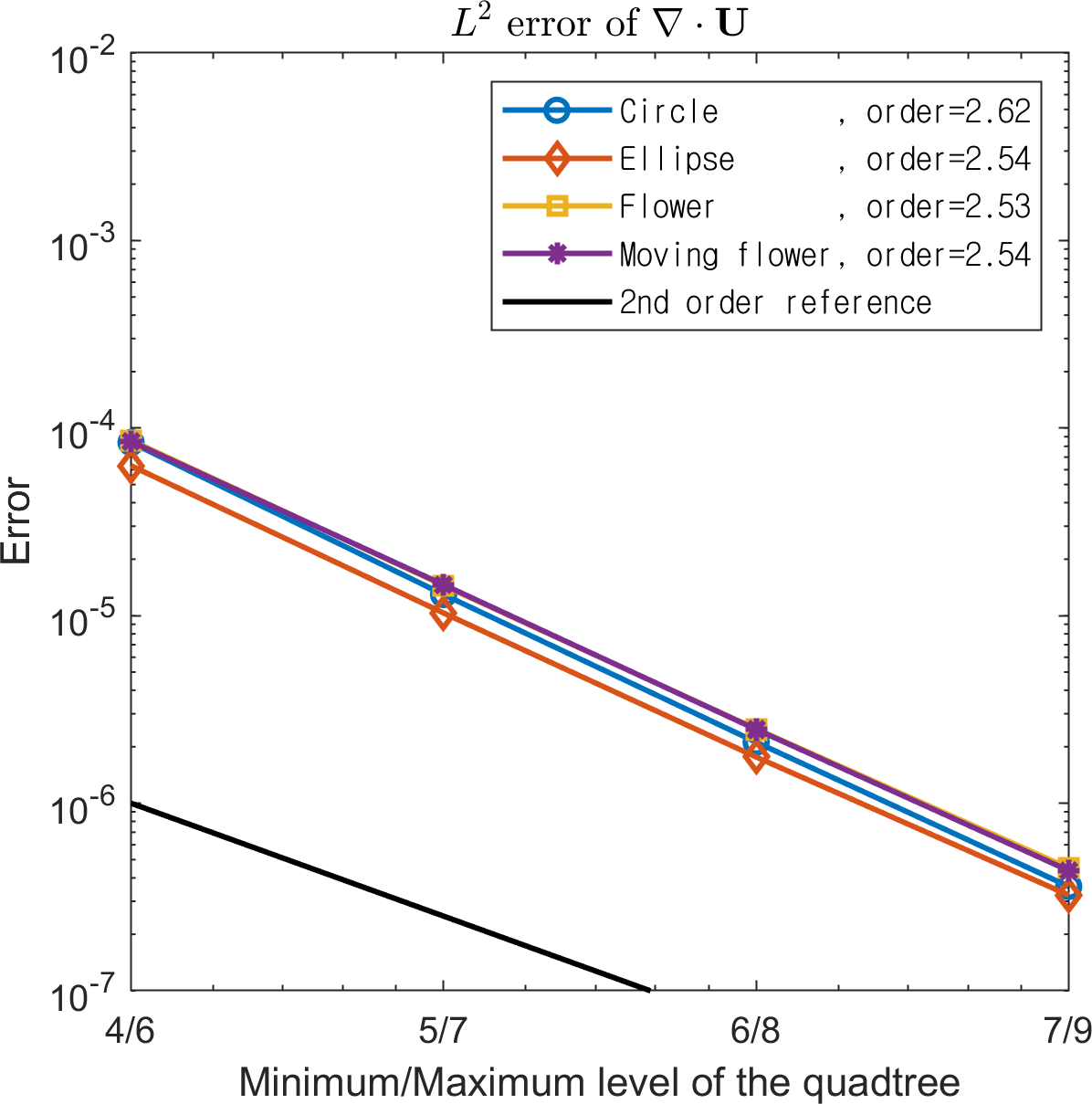}}
		}
		\caption{$L^\infty$ and $L^2$ errors of the velocity, pressure and divergence on various irregular domains.}\label{fig:irregular_domain_error}
	\end{figure}

	\begingroup
	\centering
	\begin{tabular}{lll}
		\it Circle       \rm &:& $\phi(x,y)=\sqrt{x^{2}+y^{2}}-1$      .                                                                                                                             \\
		\it Ellipse      \rm &:& $\phi(x,y)=\sqrt {\frac{\left(x \cos\frac{\pi}{6}-y\sin\frac{\pi}{6} \right)^2}{ 1.5^2}+\frac{\left( x\sin\frac{\pi}{6}+ y\cos\frac{\pi}{6}\right)^2}{0.5^2}} - 1$ .\\
		\it Flower       \rm &:& $\phi(r,\theta,t)=r-1-0.3 \cos\left(5\theta\right)$    .                                                                                                            \\ 
		\it Moving Flower\rm &:& $\phi(r,\theta,t)=r-1-0.3 \cos\left(5(\theta-t)\right)$ .                                                                                                          
	\end{tabular}
	\\
	\endgroup
	$ \quad$\\
	\noindent
	Note that the domains \it Circle, Ellipse, \rm and \it Flower \rm do not change with respect to time, whereas the other domain, \it Moving Flower\rm, rotates counterclockwise. The domains and the generated quadtree grid of level $5/7$ at $t=\pi$ are presented in Figure \ref{fig:TestDomains}. Time step size is set to be $\Delta t= \Delta x_s$. Second-order convergence of the velocity, pressure, and divergence in both the $L^\infty$ and $L^2$ norms are observed in 
	Figure \ref{fig:irregular_domain_error}. The results indicate that the proposed method can handle irregular domains and moving domains.
	
	In Table \ref{tab:pressure_grad}, we present the results obtained by discretizing $\nabla p$ in the finite volume approach introduced in Section \ref{subsubsec:pressure_fv}, which are compared with those obtained by the finite difference approach.
	The approximation errors of the velocity measured in the infinity norm and the order of accuracy at Reynolds numbers of $Re=1,1000$ on the circular domain are reported. The finite volume approach delivers second-order accuracy for $Re=1$; however, it is reduced to first-order convergence at the higher Reynolds number, $Re=1000$. Otherwise, the finite difference method shows satisfactory second-order convergence independent of Reynolds number, and the magnitude of the errors is lower at all resolutions. 	
	
\begin{figure}
		\centering{}
		\mbox{
			\subfigure[][]{	 			\includegraphics[width=0.45\textwidth]{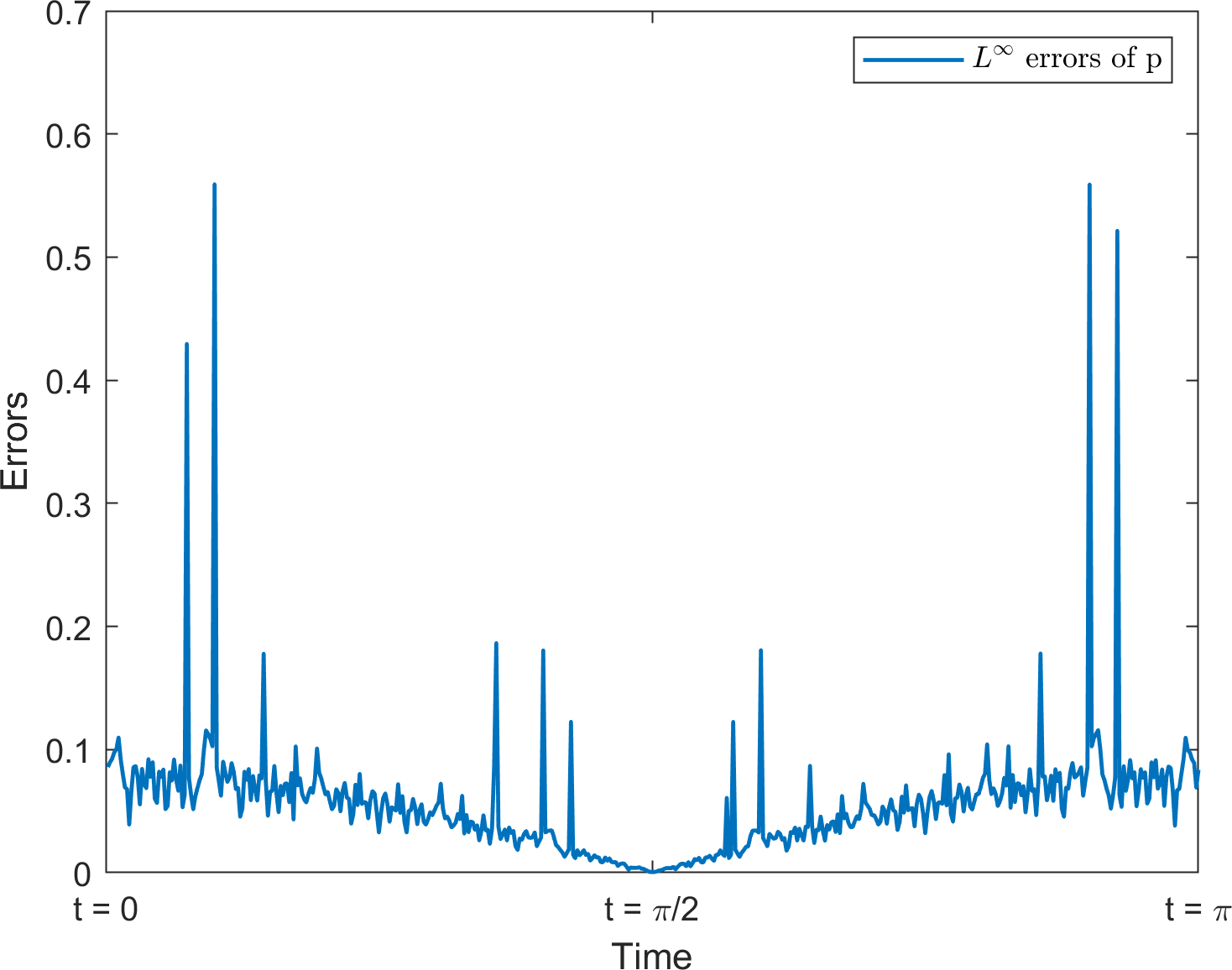} } $\ $
			\subfigure[][]{ 				\includegraphics[width=0.45\textwidth]{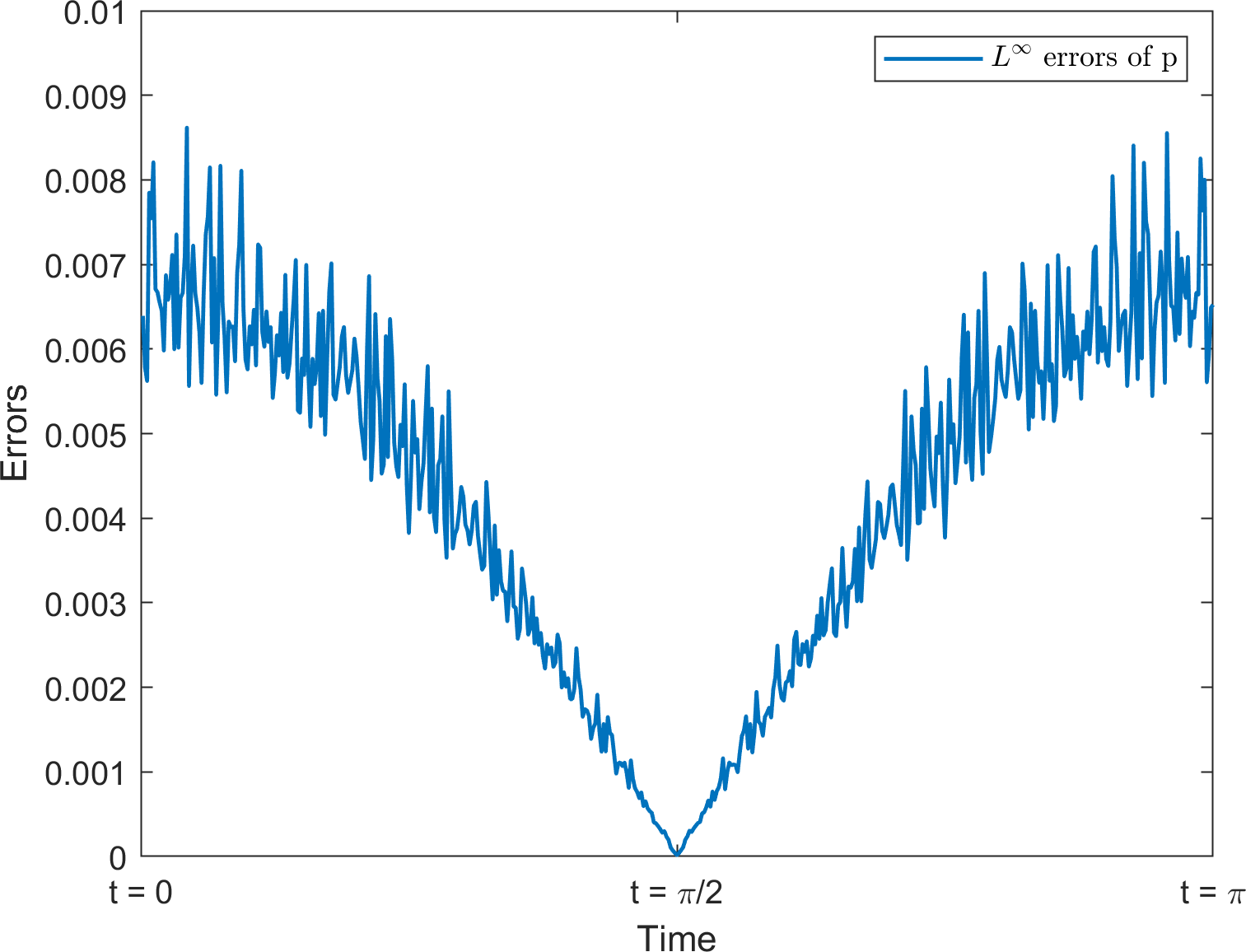} }$\ $
		}
		
		\caption{$L^\infty$ errors of pressure when (a) stabilizer is not applied and (b) stabilizer is applied for example \ref{subsubsec:irregulardomain}.}\label{fig:stabilizer_comparison}
	\end{figure}

	We demonstrate the role of pressure stabilization by comparing the results with those of an unstabilized system, that is, $\epsilon=0$. 
	Table \ref{tab:stab} reports the order of convergence of the velocity and pressure on the \textit{Moving Flower} domain with Reynolds number $Re=1$. Although second-order accuracy of the velocity was achieved in both cases, the convergence rate of the pressure was degraded in the $L^{2}$ norm when the stabilizer was not applied. In addition, Figure \ref{fig:stabilizer_comparison} shows the change in the $L^\infty$ errors of the pressure on a quadtree grid of level $6/8$ with respect to time. Although both methods show oscillations in the error, the magnitudes of the oscillations and errors are smaller when the stabilizer is used. 
	
	\begin{figure}
		\centering{}
		\mbox{
			\subfigure[]{\includegraphics[width=0.3\textwidth]{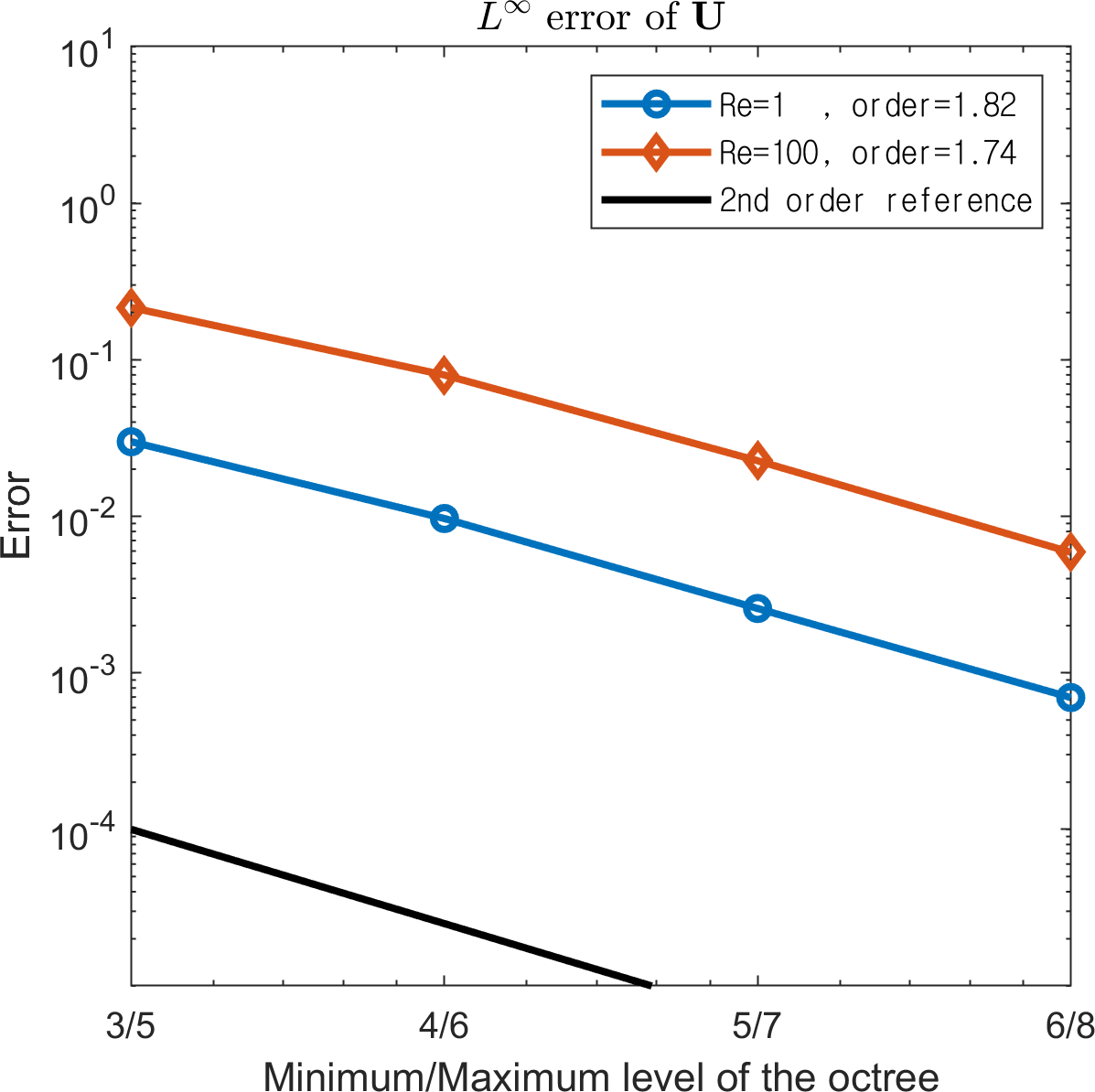}}	
			\subfigure[]{\includegraphics[width=0.3\textwidth]{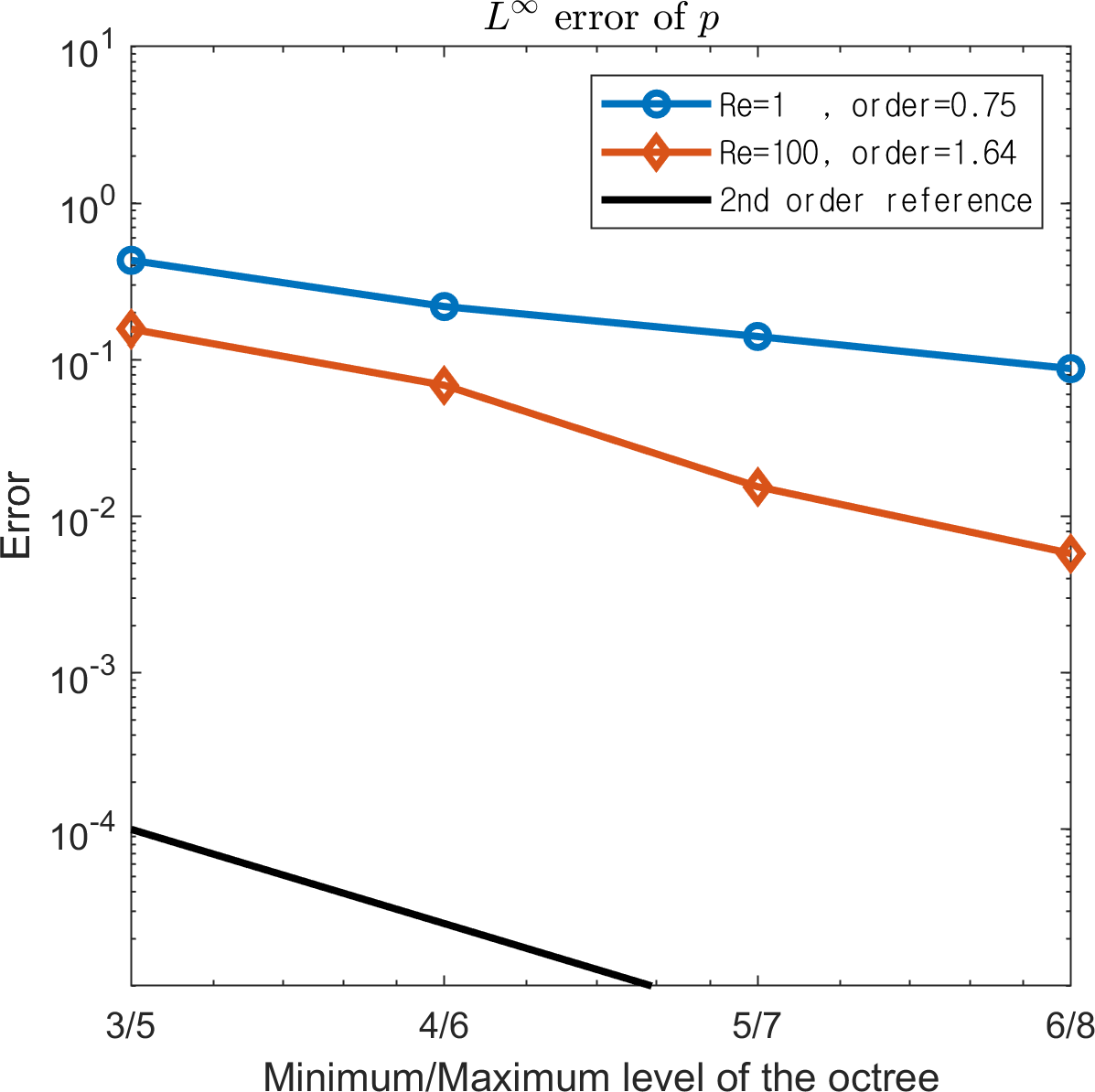}}
			\subfigure[]{\includegraphics[width=0.3\textwidth]{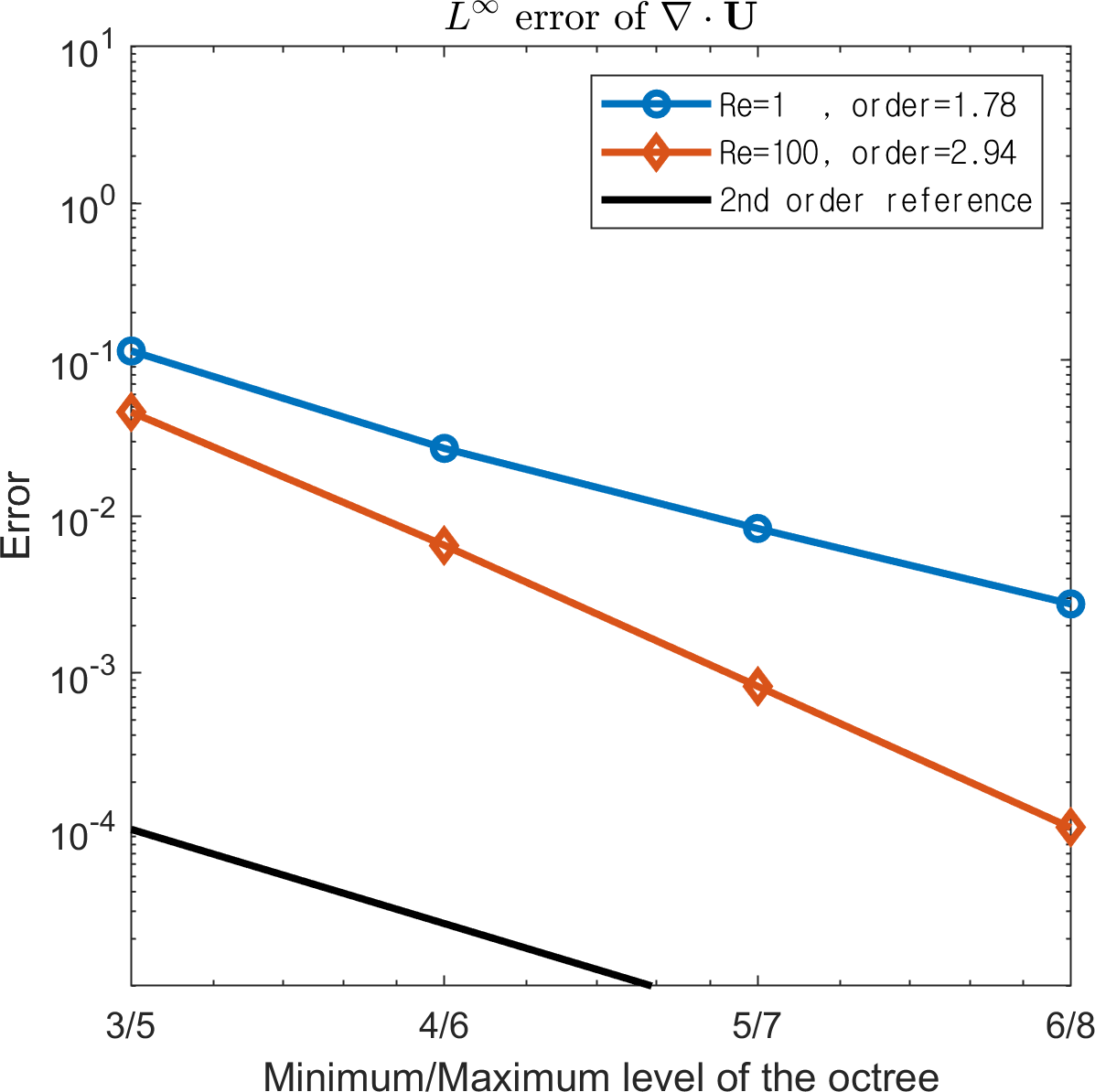}} 
		}\\
		\centering{}
		\mbox{
			\subfigure[]{\includegraphics[width=0.3\textwidth]{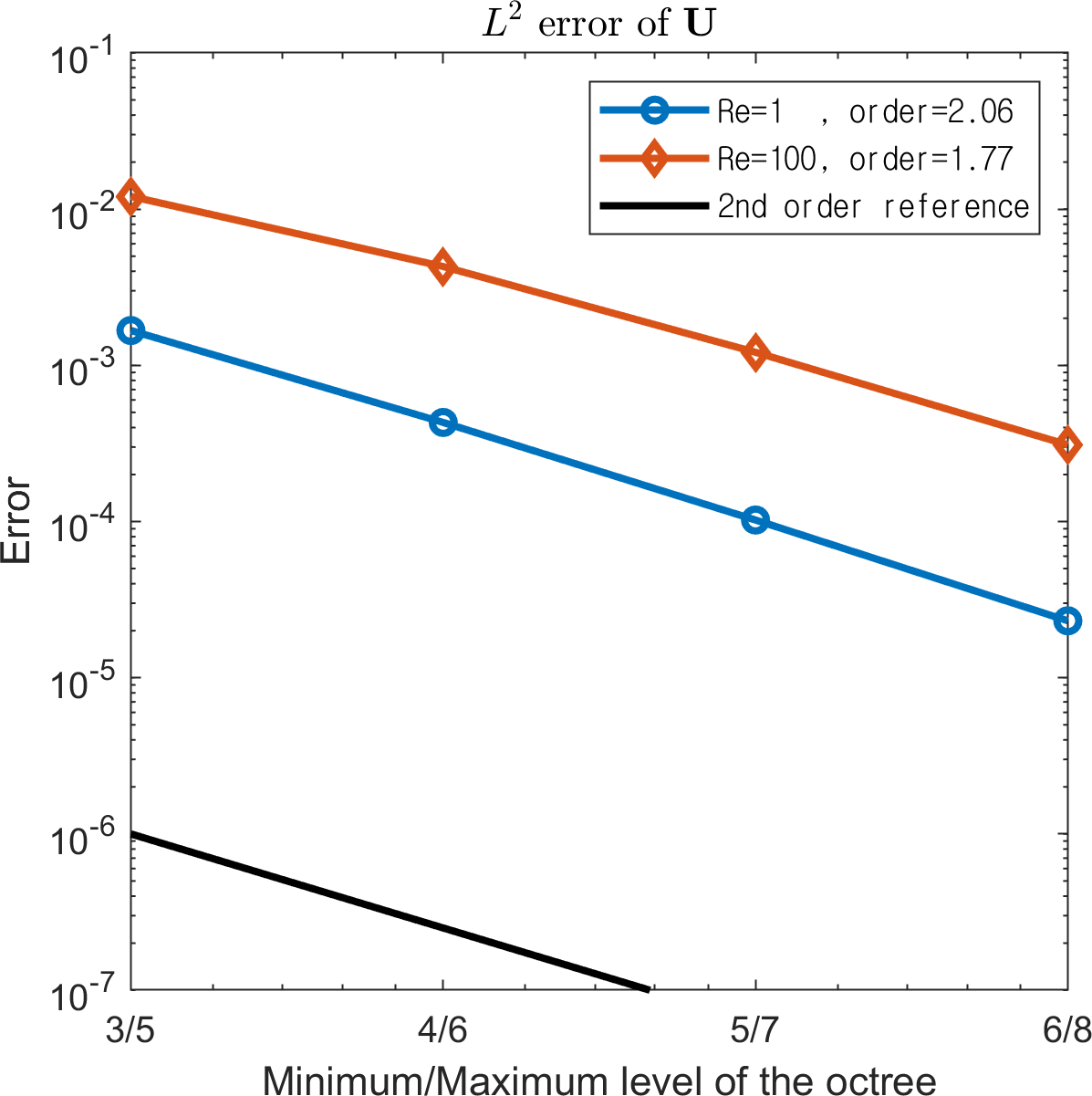}}
			\subfigure[]{\includegraphics[width=0.3\textwidth]{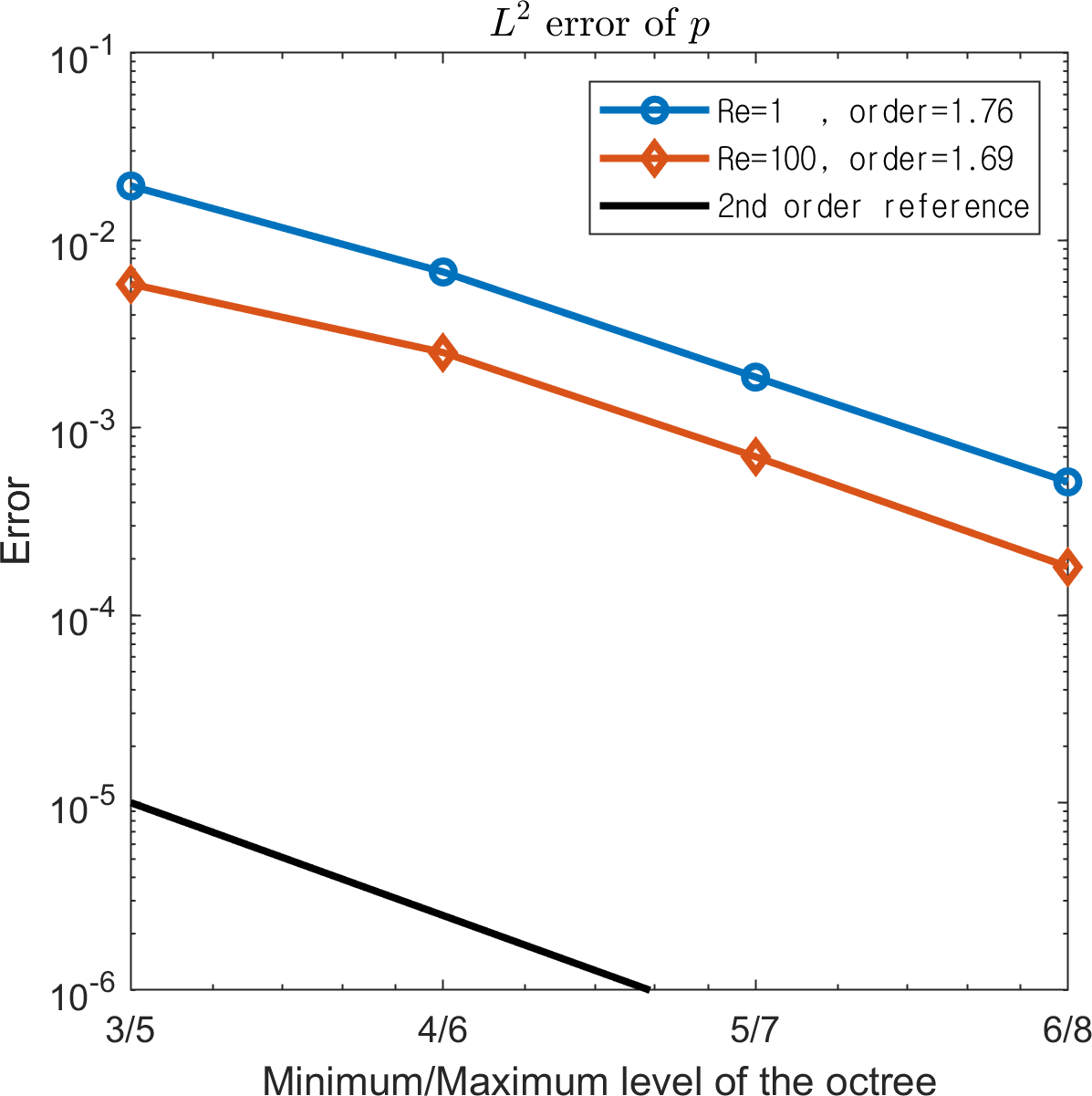}}
			\subfigure[]{\includegraphics[width=0.3\textwidth]{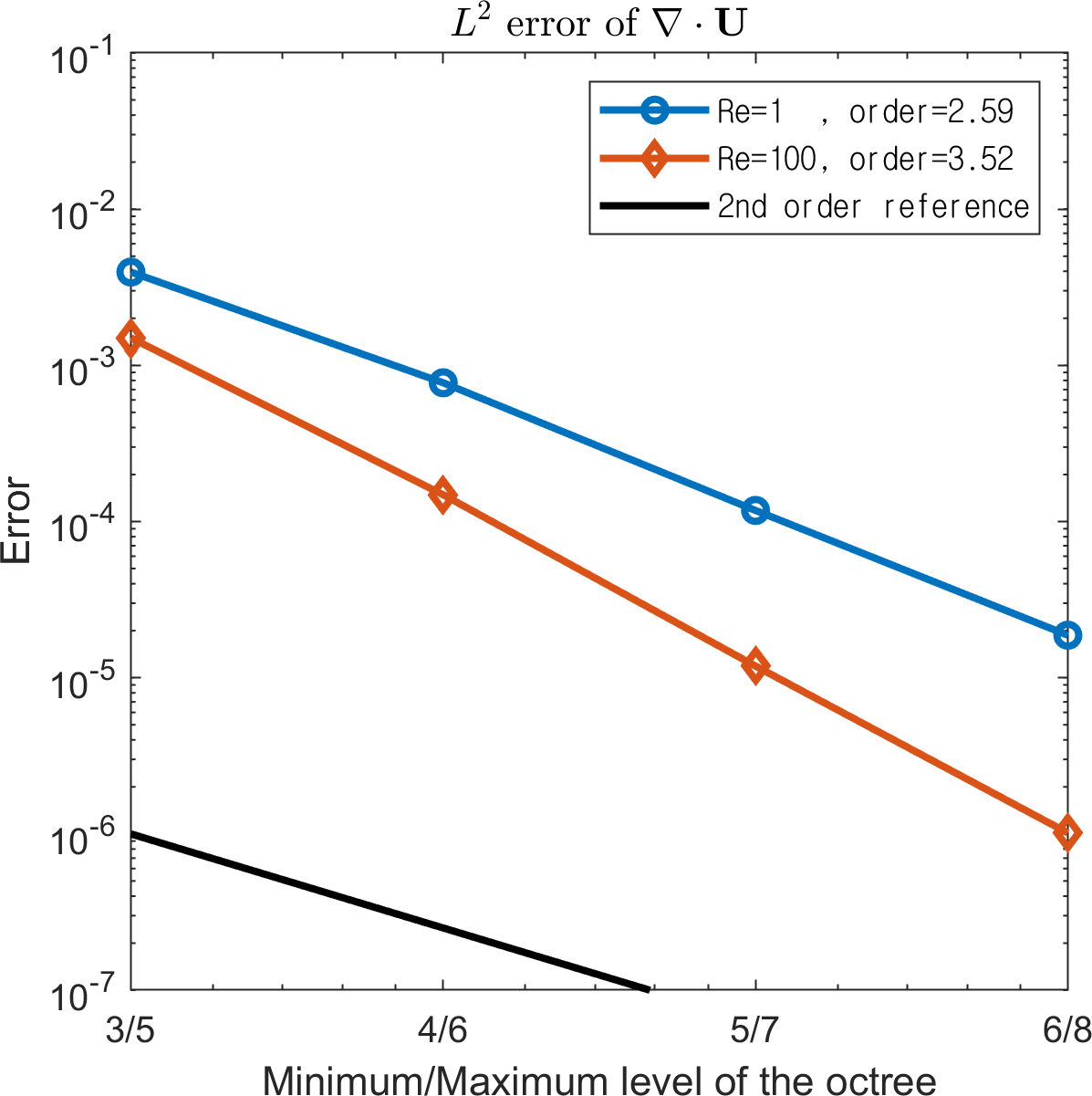}}
		}
		\caption{$L^\infty$ and $L^2$ errors of the velocity, pressure and divergence on random octree.}\label{fig:octree_error}
	\end{figure}

	\begin{table}[]
		\centering{}
		\begin{tabular}{llllllllll}
			& \multicolumn{4}{c}{ $Re=1$}               &        & \multicolumn{4}{c}{$Re=1000$}  \\ \cline{2-5} \cline{7-10} 
			& FV approach& order & FD approach& order              &        & FV approach     & order & FD approach     & order \\ \cline{2-5} \cline{7-10} 
			4/6 & 1.23.E-03&		    &6.81.E-04&	    &			&     2.12.E-02	&	 & 	5.00.E-03&			\\
			5/7 & 3.09.E-04&	2.00	&2.26.E-04&	1.59&		&9.66.E-03	&1.14&	1.89.E-03	&1.52	\\
			6/8 & 7.83.E-05&	1.98	&6.40.E-05&	1.82&	&5.94.E-03	&0.70&	5.48.E-04	&1.78\\
			7/9 & 2.13.E-05&	1.88	&1.71.E-05&	1.91&		&2.70.E-03	&1.14&	1.30.E-04	&2.07
		\end{tabular}
		\caption{Comparison of $L^\infty$ errors in velocities between finite volume and finite difference discretization for pressure gradient.}\label{tab:pressure_grad}
	\end{table}
		
	\begin{table}[]
		\centering{}
		\begin{tabular}{llllllllll}
			& \multicolumn{4}{c}{ Not stabilized}               &        & \multicolumn{4}{c}{Stabilized}  \\ \cline{2-5} \cline{7-10} 
			& $L^\infty$ errors of $\uvec$     & order & $L^2$ errors of $p$     & order              &        & $L^\infty$ errors of $\uvec$ &order & $L^2$ errors of $p$     & order \\ \cline{2-5} \cline{7-10} 
			4/6 & 7.73.E-04 &       & 2.89.E-02 &                 &           & 7.82.E-04 &       & 8.37.E-03 &       \\
			5/7 & 3.33.E-04 & 1.21  & 2.24.E-02 & 0.36             &          & 3.36.E-04 & 1.22  & 3.69.E-03 & 1.18  \\
			6/8 & 8.09.E-05 & 2.04  & 1.23.E-02 & 0.87              &         & 8.20.E-05 & 2.04  & 1.25.E-03 & 1.56  \\
			7/9 & 2.16.E-05 & 1.90  & 1.32.E-02 & -0.10 & &                   2.16.E-05 & 1.93  & 5.31.E-04 & 1.24 
		\end{tabular}
		\caption{Convergence comparison between results with/without the stabilizer for example \ref{subsubsec:irregulardomain}.}\label{tab:stab}
	\end{table}
	
	\subsubsection{Random octree}
	As a three-dimensional example, we consider the following analytical solution on the cubical domain $\Omega =\left[-\frac{\pi}{2},\frac{\pi}{2}\right]^{3}$:
	\begin{equation}\label{eq:single_vortex_3d}
		\begin{aligned}
			u(x,y,z,t) & =-2\cos x\sin y\sin z\cos t,\\
			v(x,y,z,t) & =\sin x\cos y\sin z\cos t,\\
			w(x,y,z,t) & = \sin x\sin y\cos z\cos t,\\
			p(x,y,z,t) & =\frac{1}{4}\cos^{2}t\left(2\cos 2x+\cos 2y+\cos 2z \right).
		\end{aligned}
	\end{equation}
	Simulations were conducted up to $t=\pi$ with $\Delta t = 4\Delta x_s$ on a random octree grid of levels ranging from $3/5$ to $6/8$. Two Reynolds numbers, $Re=1,100$, are considered. Figure \ref{fig:octree_error} presents the $L^\infty$ and $L^2$ errors, which are compared to those of the analytical solutions. Although the convergence rate of the pressure attained at $Re=1$ is lower than 1, we observe overall second-order convergence rates even for a random octree grid.

	\subsection{Driven cavity flows}
	
	In this section, we consider the  lid-driven cavity flow problem. Because there is no exact solution of cavity flow, we compare the simulation results to those reported in previous works. We first simulated cavity flow on a rectangular domain and compared it with the reference results of Ghia et al. \cite{ghia1982high}. The numerical simulation of driven cavity flow around a flower-shaped object introduced by Coco \cite{coco2020multigrid} follows. In this example, we adaptively refine the cell $\mathbf{C}$ when $\max_{\mathbf{x} \in Vertices(\mathbf{C})}\frac{ \| \nabla \uvec (\mathbf{X}) \|_2}{\|\uvec\|_\infty } \Delta x > \varepsilon$ for a tolerance $\varepsilon$ and cell size $\Delta x$. As a two-dimensional example, we set the tolerance to $\varepsilon=0.05$.
	
	\subsubsection{Lid-driven cavity on a rectangular domain}\label{subsubsec:cavitysquare}
	We solve the lid-driven cavity flow problem on the computational domain $\Omega= \left[0,1\right]^2$ with the boundary condition 
	\[ \uvec=\begin{cases}
		(1,0) \text{ if $y=1$, $0<x<1$}\\
		(0,0) \text{ if $y=0$, $x=0,1$}\\
	\end{cases} . \]
	Numerical simulations were conducted for two Reynolds numbers, $Re=100$ and $1000$, on quadtree grids of level $5/7$ and $6/8$, respectively.

\begin{figure}
		\centering{}
		\mbox{
			\subfigure[][$Re=100$]{	 			\includegraphics[width=0.45\textwidth]{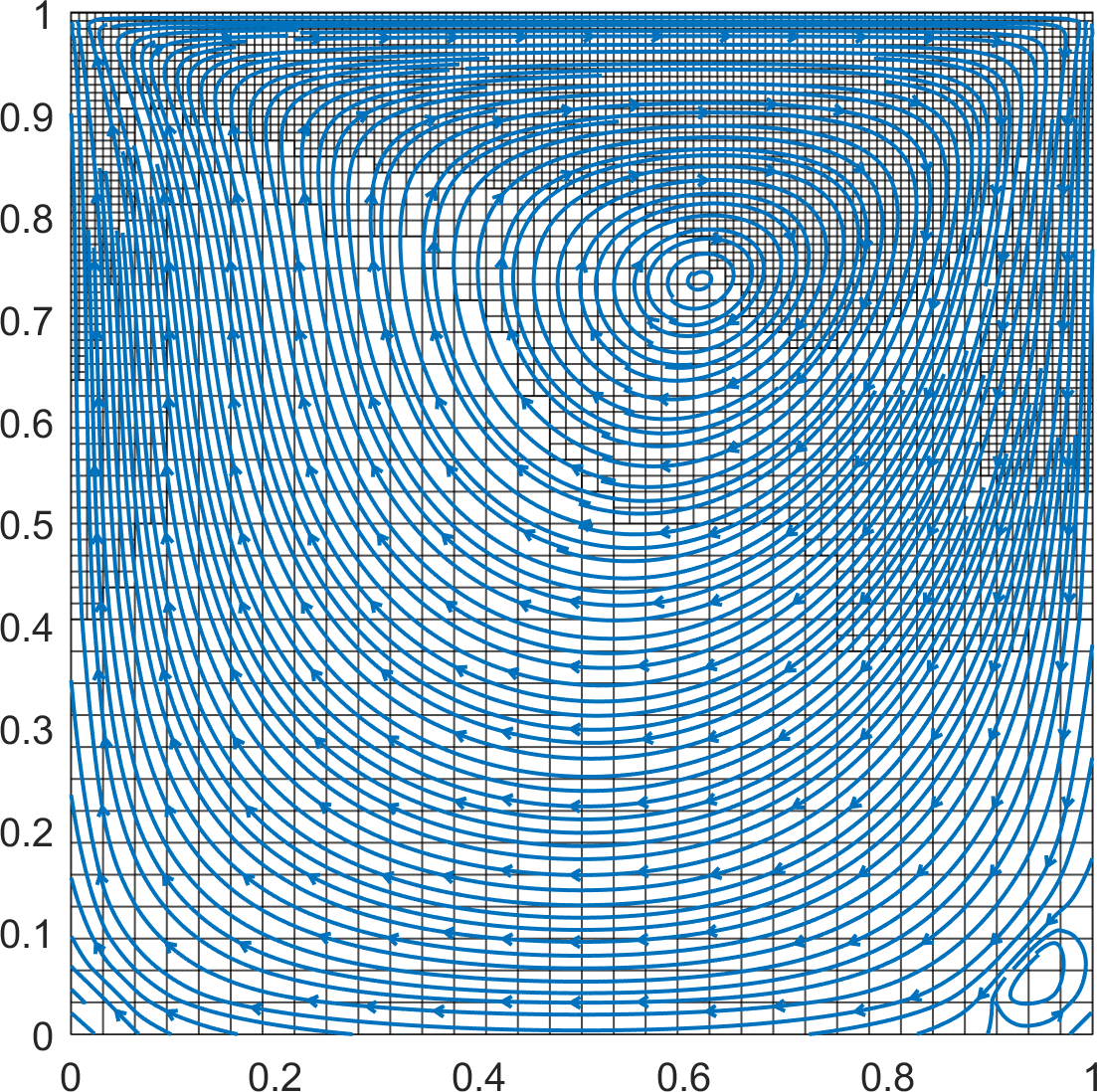} } $\ $
			\subfigure[][$Re=1000$]{ 				\includegraphics[width=0.45\textwidth]{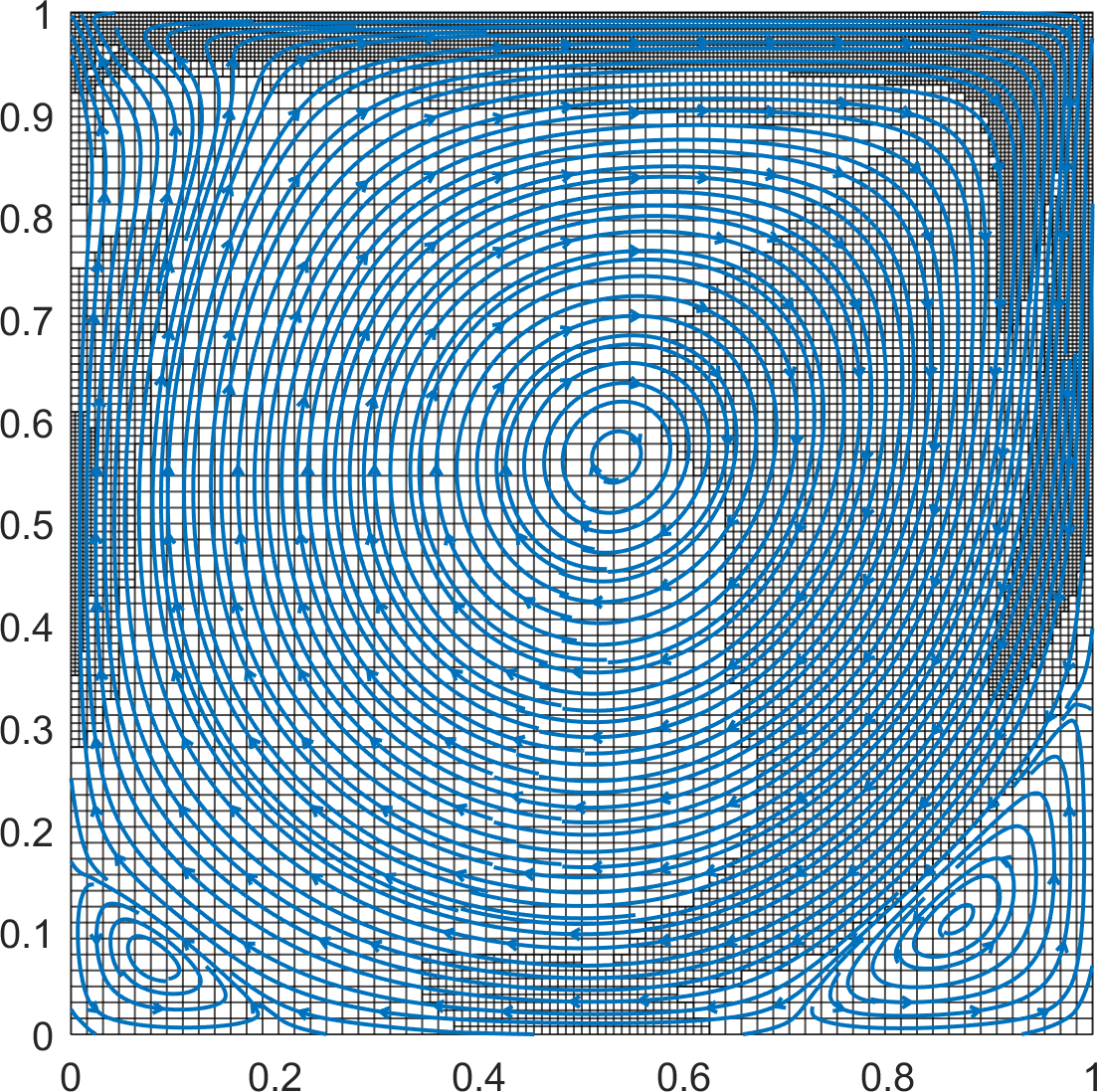} }$\ $
		}
		\caption{Streamlines and quadtree grid at $t=40$ for lid driven cavity flow.}\label{fig:cavity_grid}
	\end{figure}

\begin{figure}
		\centering{}
		\subfigure[]{\includegraphics[width=0.45\textwidth]{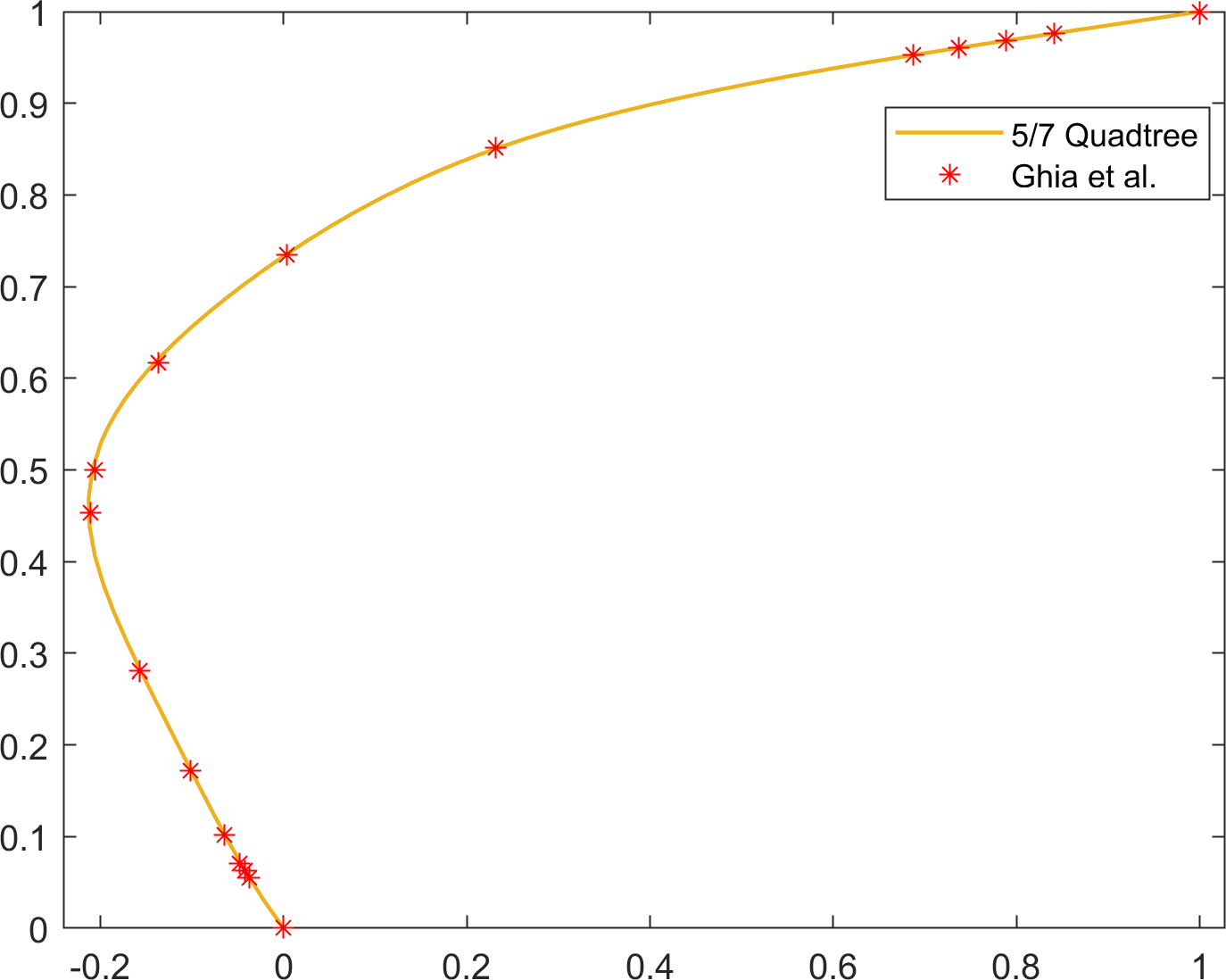}}
		\subfigure[]{\includegraphics[width=0.45\textwidth]{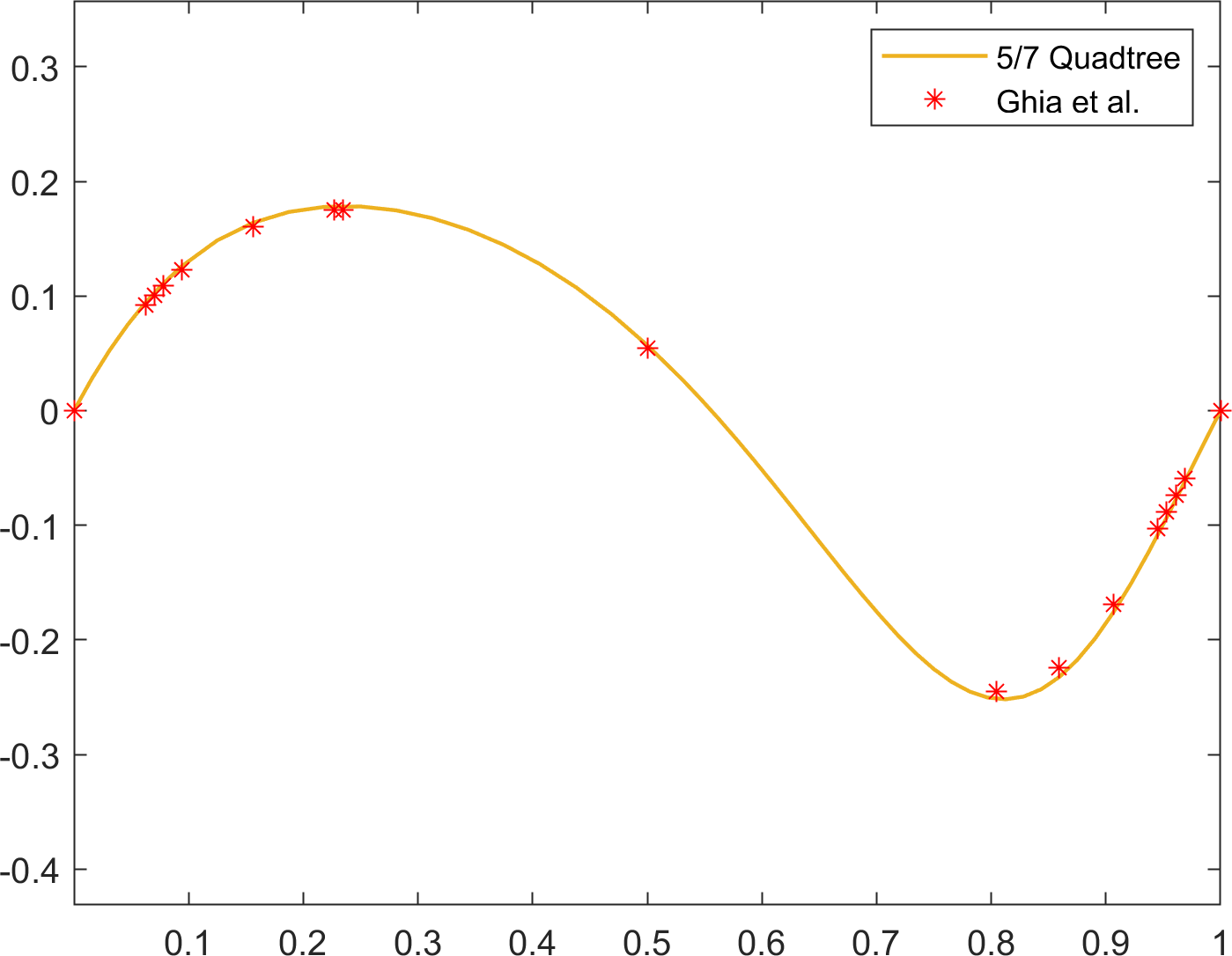}}\\
		\subfigure[]{\includegraphics[width=0.45\textwidth]{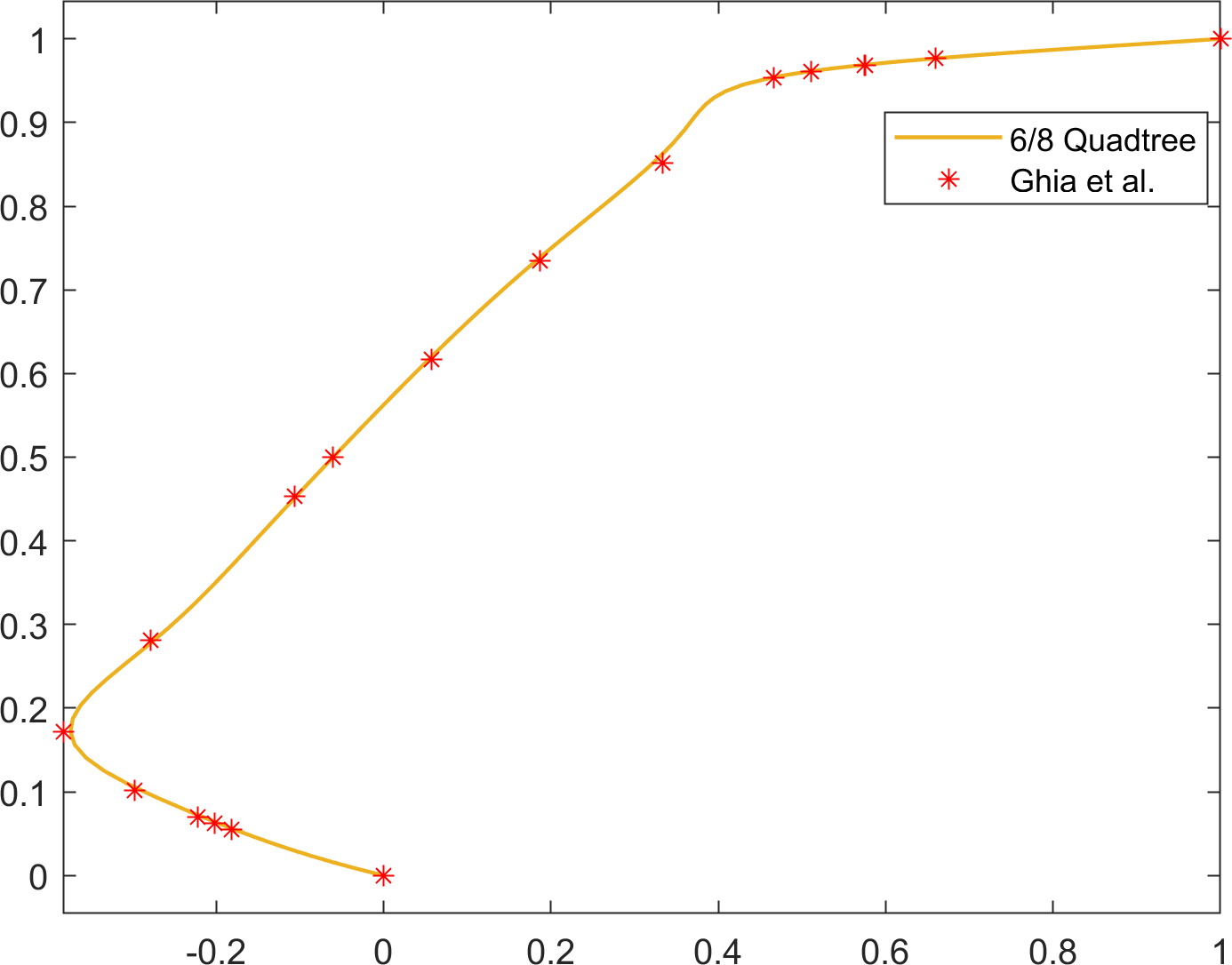}}
		\subfigure[]{\includegraphics[width=0.45\textwidth]{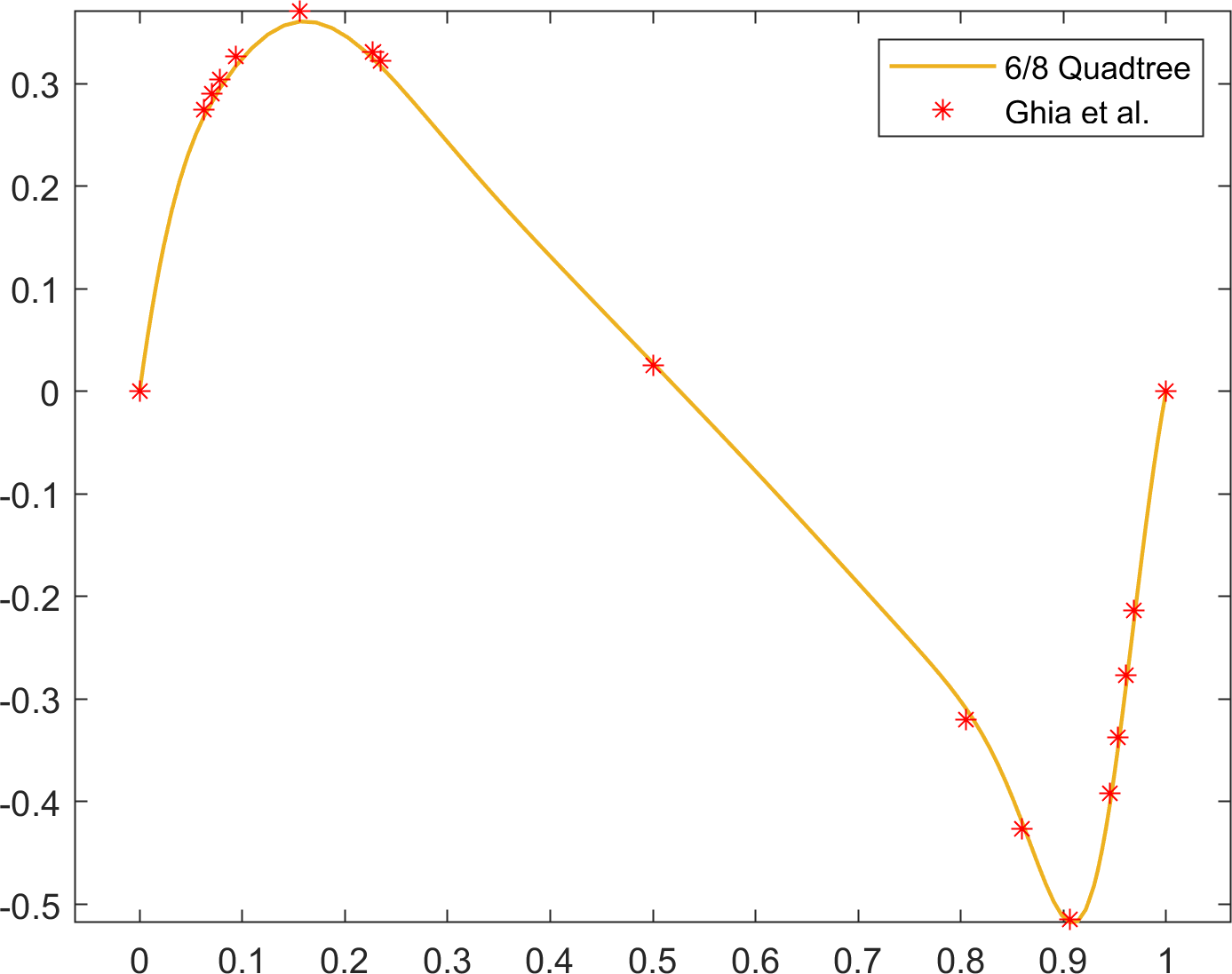}}	
		\caption{Comparison between velocity field in the driven cavity from example \ref{subsubsec:cavitysquare} and the reference results of Ghia et al. \cite{ghia1982high} of $(a),(b)$ $Re=100$, and $(c),(d)$ $Re=1000$.} \label{fig:ghia100}
	\end{figure}

	In both cases, the simulations were run until $t=40$. Streamlines and the quadtree grids are shown in Figure \ref{fig:cavity_grid}. In Figure \ref{fig:ghia100}, the results are quantitatively compared with the reference solution of Ghia et al. \cite{ghia1982high}. The results are in overall agreement with the reference solutions for all cases.
	
	\subsubsection{Lid-driven cavity flow around a flower-shaped object}\label{subsubsec:cavstar}
\begin{figure}
		\centering 
		\subfigure[$Re=1,t=10$]{\includegraphics[width=0.3\textwidth]{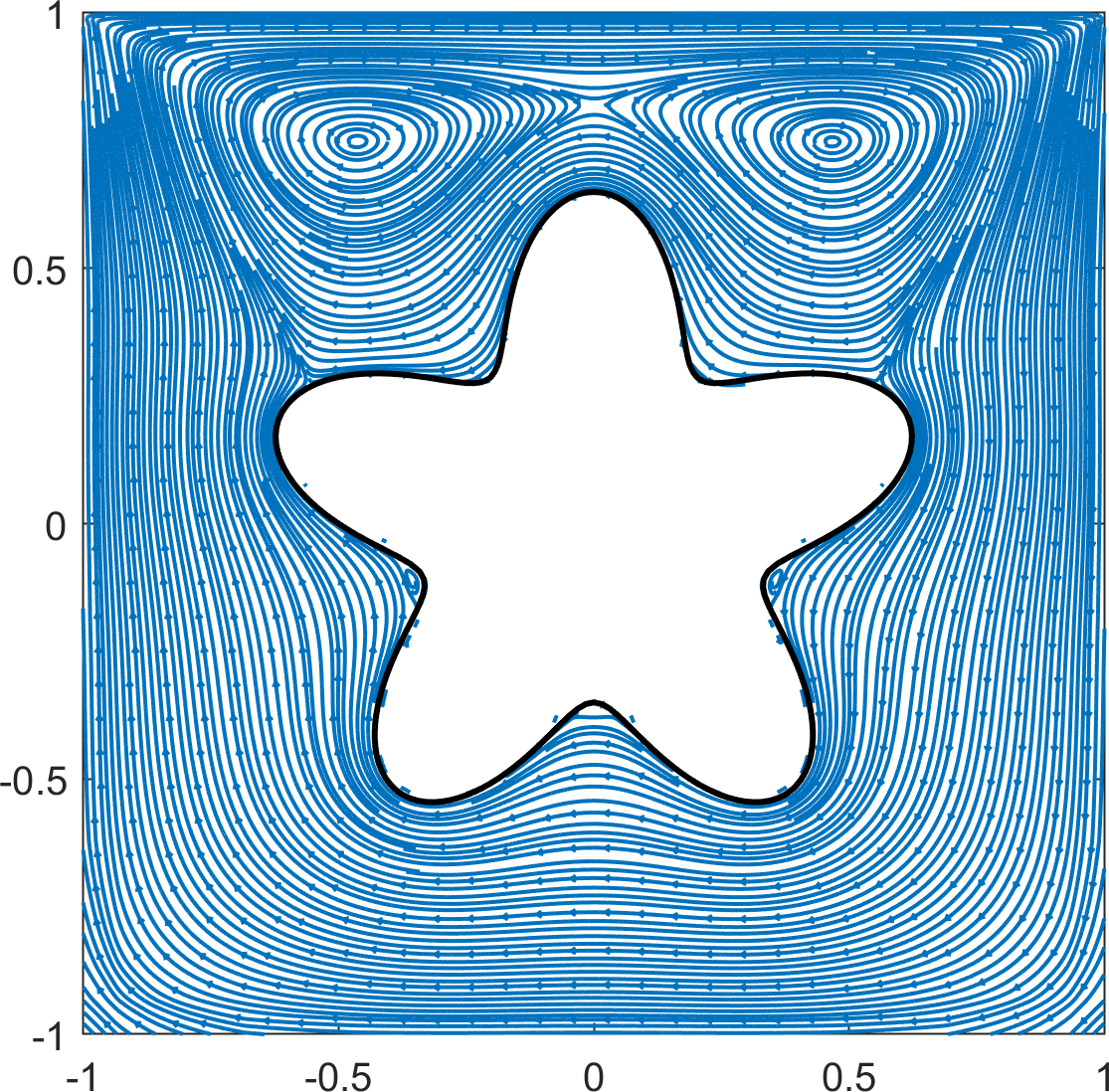}}	
		\subfigure[$Re=10,t=10$]{\includegraphics[width=0.3\textwidth]{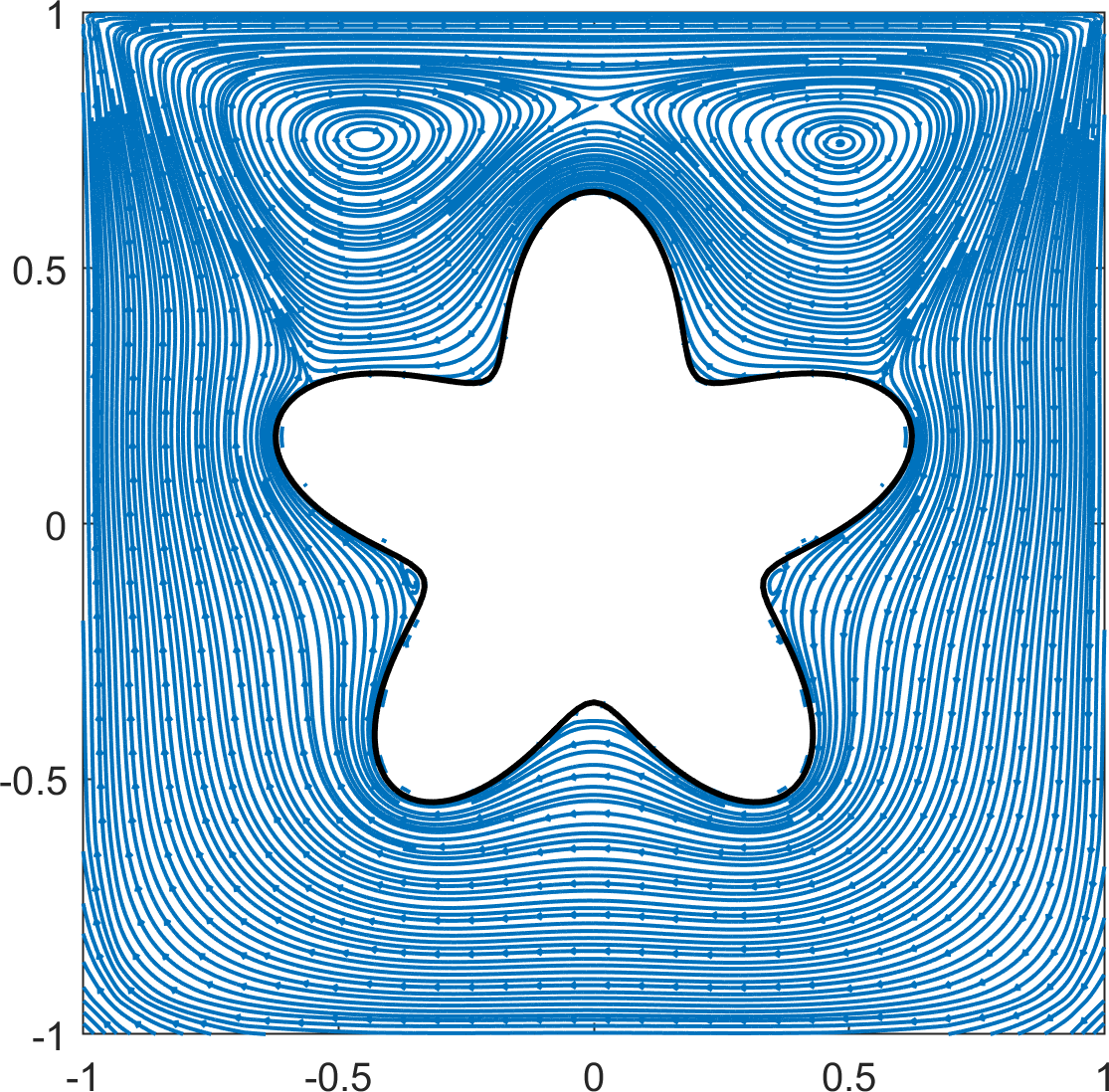}}	
		\subfigure[$Re=100,t=10$]{\includegraphics[width=0.3\textwidth]{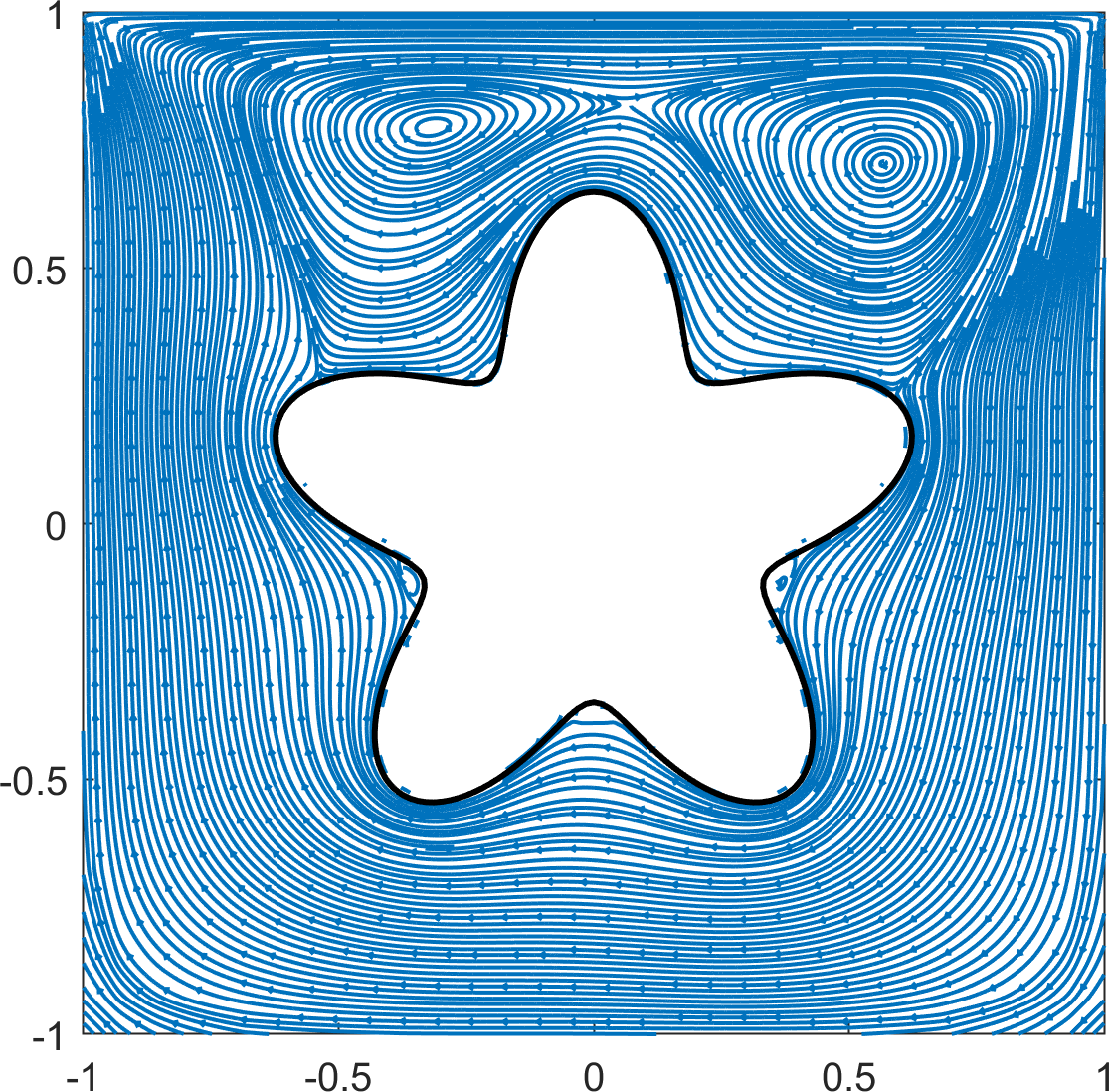}}\\ 
		\subfigure[$Re=1000,t=10$]{\includegraphics[width=0.3\textwidth]{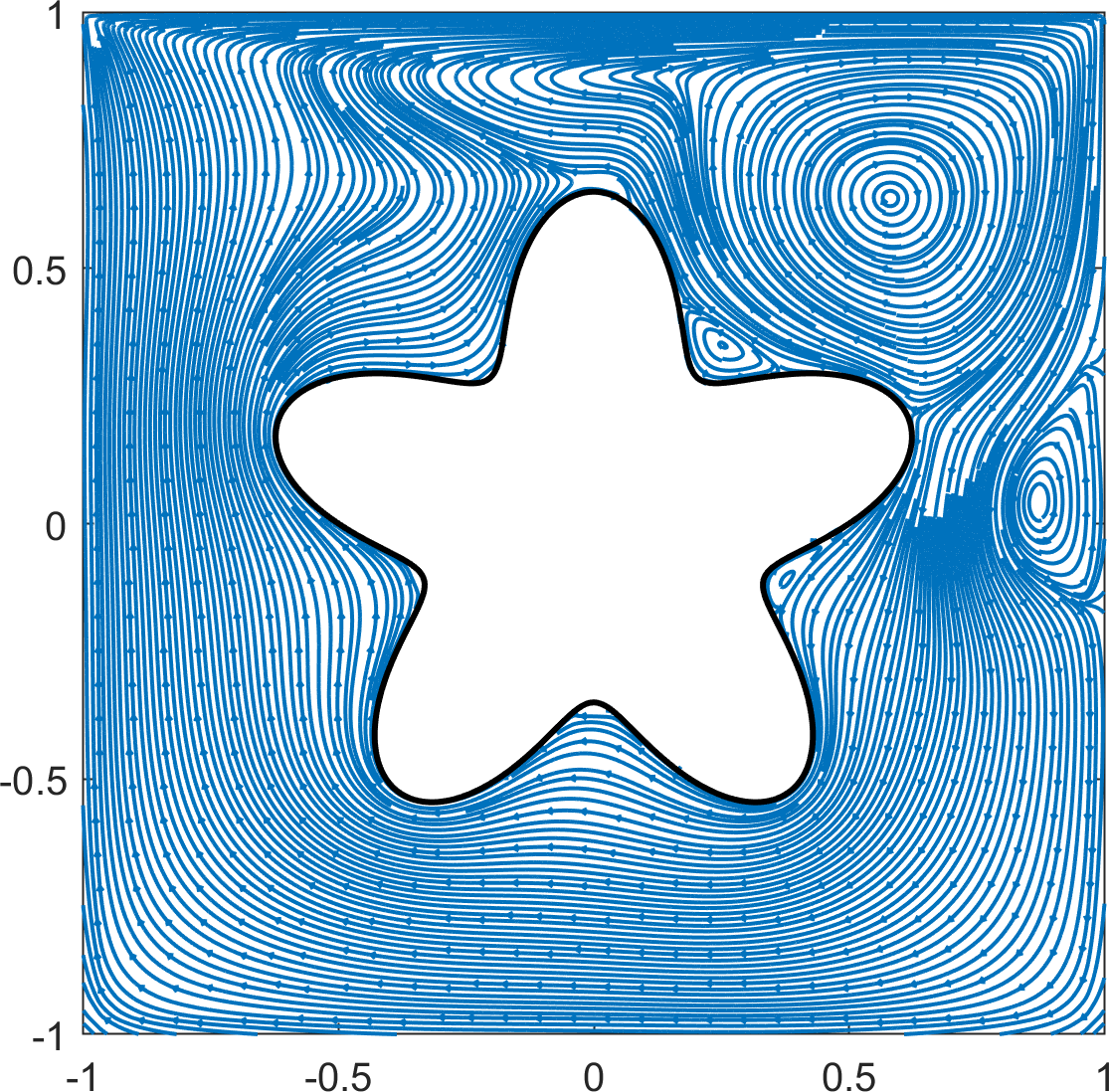}}	 
		\subfigure[$Re=1000,t=100$]{\includegraphics[width=0.3\textwidth]{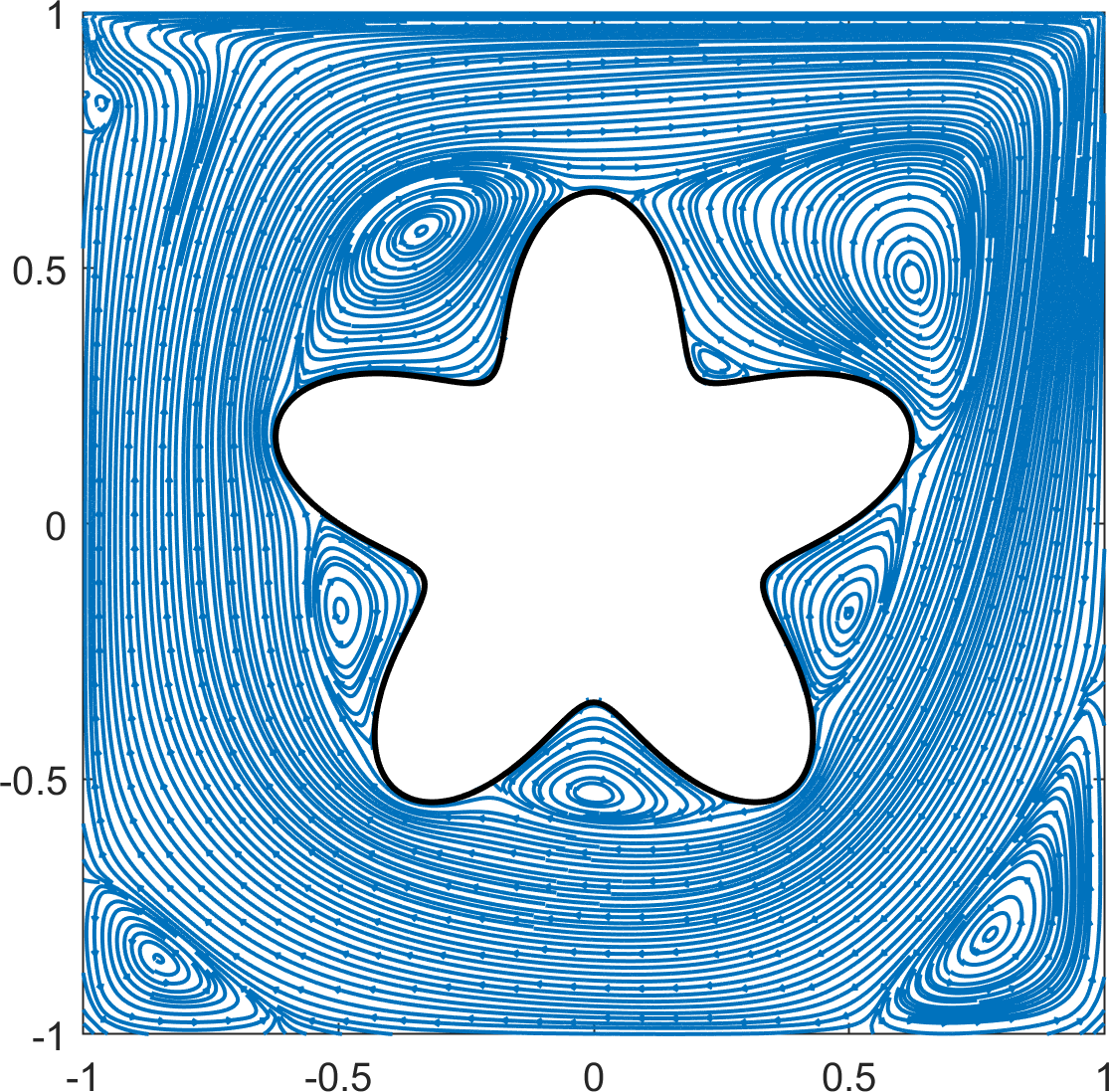}}	 
		\subfigure[$Re=5000,t=100$]{\includegraphics[width=0.3\textwidth]{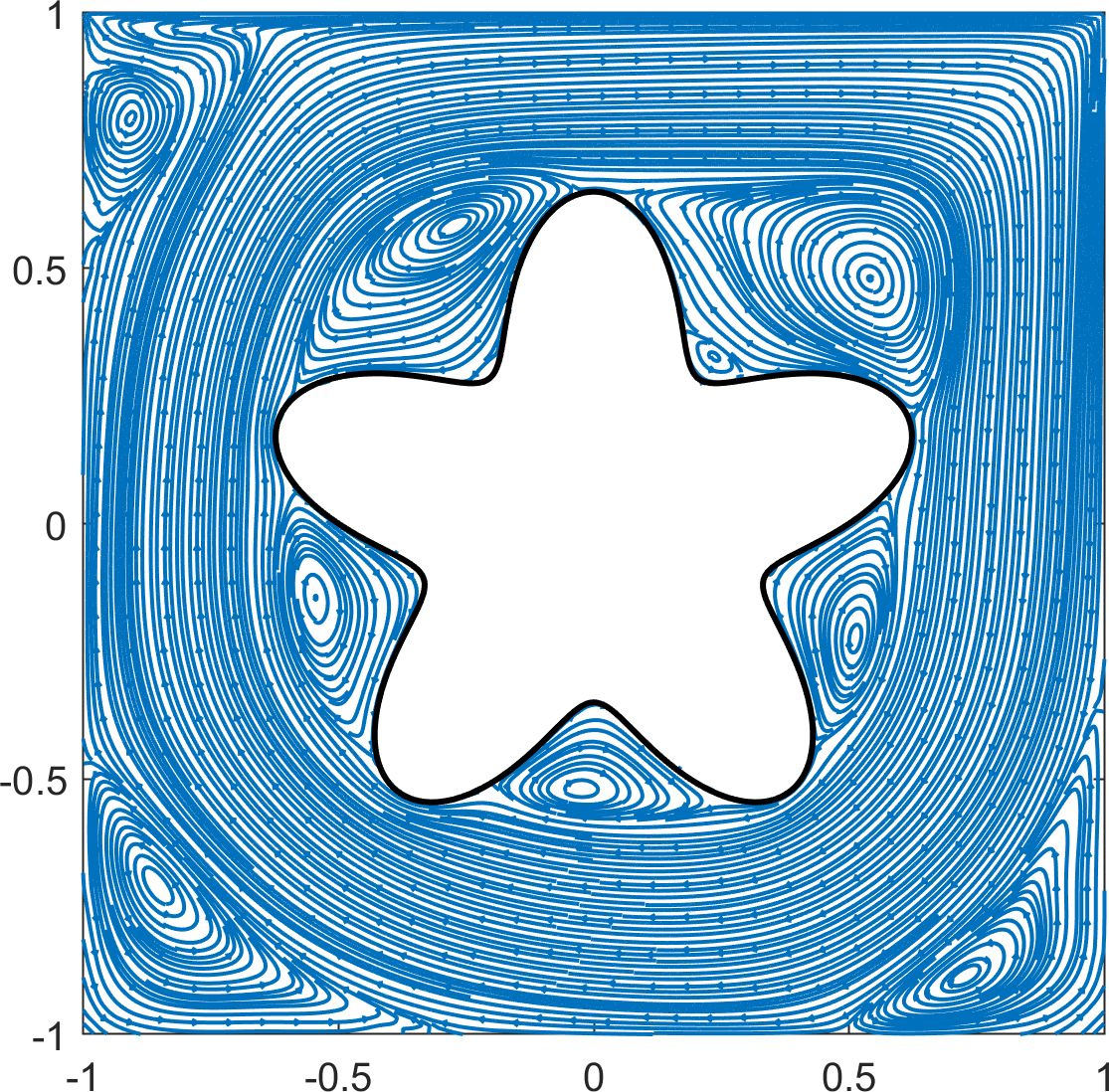}}	
		\caption{Streamlines of the results in example \ref{subsubsec:cavstar}.}\label{fig:cavstar}
	\end{figure}

	Using the experimental setup designed by Coco \cite{coco2020multigrid}, we tested the driven cavity flow on the computational domain $\Omega= \{ (x,y)\in \left[-1,1\right]^2 \mid  \phi(x,y) < 0\}$, where the level set function is given in polar coordinates as
	\[
	\phi(r,\theta)= -r +0.5+0.15\sin(5\theta).
	\]
	As in \ref{subsubsec:cavitysquare}, no-slip boundary conditions were used on all but the top wall, where the wall velocity was specified by the unit horizontal velocity, $\uvec=(1,0)$. A no-slip boundary condition was also imposed on the boundary of the flower-shaped object. 
	The cavity flow was studied for several Reynolds numbers: at $Re=1,10,100$ with a quadtree grid of level $5/7$ and at $Re=1000,5000$ with a quadtree grid of level $6/8$. Simulations were conducted up to $t=10$ and $t=100$ for low ($Re=1,10,100$) and high ($Re=1000,5000$) Reynolds numbers, respectively. Streamlines of the velocities at the terminal time are plotted in Figure \ref{fig:cavstar}. 
	
	The results for the low Reynolds numbers show good agreement with those reported in \cite{coco2020multigrid}. Two primary vortices are observed at the upper left and upper right near the flower-shaped object for $Re=1,10,100$, and the vortex at the upper left is not observed at $Re=1000$. According to \cite{coco2020multigrid}, the flower-shaped object prevents the formation of counterclockwise vortices in the lower corners. However, we found that counterclockwise vortices appear at the lower corners after a sufficiently long time at the high Reynolds numbers ($Re=1000,5000$). Furthermore, the primary vortex that appears in the standard square cavity flow problem is split into five vortices around the flower-shaped object.
	
	\subsection{Rotating flower}\label{subsec:rot_flower}
	As a third example, we consider the rotating flower-shaped object problem proposed by Coco \cite{coco2020multigrid}. A flower-shaped object rotates counterclockwise about the origin with angular velocity $\frac{2\pi}{5}$. Mathematically, the computational domain at time $t$ can be expressed as $\Omega= \{ (x,y)\in \left[-1,1\right]^2 |  \phi(x,y,t) < 0 \},$ where the level set function is given in polar coordinates by
	\[
	\phi(r,\theta,t)= -r +0.5+0.15\sin\left(5\theta -2\pi t\right).
	\]
	A no-slip boundary condition $\uvec=(0,0)$ is imposed on all four walls, and the boundary condition $\uvec= \frac{2}{5}\pi \left(-y,x\right)$ is imposed on the flower-shaped object, which agrees with the movement of the interface.
	
\begin{figure}
		\centering{}
		\mbox{
			\subfigure[$t=0$]{\includegraphics[width=0.3\textwidth]{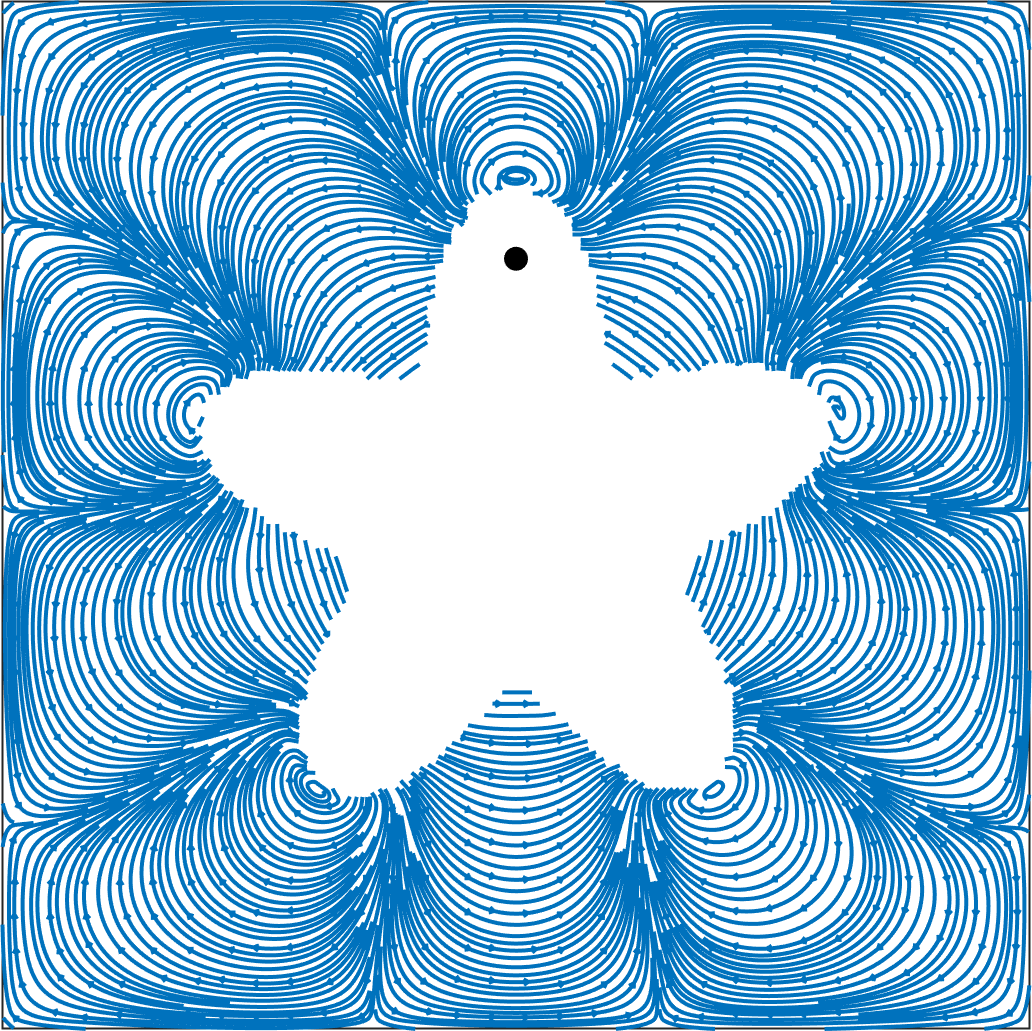}} $\ $
			\subfigure[$t=1.25$]{\includegraphics[width=0.3\textwidth]{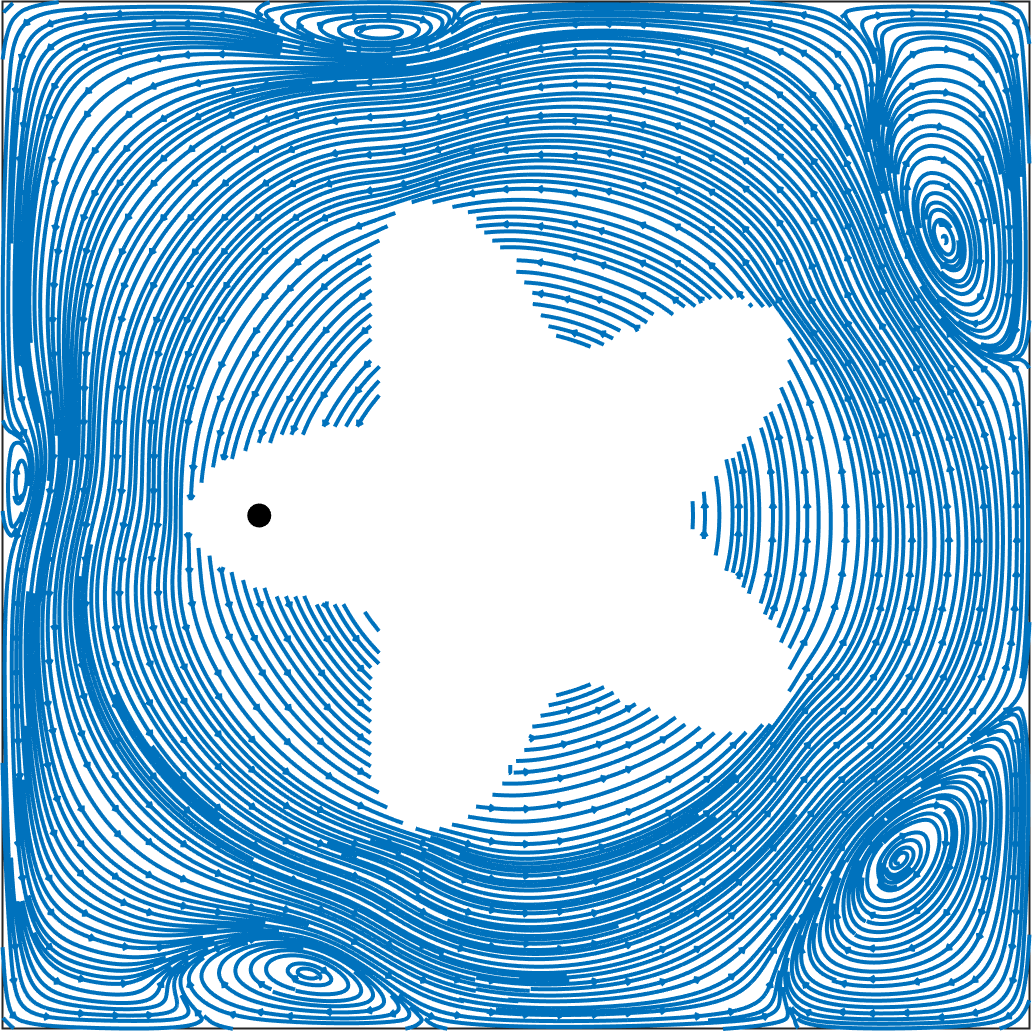}}$\ $
			\subfigure[$t=2.5$]{\includegraphics[width=0.3 \textwidth]{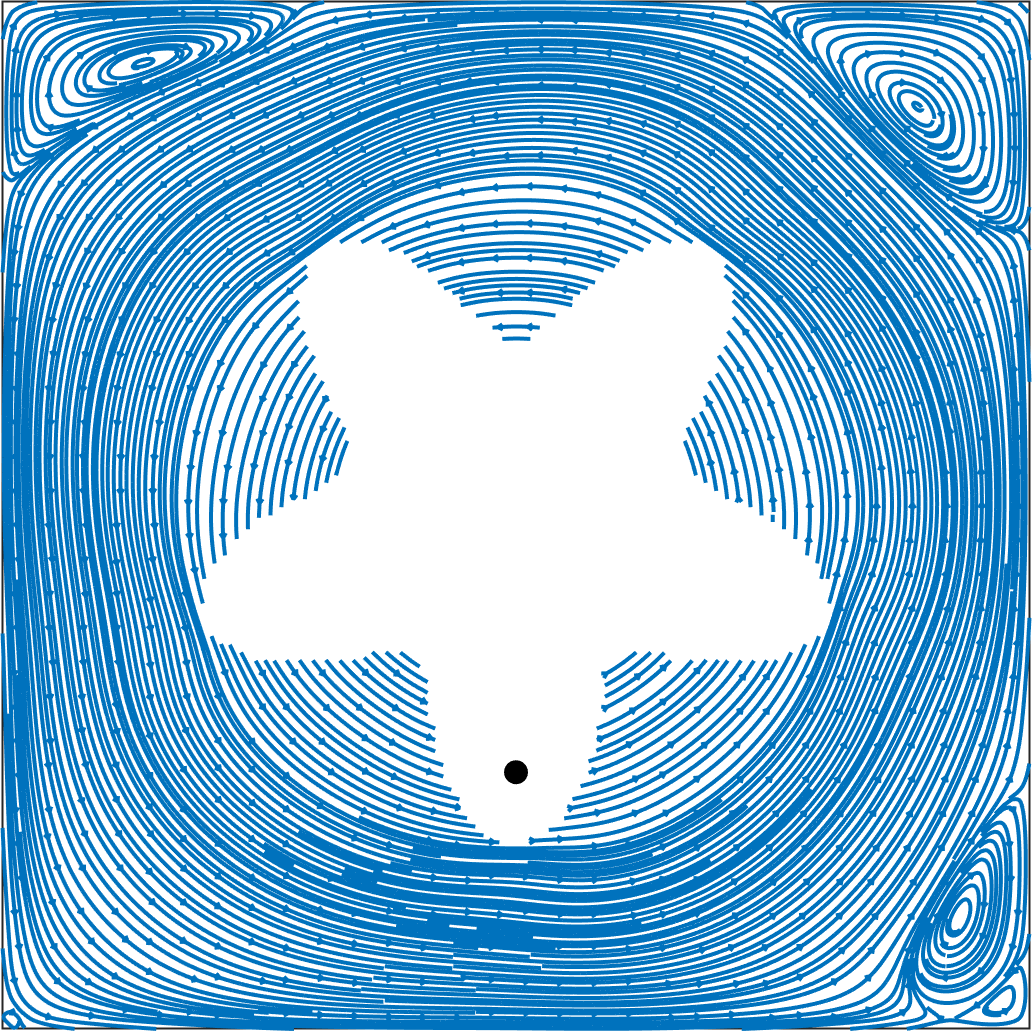}}$\ $
		}
		
		\centering{}
		\mbox{
			\subfigure[$t=3.75$]{\includegraphics[width=0.3\textwidth]{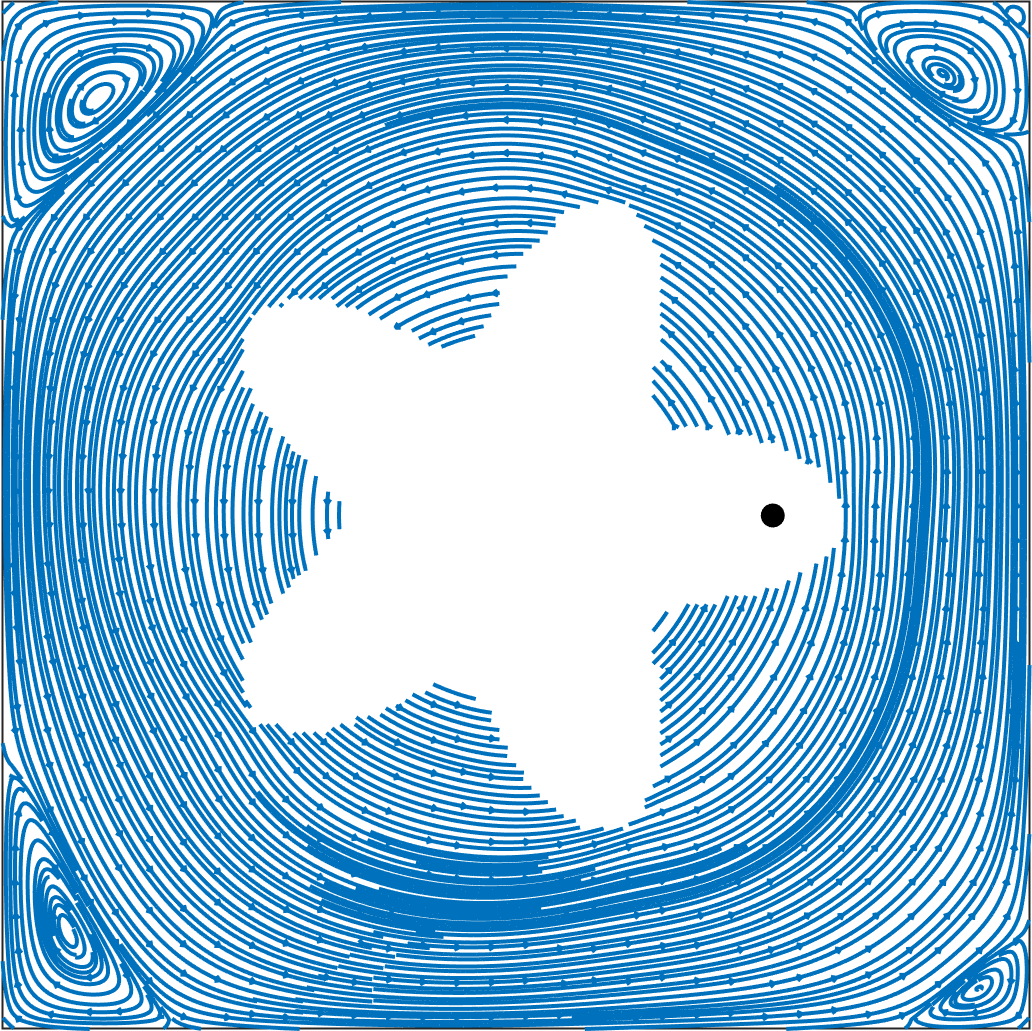}} $\ $
			\subfigure[$t=5$]{\includegraphics[width=0.3\textwidth]{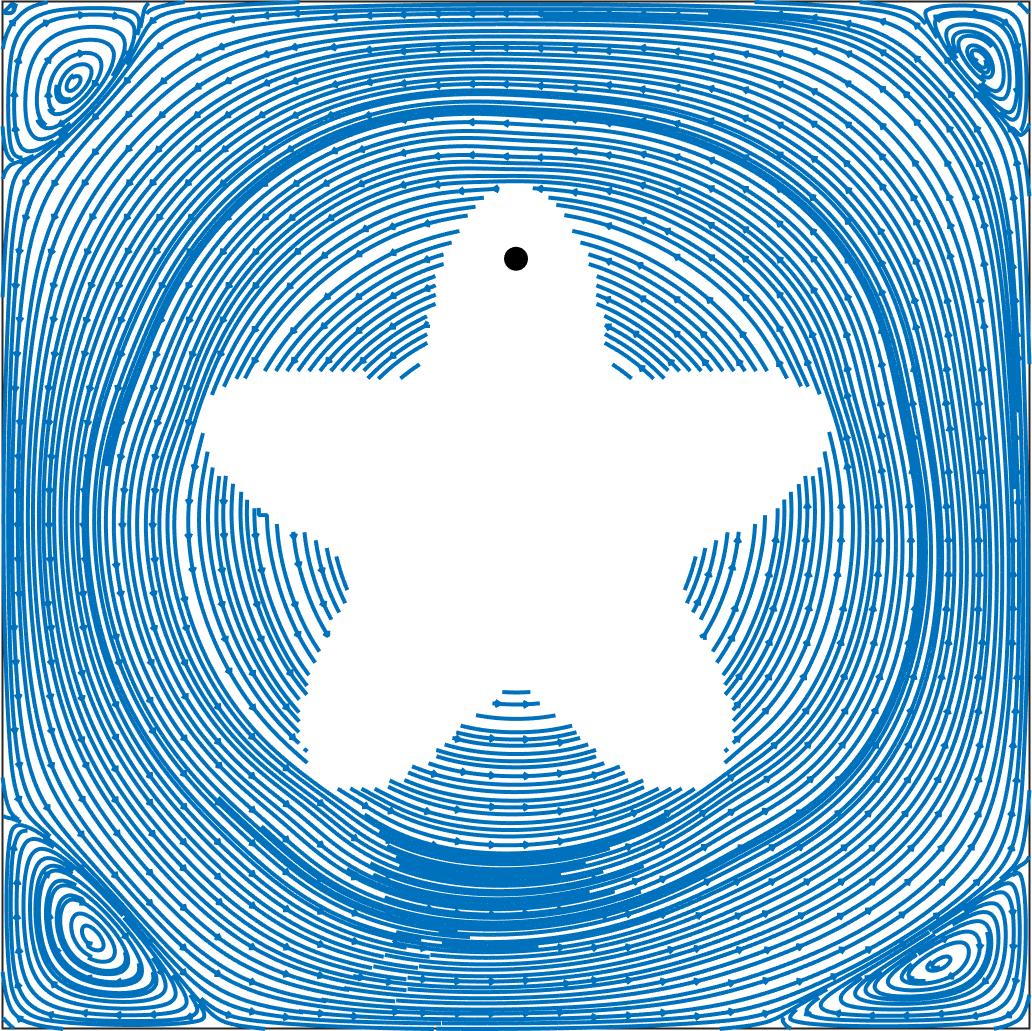}}$\ $
			\subfigure[$t=6.25$]{\includegraphics[width=0.3\textwidth]{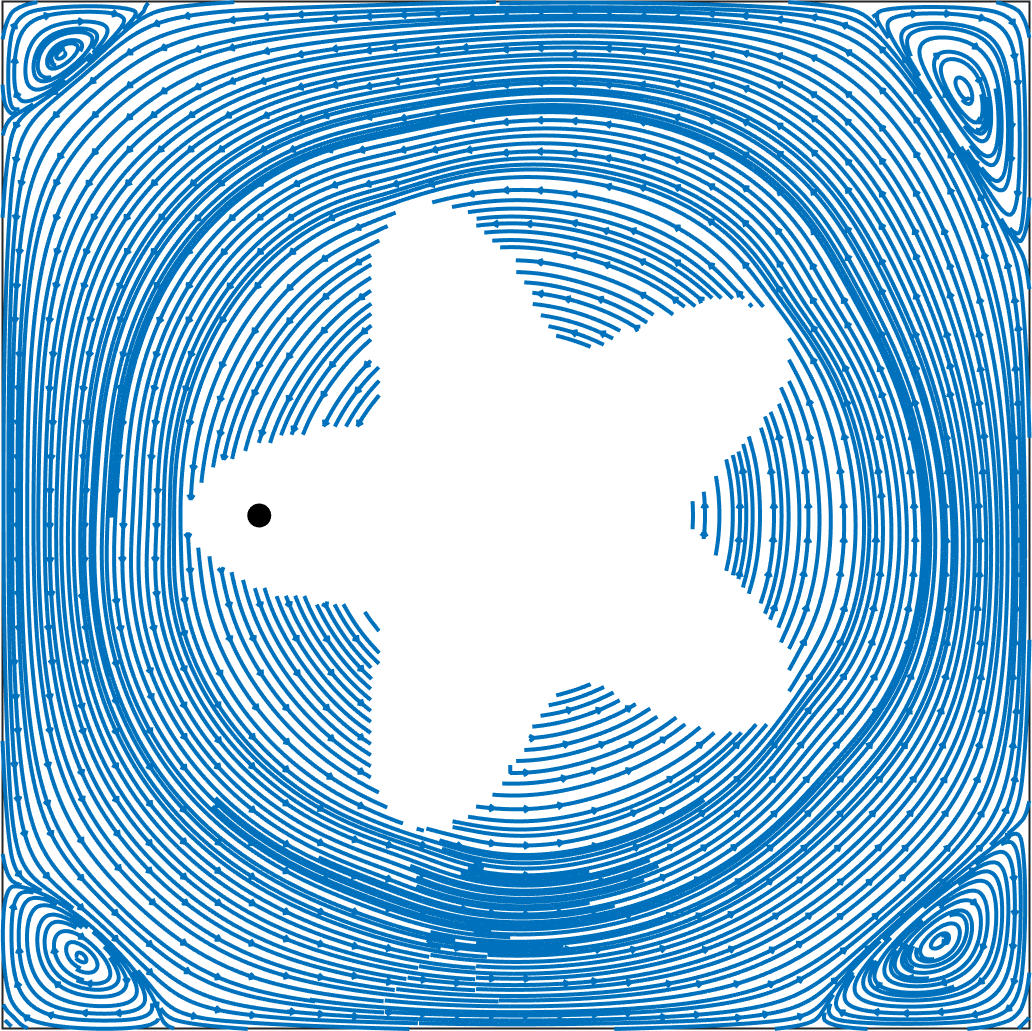}}$\ $
		}
		
		\centering{}
		\mbox{
			\subfigure[$t=7.5$]{\includegraphics[width=0.3\textwidth]{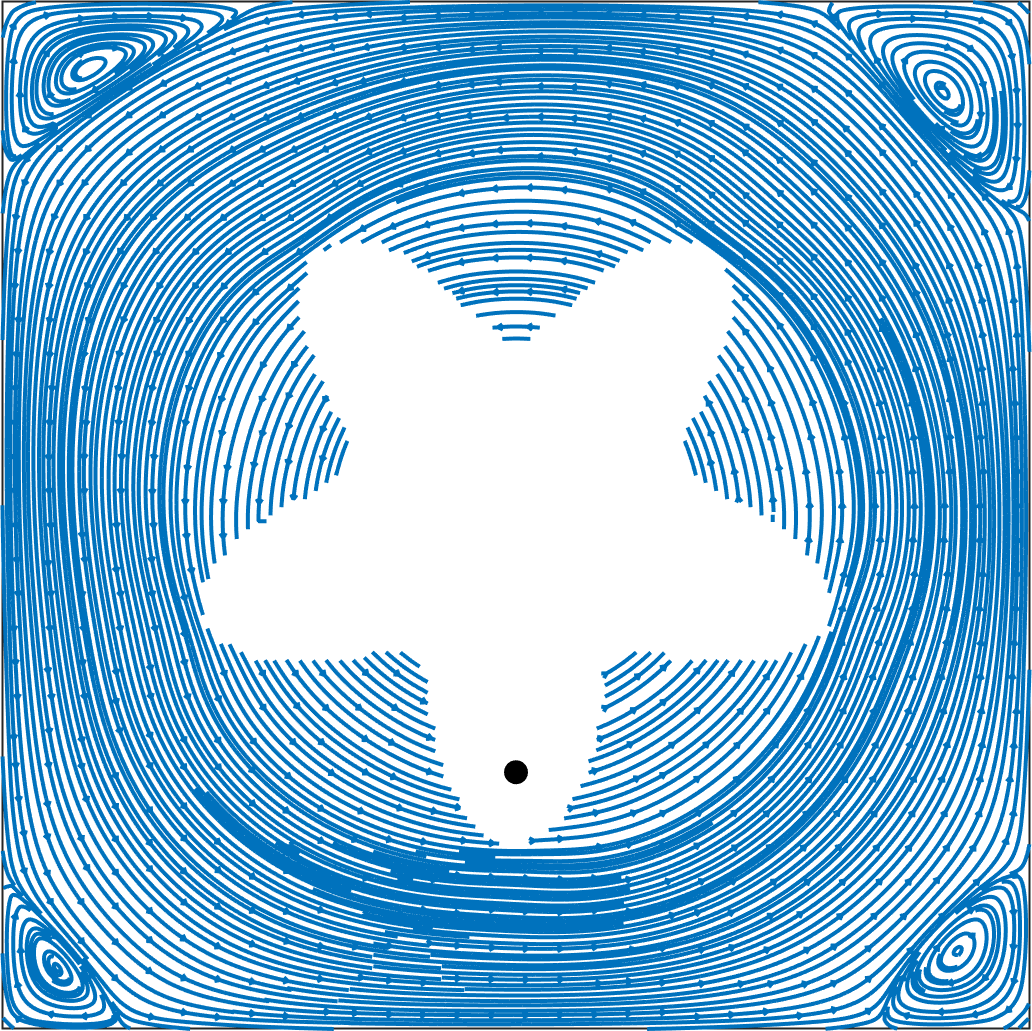}} $\ $
			\subfigure[$t=8.75$]{\includegraphics[width=0.3\textwidth]{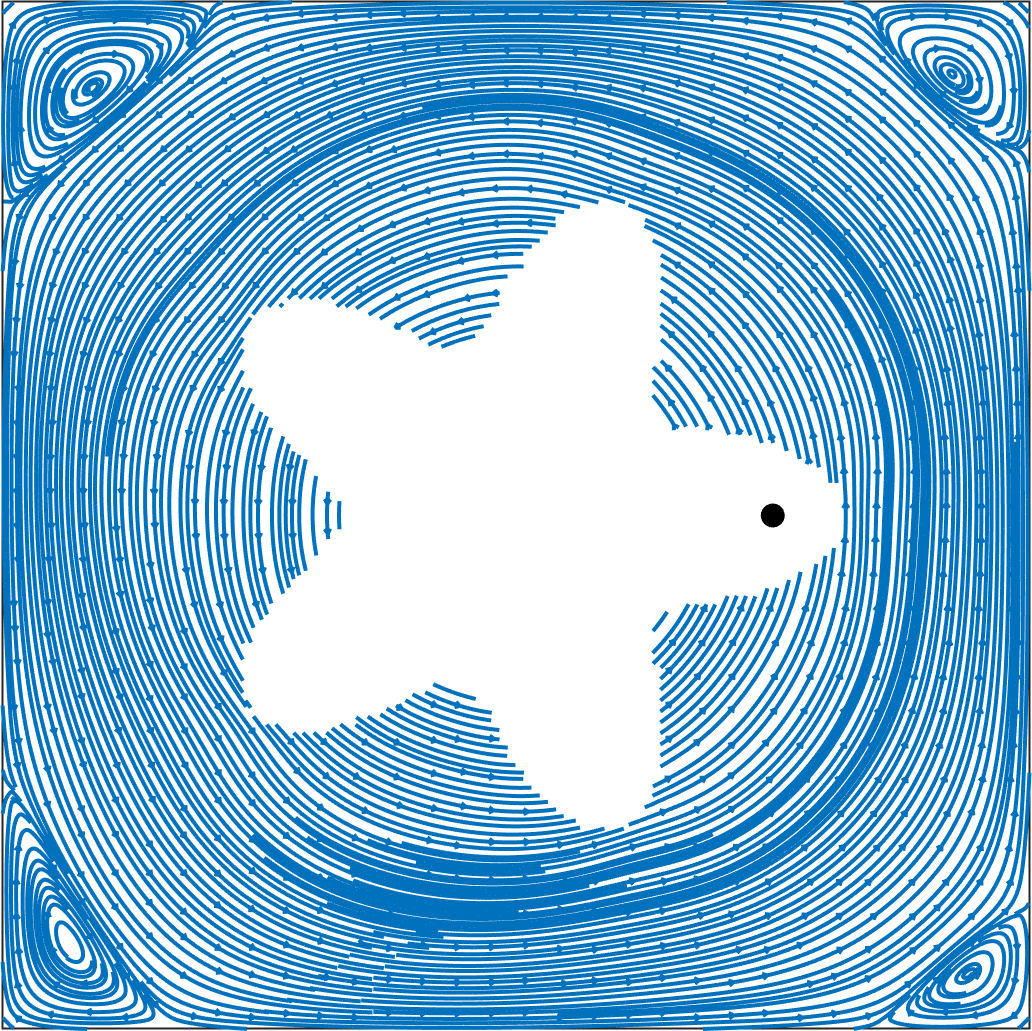}}$\ $
			\subfigure[$t=10$]{\includegraphics[width=0.3\textwidth]{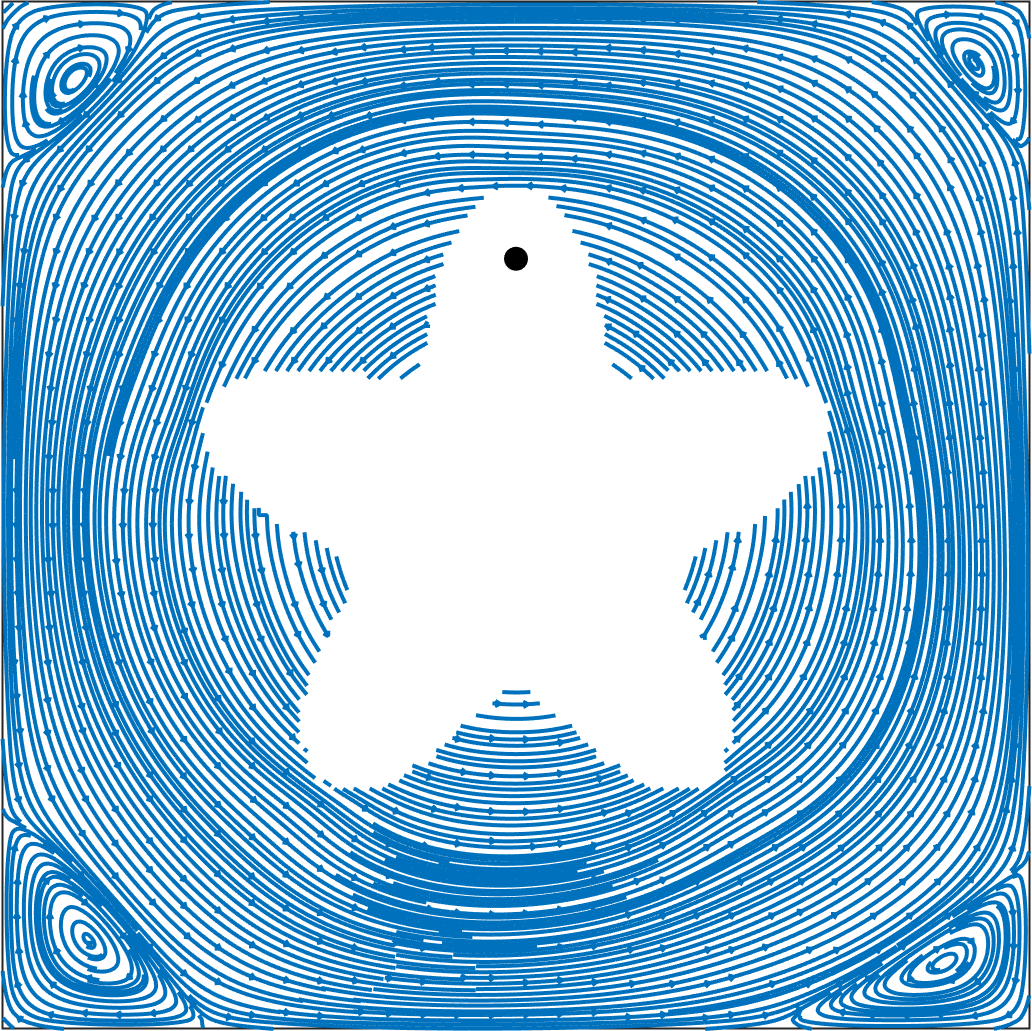}}$\ $
		}
		\caption{Streamlines of the results in example \ref{subsec:rot_flower} for $Re=100$.}\label{fig:rotstar100}
	\end{figure}

\begin{figure}
		\centering{}	
		\mbox{
			\subfigure[$t=0$]{\includegraphics[width=0.3\textwidth]{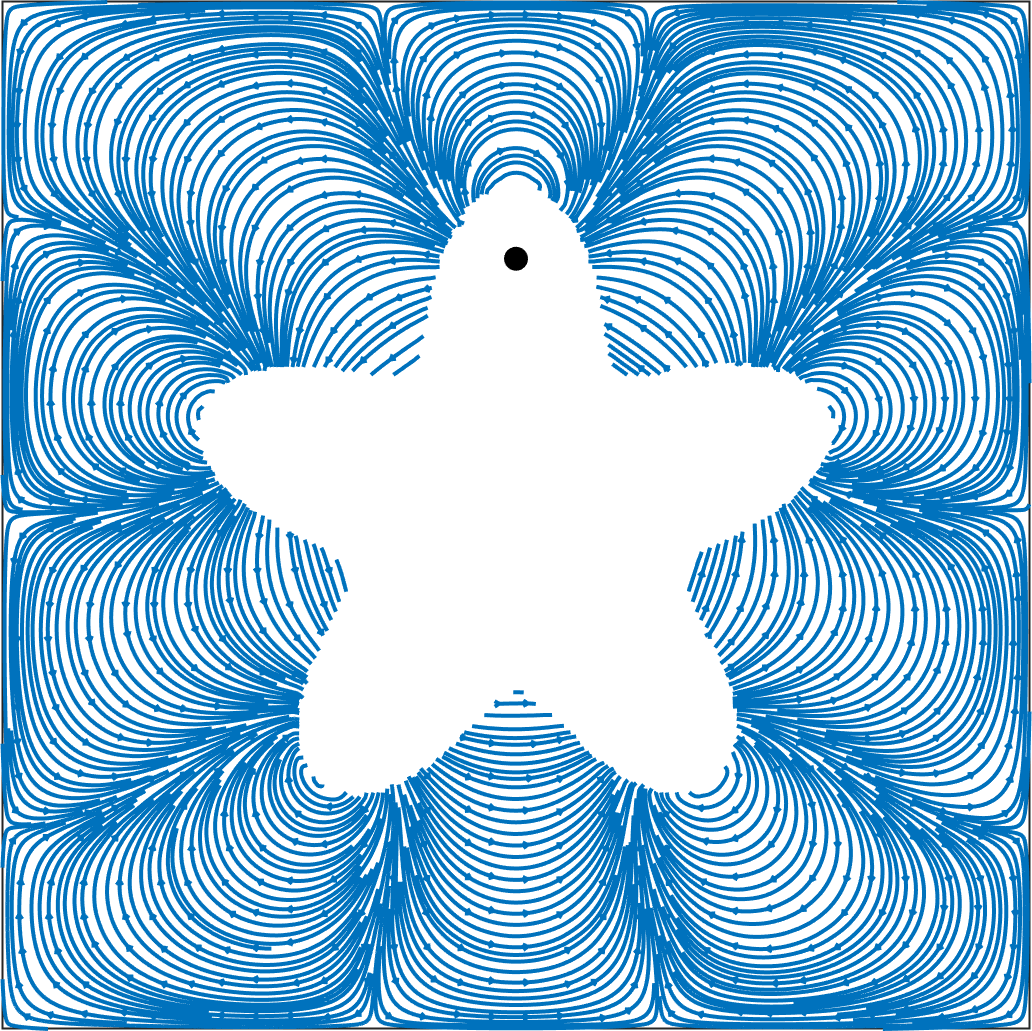}}  $\ $	
			\subfigure[$t=10$]{\includegraphics[width=0.3\textwidth]{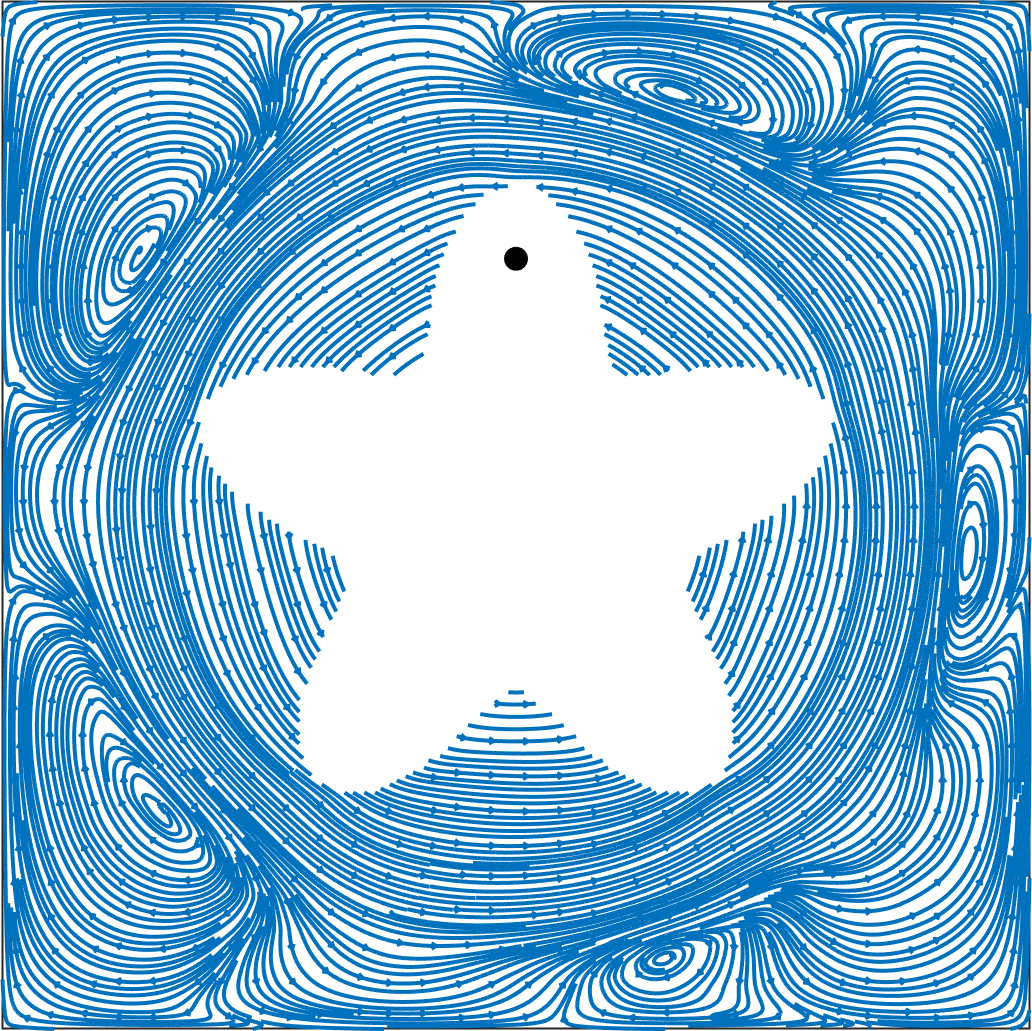}}$\ $
			\subfigure[$t=20$]{\includegraphics[width=0.3\textwidth]{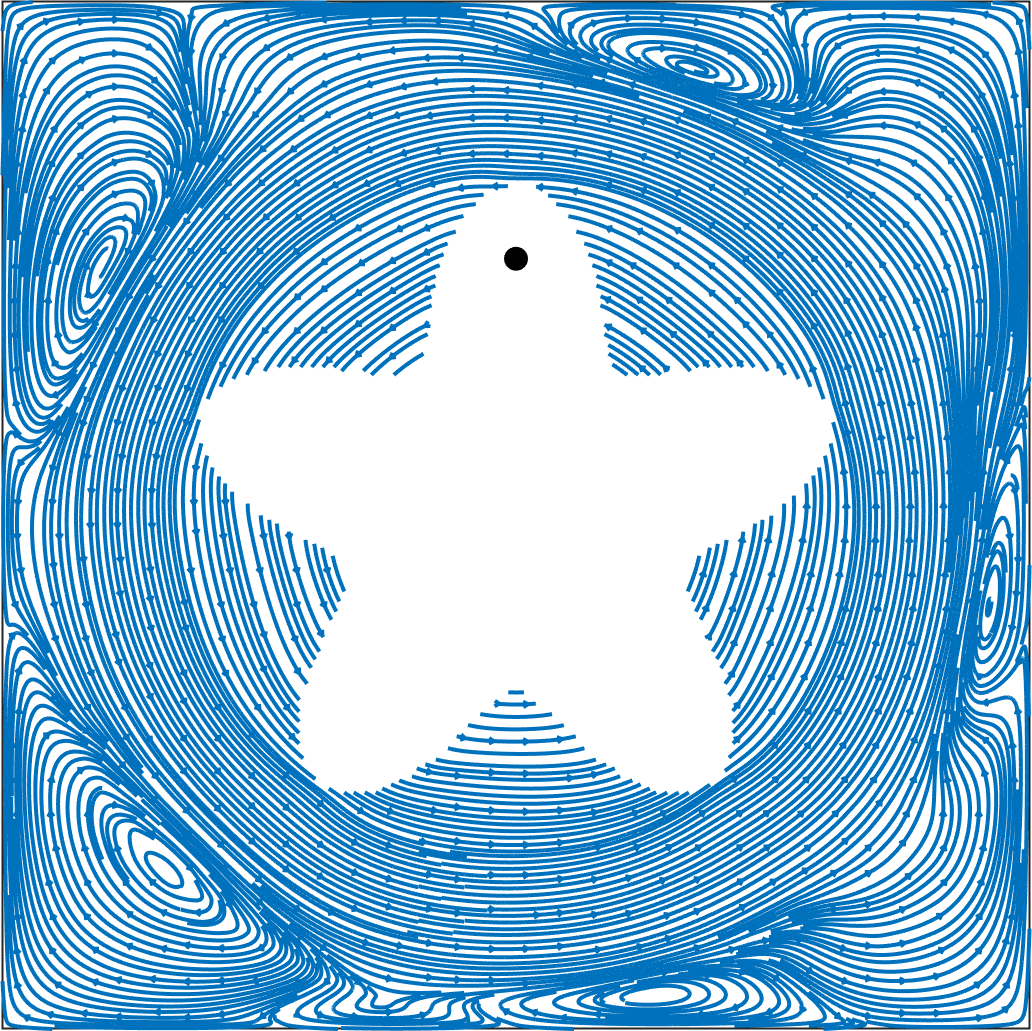}}$\ $
		}
		
		\centering{}	
		\mbox{
			\subfigure[$t=30$]{\includegraphics[width=0.3\textwidth]{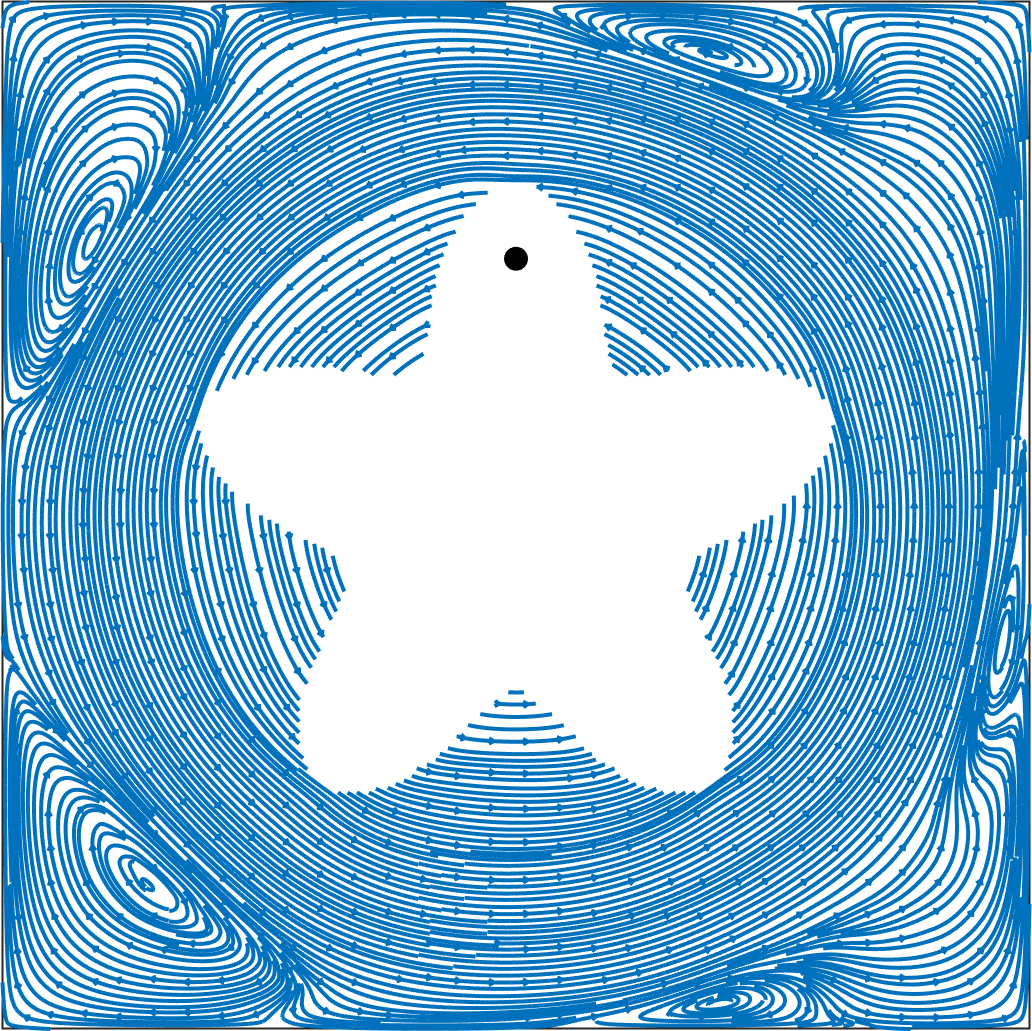}}	 $\ $		
			\subfigure[$t=40$]{\includegraphics[width=0.3\textwidth]{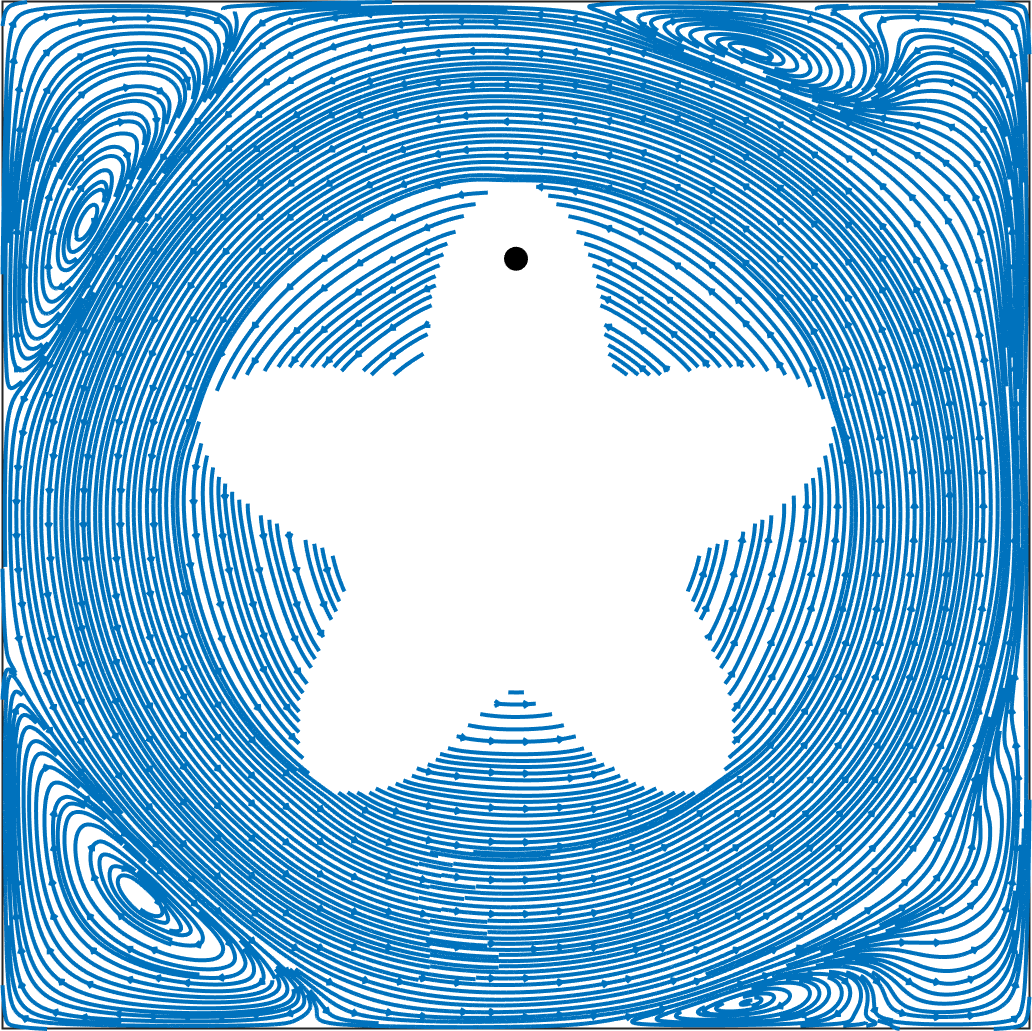}}  $\ $	
			\subfigure[$t=50$]{\includegraphics[width=0.3\textwidth]{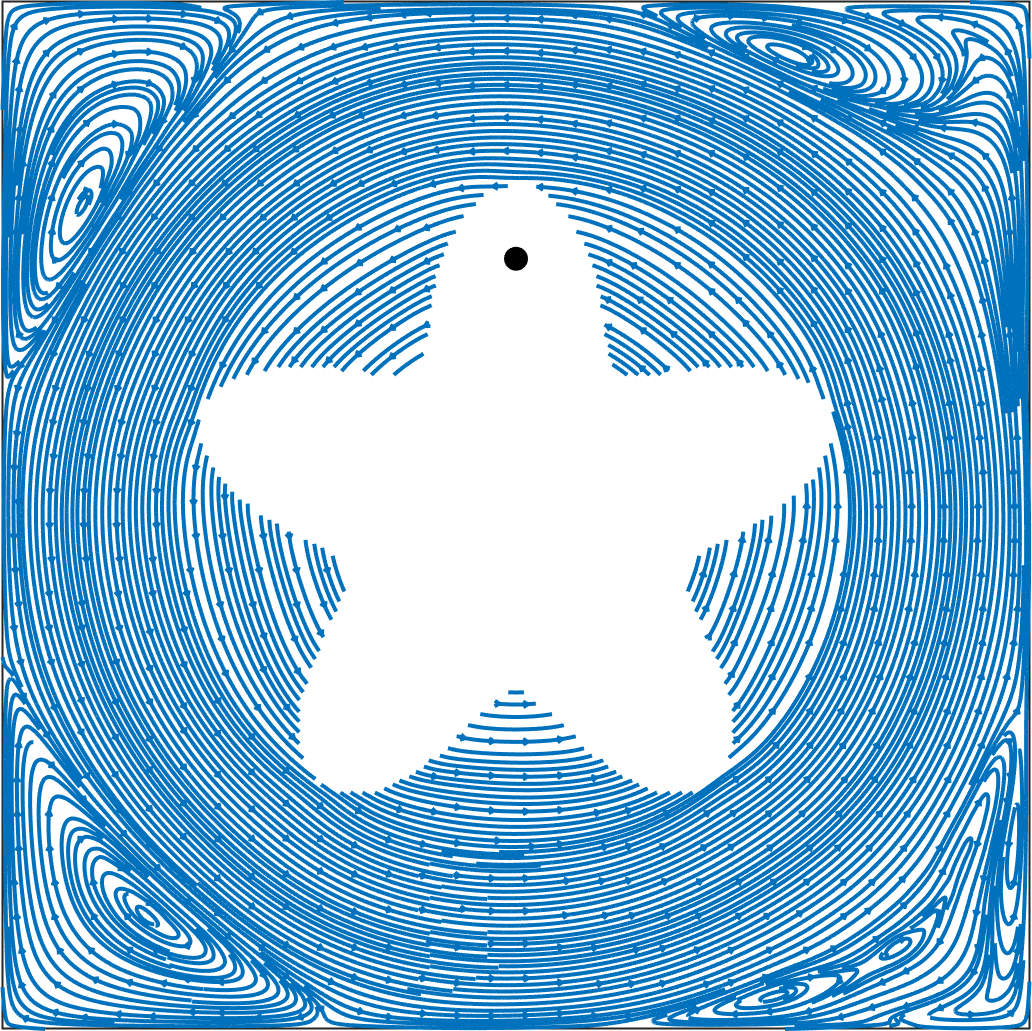}}$\ $	
		}
		
		\centering{}	
		\mbox{
			\subfigure[$t=60$]{\includegraphics[width=0.3\textwidth]{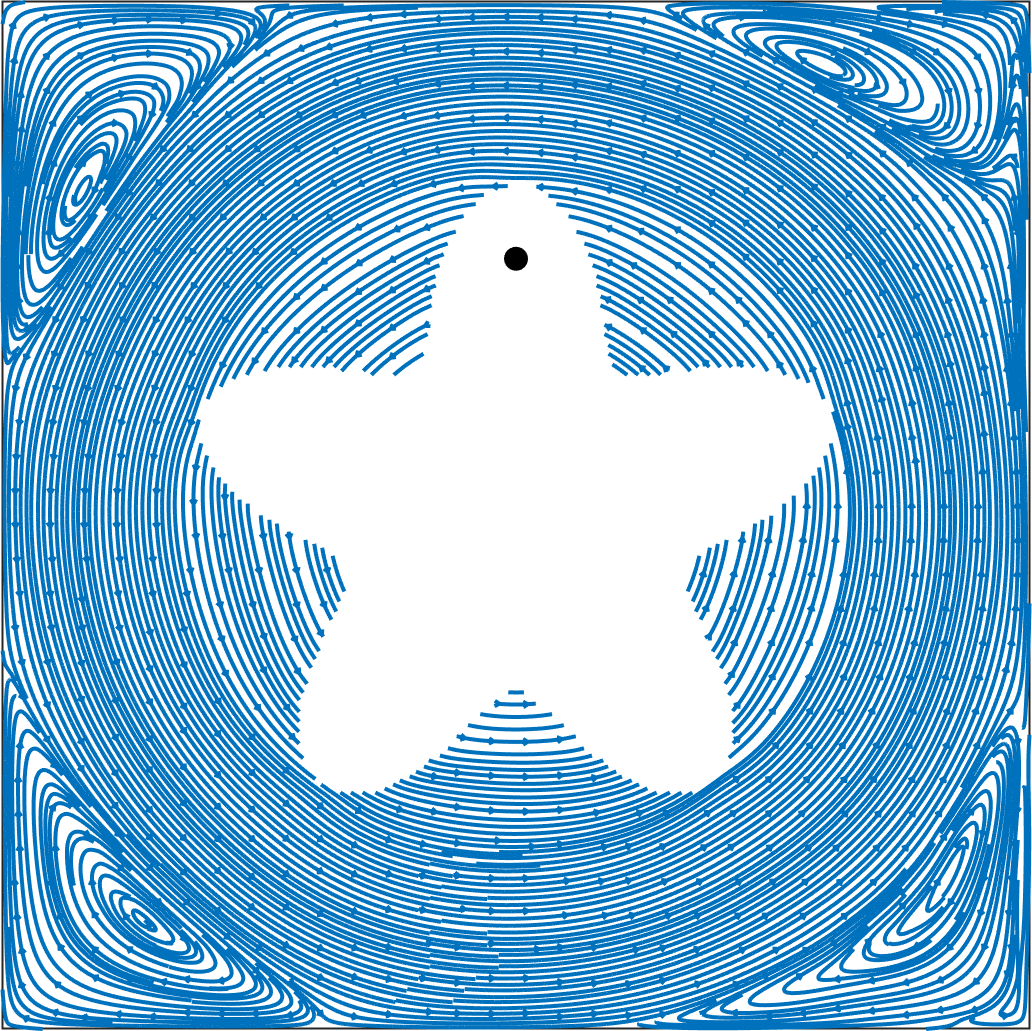}}$\ $		
			\subfigure[$t=70$]{\includegraphics[width=0.3\textwidth]{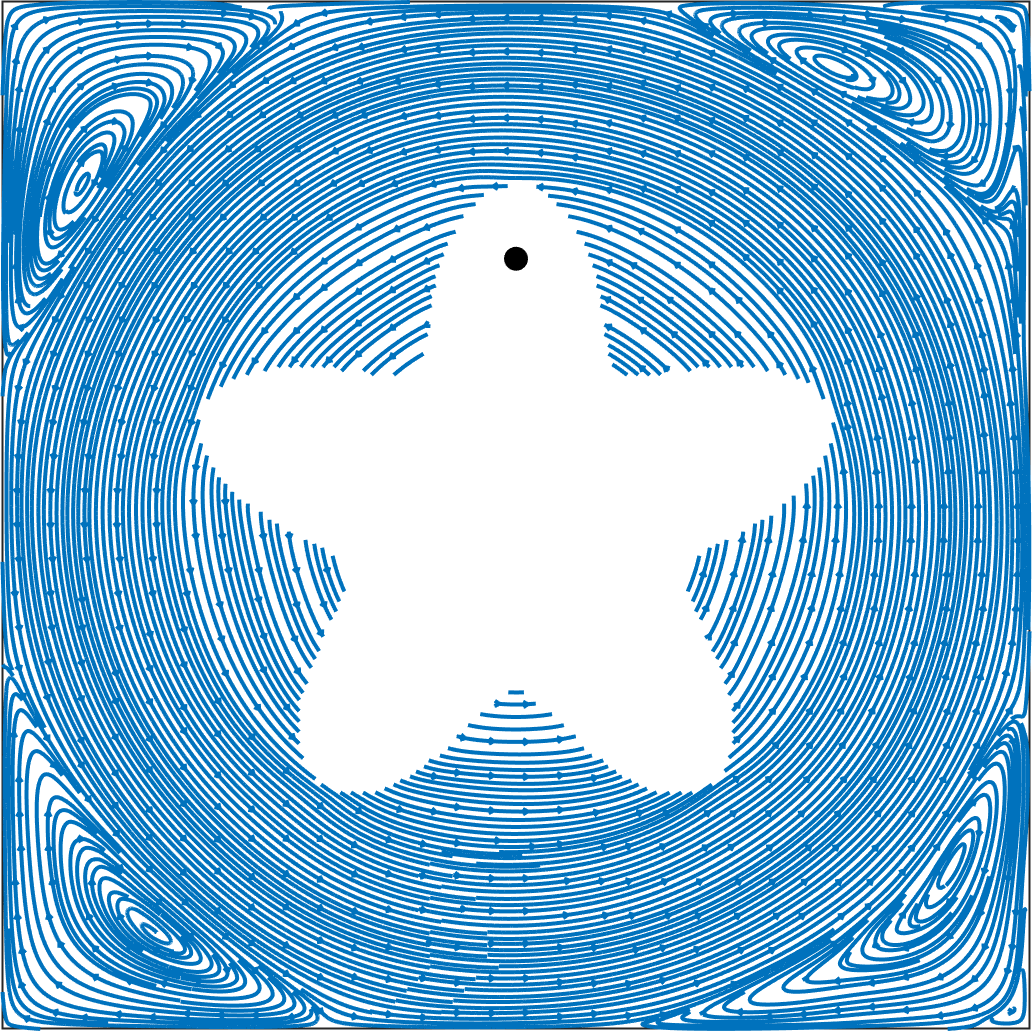}}$\ $
			\subfigure[$t=80$]{\includegraphics[width=0.3\textwidth]{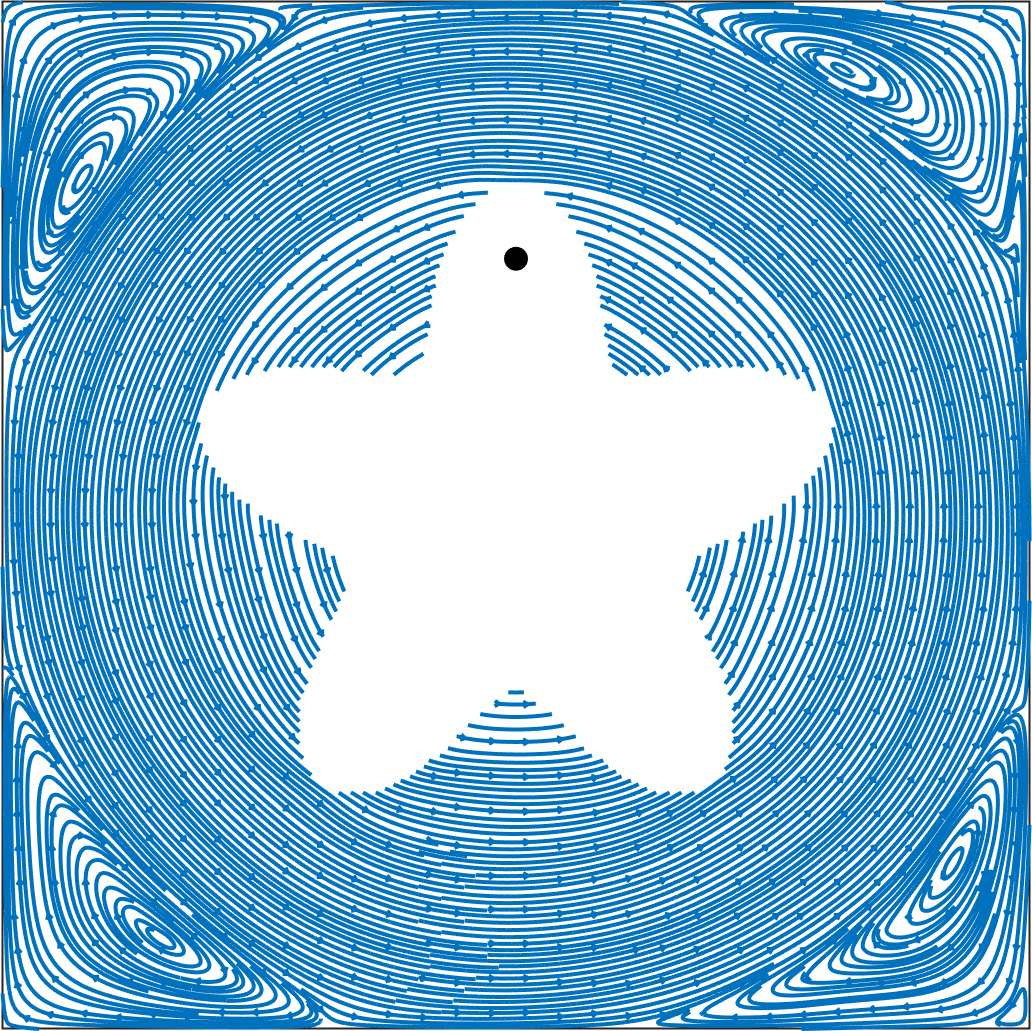}}$\ $
		}
		\caption{Streamlines of the results of example \ref{subsec:rot_flower} for $Re=10000$.}\label{fig:rotstar10000}
	\end{figure}

	Numerical experiments were conducted on a quadtree grid of level $6/8$ for two Reynolds numbers, $Re=100$ and $10000$, whereas Coco performed experiments using only $Re=100$. As in \cite{coco2020multigrid}, we plotted streamlines of the resulting fluid flow in the test with $Re=100$ every quarter rotation, that is, every 1.25 s, as shown in Figure \ref{fig:rotstar100}. The black dot is plotted on the same petal. Figure \ref{fig:rotstar100} shows that the positions of the vortices agree with those presented in \cite{coco2020multigrid}. After a full rotation, we observe four clockwise vortices in the corners. For the high Reynolds number, $Re=10000$, streamlines appear every 10 s in Figure \ref{fig:rotstar10000}. In contrast to the $Re=100$ case, vortices are formed in each corner around $t=60$. The results from $t=60$ to $80$ shows that the fluid converges to a steady-state solution.

	\subsection{Application to two-dimensional rising bubble }\label{subsubsec:bubble}
	Although the proposed method is applicable only to single-phase flow, it can be easily extended to two-phase flow problems. The level set function $\phi$ divides $\Omega$ into $\Omega^+=\{(x,y)\in \Omega\mid \phi > 0 \}$ and $\Omega^-=\{(x,y)\in \Omega\mid \phi <
	0 \}$. The densities $\rho^\pm$ and viscosities $\mu^\pm$ for $\Omega^\pm$ may differ. According to the delta function equation, two-phase flow with surface tension can be described as 
	\begin{equation}\label{eq:NS_multiphase}
		\begin{gathered}
			\rho^\pm (\uvec + \uvec \cdot \nabla \uvec)= \mu^\pm \lap \uvec -\nabla p +\sigma \delta \kappa \normal+\rho \mathbf{g},\\
			\nabla \cdot \uvec =0.
		\end{gathered}
	\end{equation}
	Here, \(\normal \) represents the outward normal pointing from $\Omega^+$ to the interface, $\mathbf{g}=(0,-g)$ represents gravity, $\kappa$ denotes the mean curvature of the interface, and $\sigma$ is the surface tension coefficient. 
	A time-step procedure similar to that for the single-phase flow described in Section \ref{sec:NumericalMethod} is used to simulate two-phase flow, as follows:
	\begin{enumerate}
		\item From time level $t^n$, compute the appropriate $\Delta t_n$, and define $t^{n+1}= t^n + \Delta t_n$. Advect the level set to obtain $\phi^{n+1}$ at time level $t^{n+1}$, and then reinitialize it.
		\item Generate a quadtree grid with respect to the level set $\phi^{n+1}$.
		\item Solve for $\uvec^{n+1}$ and $p^{n+1}$ at time level $t^{n+1}$.
	\end{enumerate}
	For advection and the reinitialization of the level set function, we adopt the second-order finite difference method on a quadtree grid \cite{min2007second}. By using the numerical delta function, the Navier--Stokes equations can be discretized as 
	\begin{equation*}
		\begin{gathered}
			\rho^{n+1} 		\left( \alpha_0 \uvec^{n+1} +\alpha_1  \uvec_d^n + \alpha_2 \uvec_d^{n-1} \right)= \mu \lap \uvec^{n+1} - \nabla p^{n+1} + \sigma\delta(\phi^{n+1})\kappa +\rho^{n+1} \mathbf{g},\\
			\nabla \cdot \uvec^{n+1} =0,
		\end{gathered}
	\end{equation*}
	where $\alpha_0,\alpha_1,\alpha_2$ are determined according to the second-order backward differences as follows:
	\begin{equation*}
		\begin{aligned}
			\alpha_0&=  \left(2\Delta t_n + \Delta t_{n-1}\right) / ( \Delta t_{n}^2+\Delta t_n \Delta t_{n-1} ),\quad
			\alpha_2  &= \Delta t_n / (\Delta t_n  \Delta t_{n-1} + \Delta t_{n-1}^2),\quad  \alpha_1 &= -\alpha_0-\alpha_2.
		\end{aligned}
	\end{equation*}
	The discretization of $\nabla p$ and $\nabla \cdot \uvec$ does not differ from the discretization introduced in Section \ref{sec:NumericalMethod}. Using notation similar to that in Figure \ref{fig:stencil_stokes}, the terms $\rho$, $\mu \lap \uvec$, and $\sigma \delta \kappa \normal$ at $\mathbf{x}_0$ are discretized by
	\begin{equation*}
		\begin{aligned}
			\rho (\mathbf{x}_0)&= \left(\rho^+- \rho^-\right) H(\phi_0) + \rho^- ,\\	
			\mu \lap \uvec (\mathbf{x}_0)&= \left[ \mu_R \frac{\mathbf{u}_{R}-\mathbf{u}_{0}}{d_{R}}-\mu_L \frac{\mathbf{u}_{0}-\mathbf{u}_{L}}{d_{L}}\right]\frac{2}{d_{R}+d_{L}}+\left[\mu_T \frac{\mathbf{u}_{T}-\mathbf{u}_{0}}{d_{T}}-\mu_B \frac{\mathbf{u}_{0}-\mathbf{u}_{B}}{d_{B}}\right]\frac{2}{d_{T}+d_{B}},\\
			\sigma \delta \kappa \normal(\mathbf{x}_0) &=\sigma \delta(\phi_0) \kappa \normal, \\
		\end{aligned}
	\end{equation*}
	where $\mu_R= \left(\mu^+- \mu^-\right) H(\frac{\phi_0 + \phi_R}{2} ) + \mu^-$ and $\mu_L,\mu_T,\mu_B$ are defined similarly. Here, $H$ and $\delta$ are numerical Heaviside and delta functions and are defined as
	\begin{equation*}
		H(\phi)=\left\{\begin{array}{ll}
			0 & \phi<-\epsilon \\
			\frac{1}{2}+\frac{\phi}{2 \epsilon}+\frac{1}{2 \pi} \sin \left(\frac{\pi \phi}{\epsilon}\right) & -\epsilon \leq \phi \leq \epsilon, \\
			1 & \epsilon<\phi
		\end{array} \quad \delta(\phi)=\left\{\begin{array}{ll}
			0 & \phi<-\epsilon \\
			\frac{1}{2 \epsilon}+\frac{1}{2 \epsilon} \cos \left(\frac{\pi \phi}{\epsilon}\right) & -\epsilon \leq \phi \leq \epsilon \\
			0 & \epsilon<\phi
		\end{array}\right.\right.
	\end{equation*}
	with $\epsilon=1.5 \Delta x$.
	
\begin{figure}
		\centering{}
		\mbox{
			\subfigure[][]{\includegraphics[width=0.45\textwidth]{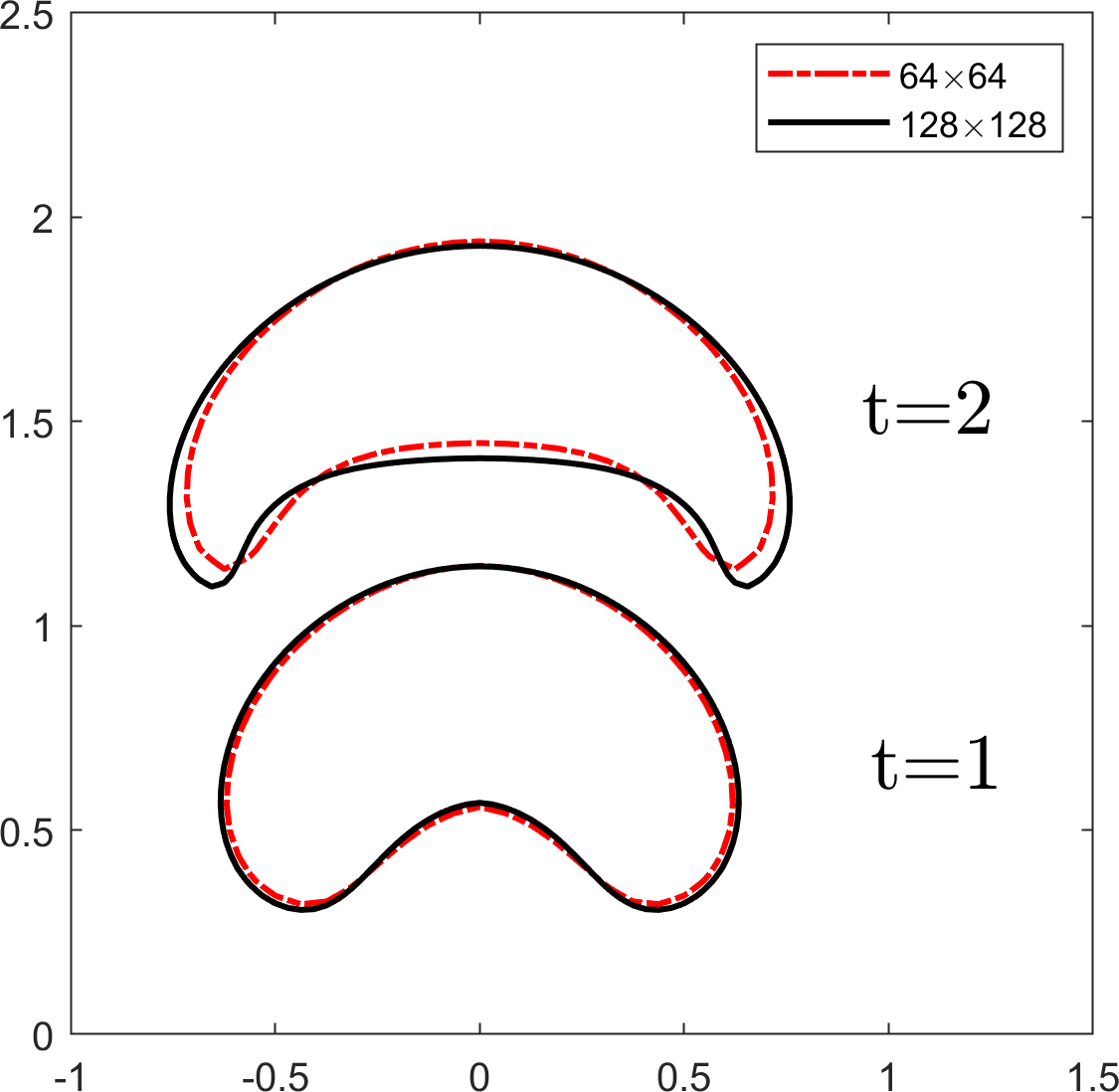}}\subfigure[][]{	\includegraphics[width=0.45\textwidth]{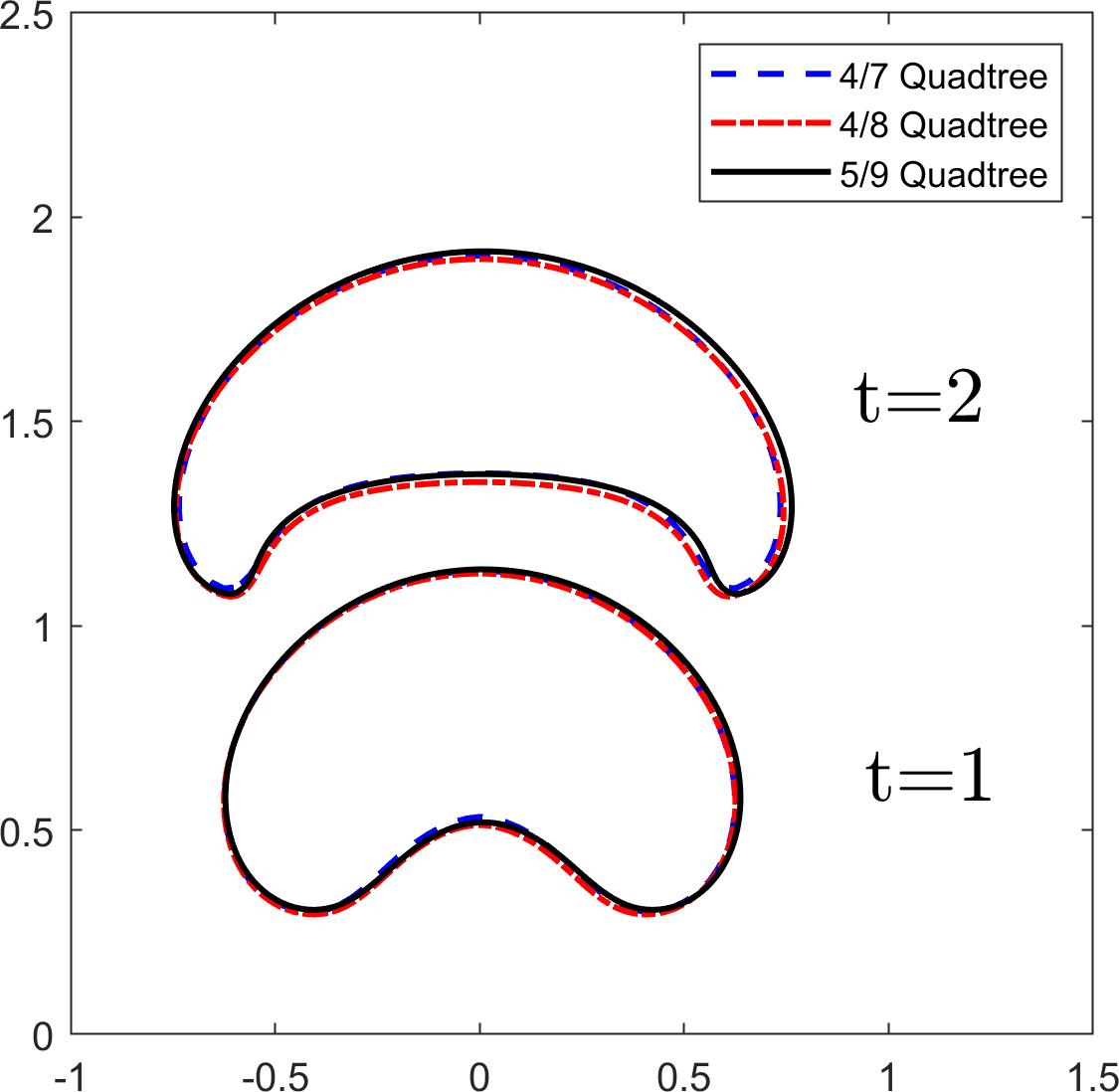}}}
		\caption{Visualization of the bubble at (a) uniform cartesian grid and (b) adaptive quadtree grid for example \ref{subsubsec:bubble}.}\label{fig:bubble_vis}
	\end{figure}

\begin{figure}
		\centering{}
		\mbox{	\subfigure[Relative area loss]{\includegraphics[width=0.3\textwidth]{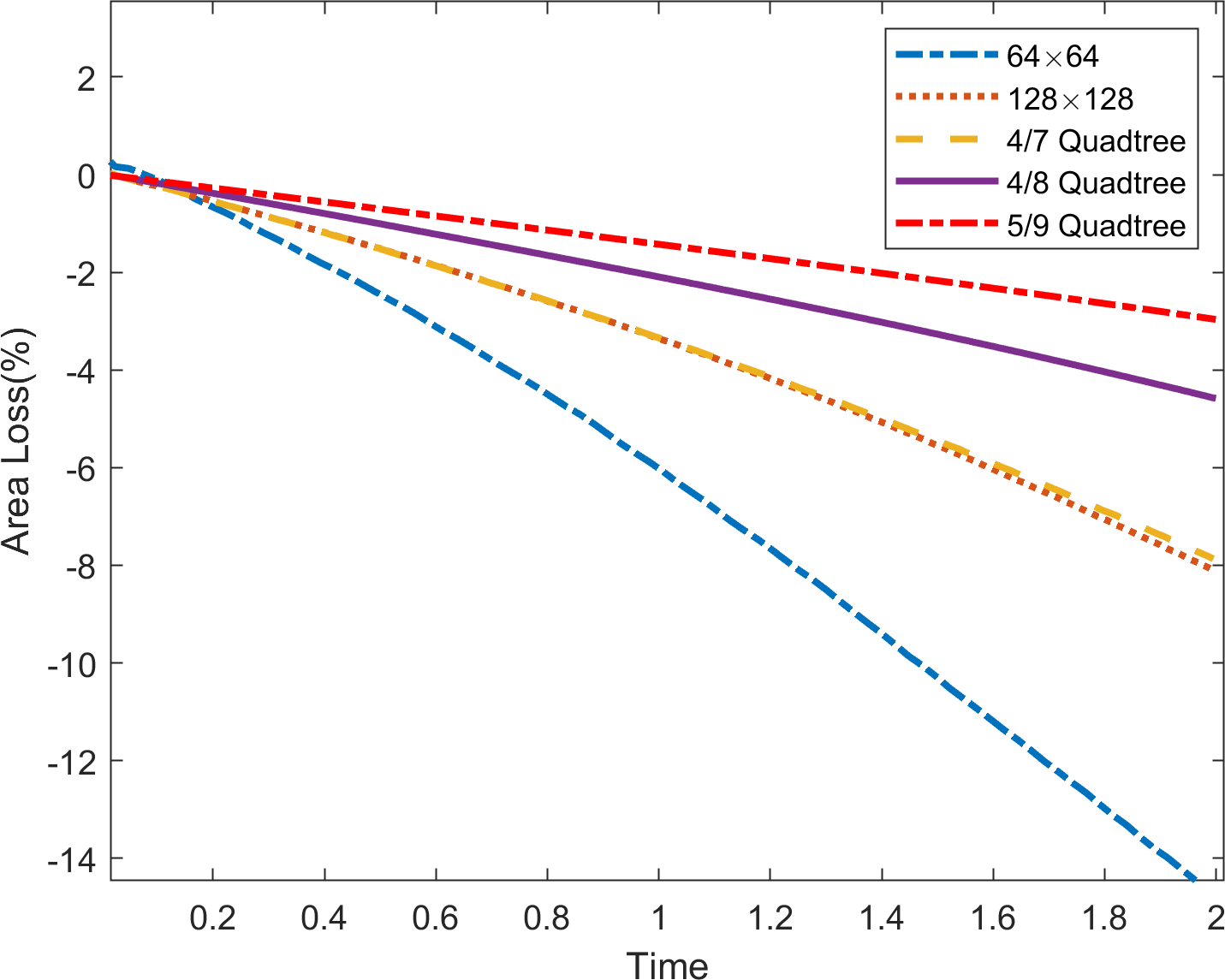}}	\subfigure[][Rising velocity]{\includegraphics[width=0.3\textwidth]{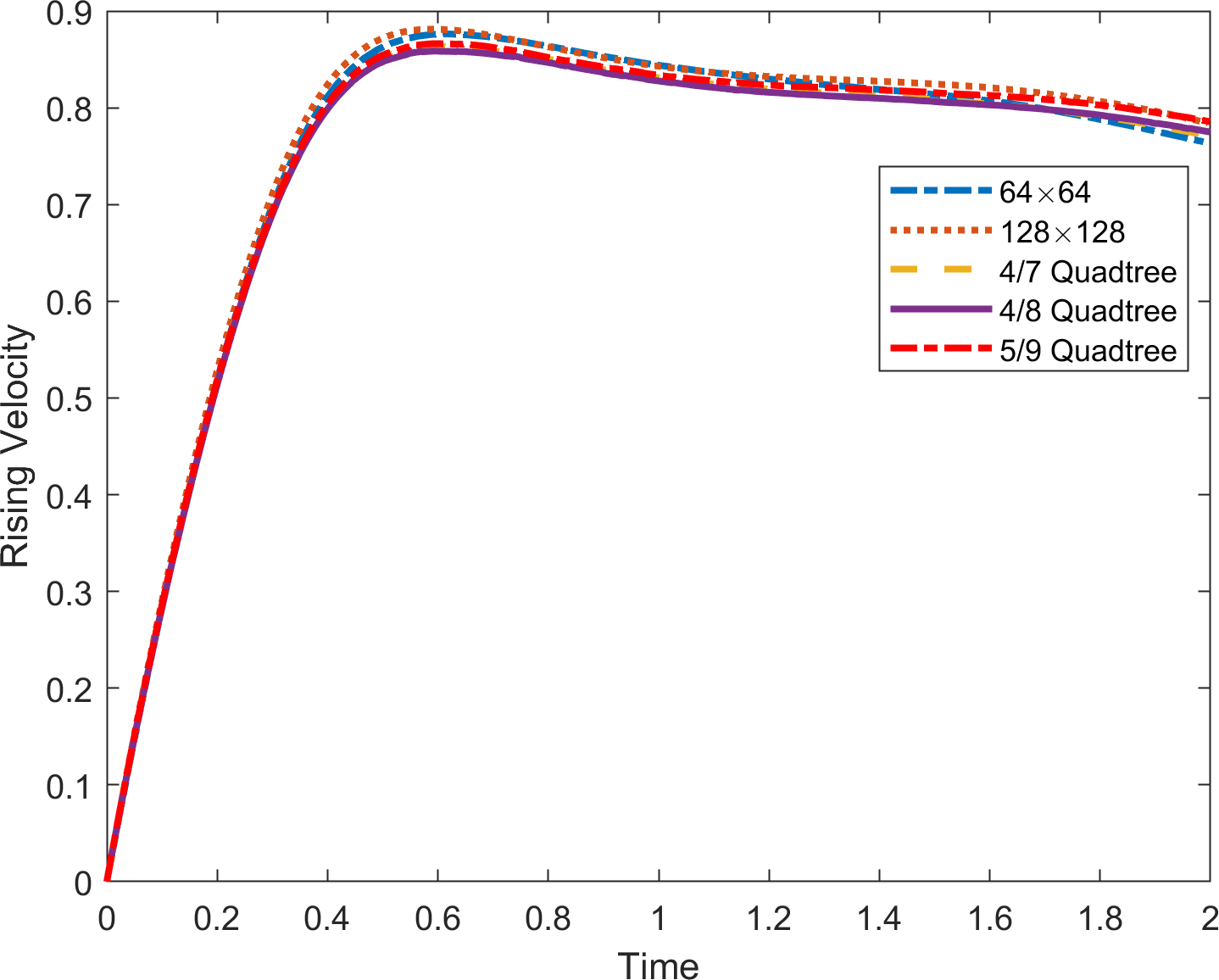}}	\subfigure[][Center of mass]{\includegraphics[width=0.3\textwidth]{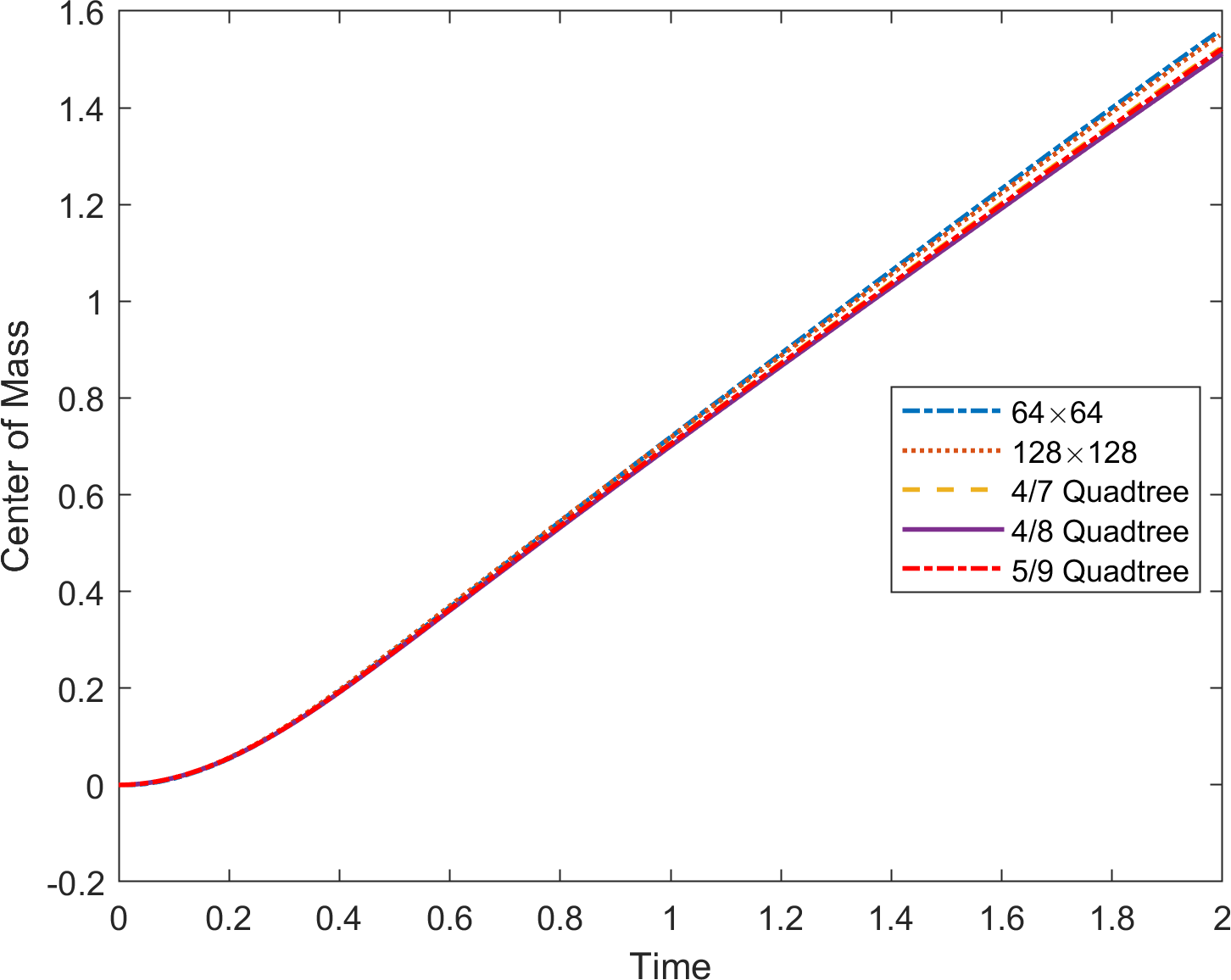} } }\\
		\mbox{\subfigure[][Rising velocity (zoomed)]{\includegraphics[width=0.3\textwidth]{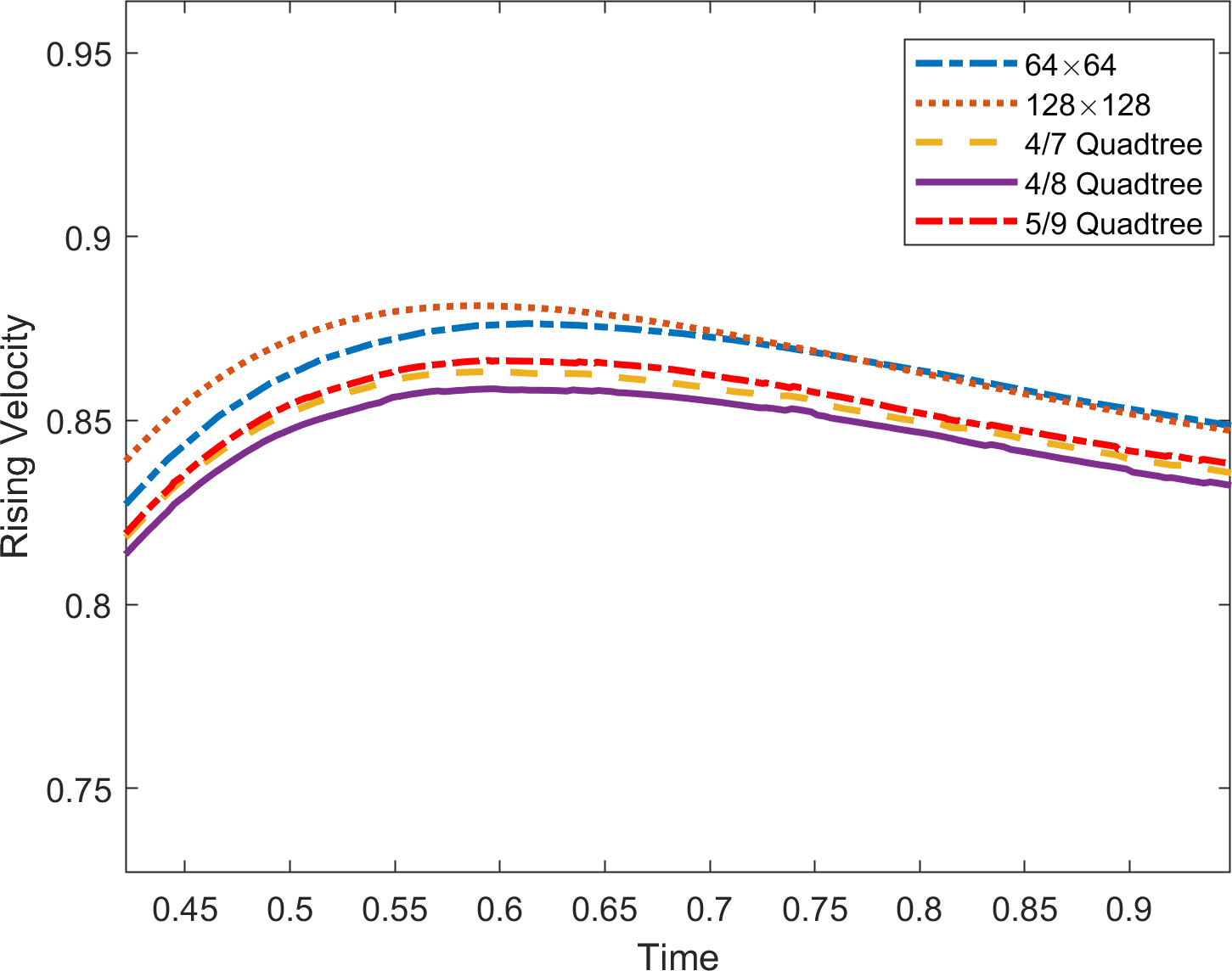}	}\subfigure[][Center of mass (Zoomed) ]{\includegraphics[width=0.3\textwidth]{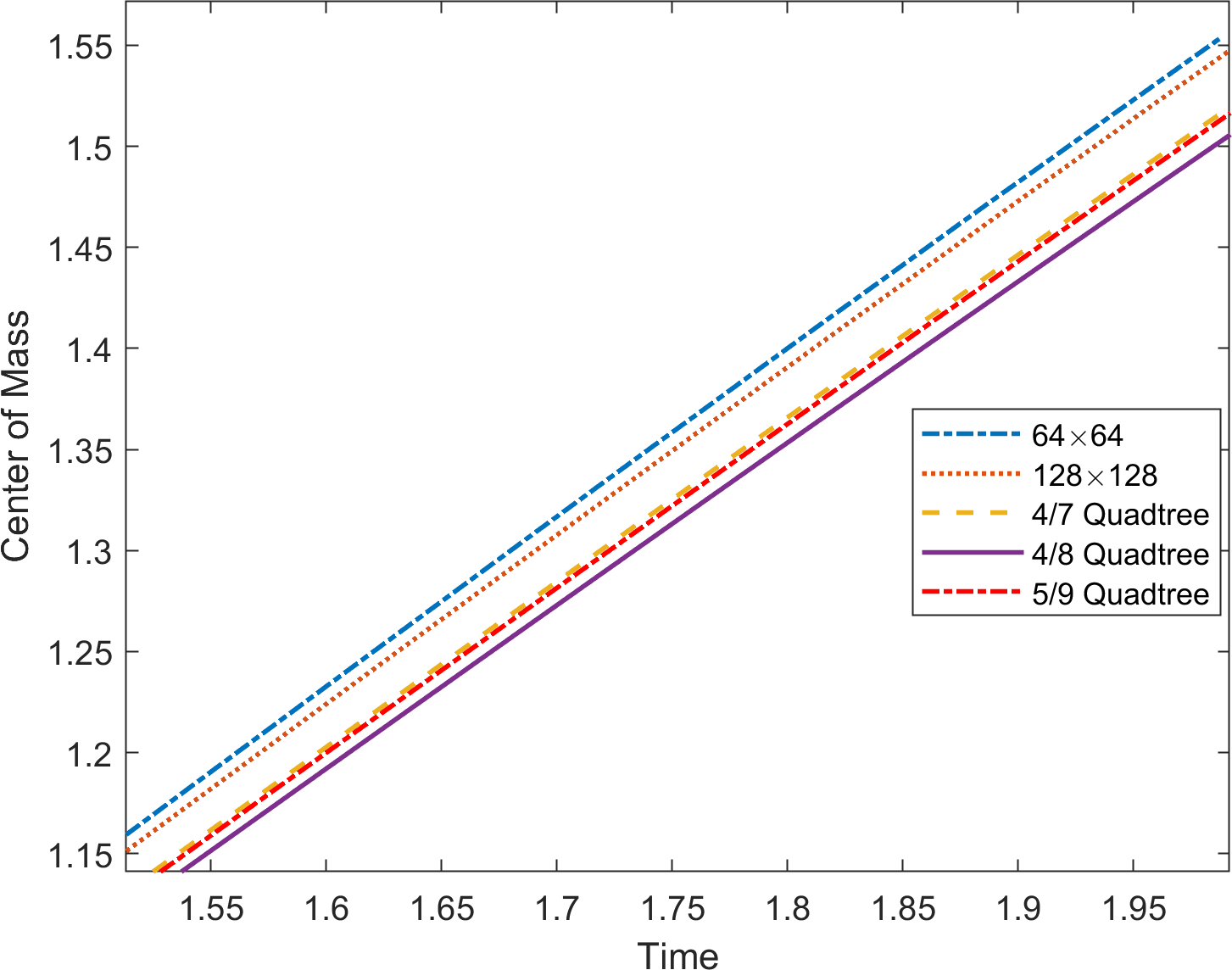}}}
		\caption{Relative area loss, rising velocity, and center of mass for rising bubble in example \ref{subsubsec:bubble}.}\label{fig:bubble_quant}
	\end{figure}

	In this section, we consider the evolution of a bubble rising in a rectangular domain $\left[-2,2\right] \times\left[-1,3\right]$. A bubble of radius $0.5$ is initially located at the origin. 
	We consider three-dimensionless quantities:
	the Reynolds number $Re=7.77$, the E\"{o}tv\"{o}s number $Eo=243$, and the Morton number $Mo=266$. Then the corresponding density, viscosity, surface tension, and gravity coefficient are determined:
	\[\rho^{+}=1, \quad \rho^{-}=0.001, \quad \mu^{+}=\frac{\rho^{+}}{R e}, \quad \mu^{-}=0.01 \mu^{+}, \quad \beta=\frac{\mu^{+^{2}}}{\rho^{+}} \sqrt{\frac{E o}{M o}}, \quad \mathrm{~g}=\left(0,-\frac{\rho^{+} \beta^{3}}{\mu^{+}}\right) .\]
	Simulations were performed up to $t=2$, and $\Delta t$ was computed under the following Courant--Friedrichs--Lewy condition:
	\[
	\Delta t = \min \left(\frac{\Delta x}{\|\mathbf{U}\|_{\infty}},0.7\sqrt{ \frac{\rho^{+}+\rho^{-}}{4 \pi \sigma} \Delta x^{3}}\right).
	\]
	In Figure \ref{fig:bubble_vis}, the terminal shapes of the bubble at time $t=2$ are shown on uniform grids of $64\times 64,128\times 128$ and quadtree grids of level $4/8,4/9$, and $5/9$. In addition, the relative area loss, rising velocity, and center of mass are plotted in Figure \ref{fig:bubble_quant}. We observe that the results on the quadtree grid are comparable to those on the uniform grid. 
	
	\setcounter{subfigure}{0}

	\section{Conclusion}
	
	We presented a second-order-accurate monolithic method for solving the two-dimensional incompressible Navier--Stokes equations on irregular domains and quadtree grids. In the proposed grid layout, velocities are located at cell vertices, and pressure variables are located at cell centers. Simple and effective second-order finite differences on the vertices of a quadtree grid \cite{min2006second} were used to discretize the advection and diffusion terms of the velocity. The pressure gradient was discretized using finite differences, and ghost points were used for cells outside the domain. The divergence of the velocity was discretized using the divergence theorem, whereas a cut-cell method was used to treat the boundary conditions. Numerical experiments included analytical and practical examples. For the analytical test, we presented the order of convergence of the velocity and pressure in both the $L^\infty$ and $L^2$ norms. The results indicate that the suggested method is second-order accurate for velocity regardless of the Reynolds number and domain shape. Second-order convergence of the pressure was observed only for relatively high Reynolds numbers, and first-order convergence was obtained for low Reynolds numbers. Classical lid-driven cavity flow and recent numerical experiment of incompressible flow on irregular domains were conducted and compared to the reference results. To demonstrate a more demanding application, numerical simulation of a rising bubble was investigated.
	
	Although an extension of the proposed method to an octree grid on a cubical domain was presented, the full three-dimensional extension of the method to irregular domains was not studied. The main difficulty in three dimensions is the discretization of the divergence-free equation using the cut-cell approach. As reported by Cheny et al. \cite{cheny2010ls}, there are many different cases, and developing consistent discretization on irregular domains may require surface integration on general polygons and curved elements. As a future work, we will extend the proposed method on an octree grid to three-dimensional irregular domains. Furthermore, the extension of the proposed method to a multiphase flow problem and fluid--structure interaction will also be an interesting future work.

	\bibliographystyle{spmpsci} 
	\bibliography{mybib}

\end{document}